\documentclass{ws-m3as}

\usepackage{graphicx}
\usepackage{amsmath}
\usepackage{amsfonts}
\usepackage{xcolor}
\usepackage{hyperref} 

\newtheorem{theo}{Theorem}[section]

\newtheorem{prop}[theo]{Proposition}
\newtheorem{cor}[theo]{Corollary}

\newcommand{\R}{{\Re}}
\newcommand{\lk}{\left\{}
\newcommand{\rk}{\right\}}

\newcommand{\g}{\gamma}

\newcommand{\Lu}{\mathcal{L}}
\newcommand{\OuO}{\Omega \cup \Omega_{\mathcal I}}

\newcommand{\mcS}{\mathcal{S}}
\newcommand{\mcK}{\mathcal{K}}

\newcommand{\mcE}{\mathcal{E}}

\newcommand{\mcI}{\mathcal{I}}

\newcommand{\mcL}{\mathcal{L}}

\newcommand{\mcT}{\mathcal{T}}
\newcommand{\mcX}{\mathcal{X}}
\newcommand{\mcO}{\mathcal{O}}

\newcommand{\eps}{\delta}

\newcommand{\del}{\delta}

\def \be{\begin{equation}}

\def \ba{\begin{array}\ds}
\def \ea{\end{array}}

\def\KKKK{K}

\def\KKddm1{{\KKKK_{\eps}^{-1}}}

\def \vb{\mathbf{v}}
\def \xb{\mathbf{x}}

\def \pb{\mathbf{p}}

\def \yb{\mathbf{y}}

\def \wuh{\widehat u_h}

\def\QQQ{{K}}

\def\acancel{{\alpha}}

\def\hangfour{\vskip0pt\hangindent=35pt\hangafter=1}
\def\hangseven{\vskip0pt\hangindent=29pt\hangafter=1}

\def\hangthree{\vskip0pt\hangindent=48pt\hangafter=1\qquad}

\def\hangfive{\vskip0pt\hangindent=15pt\hangafter=1\noindent}

\def\beq{\begin{equation}\label}
\def\eeq{\end{equation}}
\def\bal{\begin{aligned}} 
\def\eal{\end{aligned}}

\DeclareMathOperator{\argmin}{\mbox{argmin}}

\begin{document}

\markboth{M. D'Elia, M. Gunzburger, and C. Vollmann}{
{Finite Element Methods for Nonlocal Problems}
}

\catchline{}{}{}{}{}

\title{A COOKBOOK FOR FINITE ELEMENT METHODS FOR NONLOCAL PROBLEMS, INCLUDING QUADRATURE RULES AND APPROXIMATE EUCLIDEAN BALLS}

\author{MARTA D'ELIA}

\address{
{Computational Science and Analysis\\ Sandia National Laboratories\\ Livermore, CA, 94550, USA\\{mdelia@sandia.gov}}}

\author{MAX GUNZBURGER}

\address{Department of Scientific Computing\\ Florida State University\\ Tallahassee, FL 32306, USA\\{mgunzburger@fsu.edu}}

\author{CHRISTIAN VOLLMANN}

\address{{Department of Mathematics\\ Universit{\"a}t Trier\\ 54296 Trier, Germany\\{vollmann@uni-trier.de}}}

\maketitle

\begin{history}
\received{(Day Month Year)}
\revised{(Day Month Year)}
\accepted{(Day Month Year)}
\comby{(xxxxxxxxxx)}
\end{history}

\begin{abstract}
{ 
The implementation of finite element methods (FEMs) for nonlocal models with a finite range of interaction poses challenges not faced in the partial differential equations (PDEs) setting. For example, one has to deal with weak forms involving double integrals which lead to discrete systems having higher assembly and solving costs due to possibly much lower sparsity compared to that of FEMs for PDEs. In addition, one may encounter non-smooth integrands. In many nonlocal models, nonlocal interactions are limited to bounded neighborhoods that are ubiquitously chosen to be Euclidean balls, resulting in the challenge of dealing with intersections of such balls with the finite elements. We focus on developing recipes for the efficient  assembly of FEM stiffness matrices and on the choice of quadrature rules for the double integrals that contribute to the assembly efficiency and also posses sufficient accuracy. A major feature of our recipes is the use of approximate balls, e.g., several polygonal approximations of Euclidean balls, that, among other advantages, mitigate the challenge of dealing with ball-element intersections. We provide numerical illustrations of the relative accuracy and efficiency of the several approaches we develop.} 
\end{abstract}



\section{Introduction}\label{sec:intro}

Nonlocal models provide improved simulation fidelity in the presence of long-range forces and anomalous behaviors. Because of their integral form, they can capture long-range effects and relax the regularity requirements of classical (differential) models. For this reason, their applicability ranges over fracture mechanics (Refs. \refcite{Ha2011,Littlewood2010,Silling2000}), image processing (Refs. \refcite{Buades2010,Gilboa2007,Gilboa2008,Lou2010}), stochastic processes (Refs. \refcite{Burch2014,DElia-conv-diff,Meerschaert2012}), anomalous subsurface transport (Refs. \refcite{Benson2000,Schumer2003,Schumer2001}), multiscale and multiphysics systems (Refs. \refcite{Alali2012,Askari2008}), phase transitions (Refs. \refcite{Bates1999,Delgoshaie2015,Fife2003}), {and machine learning (Ref. \refcite{Wei2020}).}

The central difference between nonlocal models and partial differential equation (PDE) models is that for the former, interactions can occur at distance, whereas for the latter, they can only occur through contact. As a consequence, in nonlocal settings, a point in space at a time instant interacts with a neighborhood of points and with previous times instants, i.e., far away in space and far back in time.

Nonlocality raises many modeling and computational challenges. The former include the prescription of nonlocal analogues of boundary conditions (see Refs. \refcite{Cortazar2008,DEliaNeumann2019,Lischke2020}), the choice of kernel functions that characterize nonlocal operators (see, e.g., Refs. \refcite{DElia2016ParamControl,optcontrol,Gulian2019,Pang2020nPINNs,Pang2019fPINNs}), and the modeling of nonlocal interfaces (see Refs. \refcite{Alali2015,Capodaglio2019}). The computational challenges include the design of efficient quadrature rules for possibly singular kernel functions, the construction of nonlocal discrete systems, and the design of efficient nonlocal solvers. In fact, the numerical solution of nonlocal models is, relative to PDE models, intrinsically extremely expensive with respect to both assembling and solving discrete systems (see Ref. \refcite{DElia-ACTA-2020}). 
 
Meshfree, in particular particle-type methods, provide a popular means for discretizing nonlocal equations; see, e.g., Refs. \refcite{parks2012peridigm,parks2010lammps}. Here, however, we are interested in variational methods, and in particular finite element methods, because of the ease they provide for dealing with complicated domains, {for obtaining approximate solutions that have higher-order convergence rates,} and for defining adaptive meshing methods that can resolve solution misbehaviors such as jump discontinuities and steep gradients, the latter also arising in the PDE setting. In addition, casting the nonlocal problem into a variational framework used to define finite element methods allows for a rigorous mathematical treatment of operator and solution properties, well posedness, and stability and convergence of approximate solutions. 

In this paper, we focus on some of the computational challenges one must face in the design of efficient finite element methods in the nonlocal setting. We summarize the main contributions of this paper.

\smallskip\noindent
{\bf 1.} This is the first work where nonlocal finite element formulations and associated implementation tasks are thoroughly and rigorously addressed and illustrated. In fact, not only do we describe the assembly procedure in detail, but we also provide guidance about the choice of quadrature rules for the outer and inner integrals\footnote{As opposed to finite element methods for PDEs for which the weak form involves integration over the domain, finite element methods for nonlocal models require a double integration over the domain due to the integral form of nonlocal operators.} in relation to other errors incurred such as that due to finite element approximation.

\smallskip\noindent
{\bf 2.} We introduce approximate nonlocal neighborhoods that facilitate the assembly procedure and mitigate the computational effort. For each of them, we describe the geometric approximation and discuss the errors they incur. Again, we provide guidance about the choice of quadrature rules to use for each specific neighborhood approximation so that the overall accuracy is not compromised.

\smallskip\noindent
{\bf 3.} Among such neighborhood approximations, we provide numerical evidence, in two dimensions, that particularly inexpensive and easy-to-implement approximations preserve optimal accuracy, while significantly reducing  computational costs, making those approaches also the best candidates for three-dimensional simulations. Those techniques could potentially make variational methods as efficient as meshfree methods and, hence, become preferable alternatives.

\smallskip
In Sec. \ref{sec:notation} we introduce the strong form of the nonlocal problem and in so doing we define nonlocal operators, kernels, and domains. In Sec. \ref{sec:weak} we discuss the most straightforward variational formulation and review relevant elements of the nonlocal vector calculus developed in Ref.  \refcite{dglz2}. In Sec. \ref{fem} we describe finite element discretizations by providing their formulation, recipes for the assembly of discrete systems, accuracy results, and several useful tips and remarks. In Sec. \ref{approxballs} we introduce several geometric approximations of the nonlocal neighborhood that is in ubiquitous use in nonlocal modeling, namely Euclidean balls. By rigorously estimating the difference between approximated variational forms defined by the approximate balls to that for the exact ball, we show how such approximations (in combination with quadrature rules) affect the discretization error. In Sections \ref{apinner} and \ref{apouter} we describe quadrature rules for the double integral that appears in the weak formulation, highlight the desired properties one would want them to have, provide guidance about the choice of quadrature points and weights, and discuss how those choices affect accuracy. In Sec. \ref{efficient} we show how the quadrature rules lead to  fully-discrete finite element formulations for which we discuss efficient assembly procedures. In Sec. \ref{numerics} we illustrate the theoretical findings with several two-dimensional numerical tests and then, in Sec. \ref{conclusion}, provide some concluding remarks.

\subsection{The problem setting}\label{sec:notation}

Consider the nonlocal Dirichlet problem
\beq{prob}
\left\{\bal
-\Lu u(\xb) = f(\xb) \qquad & \mbox{for $\xb\in \Omega$}\\
     u(\xb) = g(\xb) \qquad & \mbox{for $\xb\in \Omega_{\mathcal I}$}, 
\eal\right.
\eeq
where $\Omega \subset \R^d$ denotes an open bounded domain,
\beq{oper}
\Lu u(\xb) = 2\int_{\OuO}\big(u(\yb)-u(\xb)\big)\g(\xb,\yb)d\yb 
\eeq
denotes a {\em nonlocal operator}, and $\g(\xb,\yb): \R^d\times\R^d \to \R$ denotes a nonnegative and symmetric function, i.e., $\g(\xb,\yb)=\g(\yb,\xb)$ for all $\xb$ and $\yb$, which we refer to as the \emph{kernel}.\footnote{For a discussion on nonpositive kernels and nonsymmetric kernels, see Ref.  \refcite{Mengesha-sign-changing} and Ref.  \refcite{DElia-conv-diff}, respectively.} In \eqref{prob} and \eqref{oper}, $\Omega_{\mathcal I}$ denotes the {\em interaction domain} corresponding to $\Omega$, defined to be the set of points in the complement domain $\R^d\setminus \Omega$ that interact with points in $\Omega$. More precisely, we define $\Omega_{\mathcal I}$ as
\beq{inter_intro}
  \Omega_{\mathcal I} = \big\{ \yb\in \R^d\setminus \Omega \,\,:\,\,
  \exists\, \xb\in{\overline\Omega} \,\,\,\mbox{such that}\,\,\, \g(\xb,\yb)\ne0\big\}
  \subset \R^d\setminus \Omega.
\eeq
Note that $\Omega_{\mathcal I}$ so defined is a closed domain, and, in particular, $\Omega_{\mathcal I}\cap\partial\Omega=\partial\Omega$, where $\partial\Omega$ denotes the boundary of $\Omega$. With $f(\xb): \Omega\to\R$ and $g(\xb): \Omega_{\mathcal I}\to\R$ denoting given functions, the problem \eqref{prob} determines $u(\xb): \OuO\to\R$. 

We refer to the second equation in \eqref{prob} as a {\em Dirichlet volume constraint}, with ``Dirichlet'' because the solution itself is specified on $\Omega_{\mathcal I}$ and ``volume constraint'' referring to that equation holding on a set having finite volume in $\R^d$, in contrast to the local PDE setting in which a Dirichlet constraint is applied on a $(d-1)$-dimensional surface. Hence, it is also natural to refer to problem \eqref{prob} as a nonlocal {{\em volume-constrained Dirichlet problem.}}\footnote{For the sake of economy of the exposition, we do not consider nonlocal {\em Neumann} problems. Such problems are considered in, e.g., Ref. \refcite{dglz1}.}

\begin{figure}[ht] 
\centering
\includegraphics[width=1in]{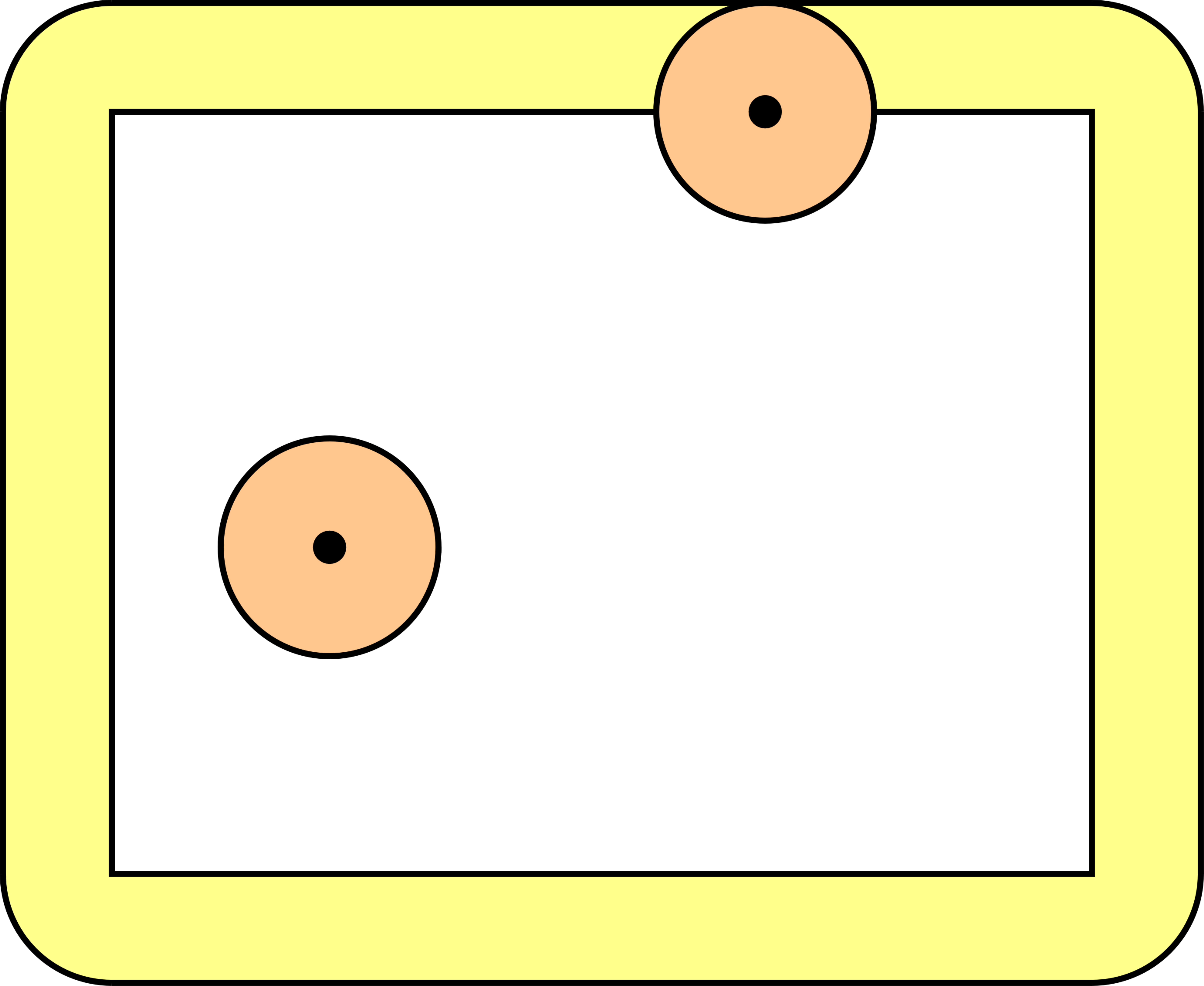}  
\caption{In white, a rectangular domain $\Omega$; in yellow, the corresponding interaction domain $\Omega_{\mathcal I}$ of thickness $\del$; in orange, two balls of radius $\del$ centered at the two points in $\Omega\cup\Omega_{\mathcal I}$ depicted by black dots, one of which is located on the boundary $\partial\Omega$ between $\Omega$ and $\Omega_{\mathcal I}$.}
\label{nonlocalsetup}
\end{figure}

The case of $\Omega=\R^d$ (so that $\Omega_{\mathcal I}=\emptyset$) could also be included as could the case $\Omega_{\mathcal I}=\R^d\setminus \Omega$ that corresponds to interactions occurring over an infinite distance. However, motivated by the fact that, in real-world applications, interactions do not occur over infinite distances, we only consider {\em kernels having bounded support} for which two points in  $\xb, \yb \in \R^d$ interact which each other, i.e., $\g(\xb,\yb)\neq 0$, only if $\yb$ is within a bounded neighborhood of $\xb$. For that neighborhood, we focus on the specific choice of closed Euclidean balls $B_\del(\xb)$ centered at $\xb$ having radius $\del$ that is in ubiquitous use in the literature;\footnote{Although we focus on Euclidean balls, the discussion and results in this paper can be extended to cover balls of other types, e.g., $\ell^\infty$-norm balls, and to even more general interaction sets.} $\del$ is often referred to as the {\em horizon} or {\em interaction radius}. Thus, we have that
\beq{inter}
  \g(\xb,\yb)= \psi(\xb,\yb)\mcX_{B_\del(\xb)}(\yb)
\eeq
for some symmetric and positive function $\psi(\xb,\yb)$ that we refer to as the {\em kernel function}, where $\mcX_{\{\cdot\}(\xb)}(\yb)$ denotes the indicator function. Note that $\g(\xb,\yb)$ given by \eqref{inter} is a symmetric function because $\mcX_{B_\del(\xb)}(\yb)$ is itself symmetric; {in fact,} if $\yb\in B_\del(\xb)$ then necessarily $\xb\in B_\del(\yb)$. Fig. \ref{nonlocalsetup} illustrates {a domain $\Omega$, its interaction domain $\Omega_{\mathcal I}$ that results from \eqref{inter}, and two balls $B_\del(\xb)$, one centered at $\xb\in\Omega$ and the other at $\xb\in\partial\Omega$.}

\section{\textbf{Weak formulation}}\label{sec:weak}

A weak formulation of the problem \eqref{prob} can be derived in the usual manner. Proceeding formally, we multiply the first equation in \eqref{prob} by a test function $v(\xb)$ to obtain\footnote{Throughout, when we encounter double integrals such as $\int\big(\int (\cdots) d\yb\big) d\xb$, we refer to $\int (\cdots) d\yb$ as the {\em inner integral} and to $\int \big(\cdots\big) d\xb$ as the {\em outer integral}.}
\beq{step0}
\bal
 0&=\int_\Omega v(\xb)\big(-\Lu u(\xb) - f(\xb)\big) d\xb 
 \\&=  
 2\int_\Omega v(\xb) \int_{\OuO}\big(u(\xb)-u(\yb)\big)\g(\xb,\yb)d\yb d\xb - \int_\Omega v(\xb)f(\xb)d\xb.
\eal
\eeq
Because the second equation in \eqref{prob} is a Dirichlet-type constraint imposed on $\Omega_{\mathcal I}$, i.e., it is a constraint on the solution $u(\xb)$ itself, we require that the test function satisfies $v(\xb)=0$ for $\xb\in\Omega_{\mathcal I}$. Then, applying Green's first identity of the nonlocal vector calculus given in Ref. \refcite{dglz2} to the first term in \eqref{step0}, we have, with $v(\xb)=0$ for $\xb\in\Omega_{\mathcal I}$,
\beq{step00}
\bal
&2\int_\Omega v(\xb) \int_{\OuO}\big(u(\xb)-u(\yb)\big)\g(\xb,\yb)d\yb d\xb
\\&\qquad= \int_{\OuO} \int_{\OuO} \big(u(\yb)-u(\xb)\big)\big(v(\yb)-v(\xb)\big)\g(\xb,\yb)d\yb d\xb.
\eal
\eeq
Combining \eqref{step0} and \eqref{step00}, we have
\beq{weak0}
D(u,v)=G(v),
\eeq
where
\beq{bform0}
D(u,v)=\int_{\OuO} \int_{\OuO}
    \big(u(\yb)-u(\xb)\big)\big(v(\yb)-v(\xb)\big)\g(\xb,\yb)d\yb d\xb
\eeq
and 
\beq{lfunc0}
   G(v) = \int_\Omega v(\xb)f(\xb)d\xb.
\eeq
Applying the volume constraint in \eqref{prob} to set $u(\xb)=g(\xb)$ on $\Omega_{\mathcal I}$ and again setting $v(\xb)=0$ on $\Omega_{\mathcal I}$, we obtain from \eqref{weak0} that
\beq{step1}
\bal
0 & = \int_{\Omega}\int_{\Omega} \big(u(\yb)-u(\xb)\big)\big(v(\yb)
    - v(\xb)\big)\g(\xb,\yb)d\yb \,d\xb 
\\&\qquad\qquad
- \int_{\Omega}
    \int_{\Omega_{\mathcal I}} 
      \big(g(\yb)-u(\xb)\big)v(\xb)\g(\xb,\yb)d\yb \,d\xb \\
  &   \qquad\qquad +\int_{\Omega_{\mathcal I}}\int_{\Omega} \big(u(\yb)-g(\xb)\big)v(\yb)\g(\xb,\yb)d\yb \,d\xb 
    - \int_\Omega f(\xb)v(\xb) d\xb\\
  & = \int_{\Omega}\int_{\Omega} \big(u(\yb)-u(\xb)\big)
      \big(v(\yb)-v(\xb)\big)\g(\xb,\yb)d\yb \,d\xb\\
  &   \qquad\qquad + \int_{\Omega}\int_{\Omega_{\mathcal I}} u(\xb)v(\xb)\g(\xb,\yb)d\yb \,d\xb 
    + \int_{\Omega_{\mathcal I}}\int_{\Omega} u(\yb)v(\yb)\g(\xb,\yb)d\yb \,d\xb\\
  &   \qquad\qquad - \int_{\Omega}\int_{\Omega_{\mathcal I}} g(\yb)v(\xb)\g(\xb,\yb)d\yb \,d\xb 
    - \int_{\Omega_{\mathcal I}}\int_{\Omega} g(\xb)v(\yb)\g(\xb,\yb)d\yb \,d\xb
 \\&\qquad\qquad    -\int_\Omega f(\xb)v(\xb) d\xb.
\eal
\eeq
Changing the order of the integration, renaming the dummy variables of integration, and using the symmetry of the kernel $\g(\xb,\yb)$, we have that
\beq{step2}
\bal
   & \int_{\Omega_{\mathcal I}}\int_{\Omega} u(\yb)v(\yb)\g(\xb,\yb)d\yb \,d\xb  
 = \int_{\Omega}u(\yb)v(\yb) \int_{\Omega_{\mathcal I}} \g(\xb,\yb)d\xb \,d\yb \\
& \qquad\qquad
= \int_{\Omega}u(\xb)v(\xb) \int_{\Omega_{\mathcal I}} \g(\yb,\xb)d\yb \,d\xb
  = \int_{\Omega}u(\xb)v(\xb) \int_{\Omega_{\mathcal I}} \g(\xb,\yb)d\yb \,d\xb
\eal
\eeq
and similarly
\beq{step3}
\bal
 &   \int_{\Omega_{\mathcal I}}\int_{\Omega} g(\xb)v(\yb)\g(\xb,\yb)d\yb \,d\xb
 = \int_{\Omega}v(\yb)\int_{\Omega_{\mathcal I}} g(\xb)\g(\xb,\yb)d\xb \,d\yb \\
&\qquad\qquad = \int_{\Omega}v(\xb)\int_{\Omega_{\mathcal I}} g(\yb)\g(\yb,\xb)d\yb \,d\xb
  = \int_{\Omega}v(\xb)\int_{\Omega_{\mathcal I}} g(\yb)\g(\xb,\yb)d\yb \,d\xb.
\eal
\eeq
Combining \eqref{bform0}--\eqref{step3}, we have
$$
\bal
0= &  \int_{\Omega}\int_{\Omega } \big(u(\yb)-u(\xb)\big)\big(v(\yb)-v(\xb)\big)\g(\xb,\yb)d\yb \,d\xb
\\&\qquad    + 2\int_{\Omega}u(\xb)v(\xb)\int_{\Omega_{\mathcal I} } \g(\xb,\yb)d\yb \,d\xb 
\\&\qquad   - 2\int_{\Omega}v(\xb)\int_{\Omega_{\mathcal I} } g(\yb)\g(\xb,\yb)d\yb \,d\xb 
    - \int_\Omega f(\xb)v(\xb) d\xb.
\eal
$$
Thus, we have that
\beq{weakf2}
A(u,v) = F(v)
\eeq
with the symmetric bilinear form 
\beq{bform2}
\bal
A(u,v) =& \int_{\Omega}\int_{\Omega } \big(u(\yb)-u(\xb)\big)\big(v(\yb)-v(\xb)\big)\g(\xb,\yb)d\yb d\xb
\\&\qquad\qquad
 + 2\int_{\Omega}u(\xb)v(\xb)\bigg(\int_{\Omega_{\mathcal I} } \g(\xb,\yb)d\yb\bigg) \,d\xb 
\eal
\eeq
and the linear functional
\beq{lfunc2}
F(v) = \int_{\Omega}v(\xb) \bigg(f(\xb) + 2\int_{\Omega_{\mathcal I} } g(\yb)\g(\xb,\yb)d\yb \bigg) d\xb.
\eeq
It is useful to note that 
\beq{equival}
\mbox{\em \eqref{weak0} along with $u(\xb)=g(\xb)$ and $v(\xb)=0$ on $\Omega_{\mathcal I}$ are equivalent to \eqref{weakf2}.}
\eeq
Throughout, we take advantage of this equivalence by using one or the other of the pairs $\{D(u,v),G(v)\}$ and $\{A(u,v),F(v)\}$ as is most convenient for describing the specific task at hand.
 
We are now in position to define a weak formulation of the problem \eqref{prob}. To this end, for functions $v(\xb)$ defined for $\xb\in\OuO$, we define the norm $|||v|||= \sqrt{A(v,v)} + \|v\|_{L^2(\Omega_{\mathcal I})}$ and the function spaces, often referred to as the (nonlocal) ``energy'' spaces,
\beq{espace}
\left\{\bal
V(\OuO) &= \{ v \in L^2(\OuO)\,\, \colon \,\, |||v|||  <\infty\} \\
V_c(\OuO) &= \{ v \in V(\OuO)\,\, \colon \,\, v=0 \mbox{ on $\Omega_{\mathcal I}$} \} .
\eal
\right.
\eeq
Because $\gamma(\xb,\yb)>0$ for $\yb \in B_\del(\xb)$ by assumption \eqref{inter}, the bilinear form $A(\cdot,\cdot)$ is positive, i.e., $A(v,v)>0$ for all $v\in V_c(\OuO)$ such that $v\ne0$. Thus, $A(u,v)$ is an inner product on $V_c(\OuO)\times V_c(\OuO)$ and $\sqrt{A(v,v)}$ is a norm on $V_c(\OuO)$. We also introduce the trace space $V_t(\OuO)=\{v|_{\Omega_{\mathcal I}}~\colon \,\, v\in V(\OuO)\}$ and denote by $V_c'(\Omega)$ the dual space whose elements are bounded linear functionals on $V_c(\OuO)$.

We define the weak formulation of \eqref{prob} as follows. {\em Given $f(\xb)\in V_c'(\Omega)$ and $g(\xb)\in V_t(\OuO)$, seek $u(\xb)\in V(\OuO)$ such that $u(\xb)=g(\xb)$ for $\xb\in\Omega_{\mathcal I}$ and $u(\xb)$ for $\xb\in\Omega$ is determined from the variational problem}
\beq{weak}
A(u,v) = F(v)\qquad \forall\, v\in V_c(\OuO) .
\eeq
The well posedness of the problem \eqref{weak} follows from {the Riesz representation theorem because $A(u,v)$ defines an inner product on $V_c(\OuO)$.}

For some specific kernels, it is known that the energy space $V(\OuO)$ is equivalent to standard function spaces. For example, for square integrable kernel functions or translationally invariant integrable kernel functions\footnote{Square integrable kernels satisfy $\int_{\OuO}\psi(\xb,\yb)^2d\yb<\infty$ for all $\xb\in\OuO$ and integrable kernel functions satisfy $\int_{\OuO}\psi(\xb,\yb)d\yb<\infty$ for all $\xb\in\OuO$. Translational invariant kernel functions are such that $\psi(\xb,\yb)=\psi(\yb-\xb)$.} $\psi(\xb,\yb)$, $V(\OuO)$ is equivalent to $L^2(\OuO)$. For non-integrable singular kernels,  $V(\OuO)$ is equivalent to function spaces of smoother functions defined on $\Omega$. For example, for kernels having the singular behavior of a fractional Laplacian kernel, $V(\OuO)$ is equivalent to $H^s(\Omega)\times L^2(\Omega_{\mathcal I})$ for an appropriate $s\in(0,1)$, where $H^s(\Omega)$ denotes the fractional Sobolev space of order $s$.

In what follows, to avoid further complications that arise in case of strongly singular kernels (e.g., non-integrable kernels) and which are not germane to the issues addressed here, we restrict our discussion to square integrable kernel functions or translationally invariant integrable kernels\footnote{Finite element discretizations, including proper choice of quadrature rules, for (non-truncated) fractional kernels have been investigated in Refs. \refcite{AinsworthGlusa2017} and \refcite{AinsworthGlusa2018}.}. However, we will briefly address the additional challenges posed by singular kernels in several remarks throughout the paper. 

\vskip5pt
\textbf{\em Sources of error.}
Of course, in practice, one implements a {\em fully-discrete} approximation of \eqref{weak}. In so doing, four types of errors can be possibly incurred:

\hangseven--~a finite element method is used to discretize \eqref{weak}; see Sec. \ref{fem};

\hangseven--~an approximate ball $B_{\del,h}(\xb)$ is used to approximate the ``exact'' ball $B_{\del}(\xb)$; see Sec. \ref{approxballs}; 

\hangseven--~a global or composite quadrature rule is used to approximate the inner integrals in \eqref{bform2} and \eqref{lfunc2}; see Sec. \ref{apinner};  

\hangseven--~a composite quadrature rule is used to approximate the outer integrals in \eqref{bform2} and \eqref{lfunc2}; see Sec. \ref{apouter}.

\noindent In principle, the four errors should be {\em commensurate}, i.e., none of the errors incurred should dominate the others and none should be dominated by any of the others. Otherwise, there would be wasteful computations involved. All of the errors listed above depend on the grid size $h$, so that, to be commensurate, all would have an error of $\mcO(h^\beta)$ as would the total error. Note that having one or more errors have a larger $\beta$ than the others cannot improve on the rate of convergence of the overall error, but could result in a smaller constant in error estimates and in smaller absolute errors in practice. 
\vskip5pt
\textbf{\em Analogies with (local) PDE problems.}
The problem \eqref{prob} with the operator \eqref{oper} is a nonlocal analogue of second-order elliptic PDE problems such as $-\Delta u =f$ in $\Omega$ and $u=g$ on the boundary $\partial\Omega$ {of $\Omega$}. The nonlocal weak problem \eqref{weak} is a nonlocal analogue to, e.g., the local weak formulation $\int_\Omega \nabla u(\xb) \cdot \nabla v(\xb) d\xb - \int_\Omega f(\xb)v(\xb)d\xb=0$ that is derived starting from $\int_\Omega v(\xb)\big( -\Delta u(\xb) -f(\xb)\big)d\xb=0$ using the classical (local) Green's first identity.

\vskip5pt
\textbf{\em Choice of weak formulation.}
In the local case, the form $\int_\Omega v(\xb)\big( -\Delta u(\xb) -f(\xb)\big)d\xb=0$ is well defined only for sufficiently smooth solutions and, in particular, it cannot be used as a weak formulation (i.e., it is not well defined) if, as is most often then case, the local energy space is chosen to be a subspace of the Sobolev space $H^1(\Omega)$. On the other hand, \eqref{step0} can be used as a nonlocal weak formulation in some settings. For example, if the kernel $\g(\xb,\yb)$ is integrable, then the nonlocal energy space is $V(\OuO) =L^2(\Omega)\times L^2(\Omega_{\mathcal I})$ for which \eqref{step0} is well defined; see Refs. \refcite{dglz1,dglz2}. In this case, \eqref{step0} with $u(\xb)=g(\xb)\in L^2(\Omega_{\mathcal I})$ for $\xb\in\Omega_{\mathcal I}$ is entirely equivalent to \eqref{weak}.

\vskip5pt
\textbf{\em Energy minimization characterization of the weak formulation.}
The weak formulations \eqref{weak0} and \eqref{weak} can also be derived from a minimization principle. Define the functional
$$
 {\mathcal J}(v;f) =  \frac12\int_{\OuO} \int_{\OuO} |v(\xb)-v(\yb)|^2\g(\xb,\yb) d\yb d\xb -
  \int_\Omega f(\xb) v(\xb) d\xb
$$
that is often referred to as a {\em nonlocal ``energy'' functional}. Then, given $f(\xb)\in L^2(\Omega)$ and $g(\xb)=L^2(\Omega_{\mathcal I})$, consider the minimization problem
$$
 u(\xb) =  \argmin\limits_{\{v\in V(\OuO):\,\,\, v=g\,\,\mbox{\em\footnotesize for}\,\, \xb\in\Omega_{\mathcal I}\}}  
 {\mathcal J}(v;f).
$$
It is easily seen that the minimizing function $u(\xb)$ is the solution of the weak formulation \eqref{weak}. We note that the approximate balls and quadrature rules discussed in this paper are also applicable to problems that cannot be characterized as minimizers of an energy functional.

\vskip5pt
\textbf{\em An advantage of weak forms over strong forms for singular kernels.}
For singular kernels $\g(\xb,\yb)$, i.e., for kernels such that $\g(\xb,\yb)\to\infty$ as $\yb\to\xb$, the integral in the definition \eqref{oper} of the operator $\mcL$ has to be interpreted in the principal-value sense. As a result, discretization of \eqref{prob} requires the use of very carefully designed quadrature rules. For the weak formulation \eqref{weak}, the first term in the bilinear form $A(u,v)$ defined in \eqref{bform2} also has a problematic integrand if the kernel is singular. However, dealing with approximations of that term is less troublesome compared to dealing with approximations of \eqref{oper}. Heuristically, both \eqref{oper} and the first term in \eqref{bform2} have to deal with a $\frac00$ for $\yb=\xb$. The zero in the denominator is the ``same'' for both cases. However, the zero in the numerator is ``stronger'' for \eqref{bform2} because it involves a double integration and the quadratic mollifying contribution $\big(u(\yb)-u(\xb)\big)\big(v(\yb)-v(\xb)\big)$ to the integrand whereas \eqref{oper} involves a single integral and a linear mollifying contribution $\big(u(\xb)-u(\yb)\big)$ to the integrand.

\section{Finite element discretization}\label{fem}

In this section we consider finite element discretizations of the weak formulation \eqref{weak} using general piecewise-polynomial bases defined with respect to a grid. However, in the remaining sections, we focus on piecewise-linear bases and only remark, in Sec. \ref{conclusion}, about extensions to higher-order piecewise-polynomial bases.

Finite element methods for nonlocal volume-constrained problems have been studied using continuous and discontinuous piecewise-linear finite element spaces and discontinuous piecewise-constant finite element spaces; see, e.g., {Refs. \refcite{xchen,TiDu13,TiDu14,feifei2,feifei1}.} These approaches have been tested on manufactured smooth solutions (e.g., polynomial solutions). If $\delta>h$, all the approaches perform well, whereas the piecewise-linear finite element spaces, {both continuous and discontinuous}, are more robust if $\delta<h$ in the sense that optimal accuracy with respect to $h$ is again obtained whereas piecewise-constant approximations fail to do so. 

As stated in Sec. \ref{sec:intro}, the central goals of this paper are dealing with difficulties arising from choosing, as is ubiquitous, the Euclidean ball $B_\del(\xb)$ as the interaction set corresponding to a point $\xb$ and also with the selection of quadrature rules that do not compromise the accuracy of finite element approximations when used for approximating the double integrals appearing in the weak formulation \eqref{weak}. However, there are other challenges that can arise when using finite element methods for nonlocal problems. Because these challenges are not germane to our goals, we only consider them in brief remarks including those that follow here.

\vskip5pt
{\em\textbf{Singular kernels.}} 
A challenge arising in the assembly process occurs if singular kernels are involved; such kernels arise in several important applications such as fractional derivative models and the peridynamics model for solid mechanics. Singular kernels induce a need for the use of sophisticated numerical quadrature rules. The implementation becomes more demanding and additional computational costs may arise. See, e.g., Ref. \refcite{DElia-ACTA-2020} for further discussions about this issue.

\vskip5pt
{\em\textbf{Solutions with jump discontinuities.}} 
Solutions with jump discontinuities are of interest because they arise in applications and because such solutions are not admissible for second-order elliptic PDE problems but are admissible for nonlocal problems with, e.g., translationally invariant integrable kernels $\g(\xb,\yb)$. All types of finite element discretizations, be they continuous or discontinuous or be they piecewise constant or linear, loose accuracy in the presence of discontinuities. For example, if one uses a uniform grid of size $h$ and piecewise-polynomial finite element spaces of any degree, in general, the best accuracy that can be achieved in the $L^2$-norm of the error is of $\mcO(h^{1/2})$; the $L^\infty$-norm of the error could be of $\mcO(1)$. However, unlike the other choices, the accuracy of discontinuous approximations can be improved by, e.g., {\em abrupt mesh refinement} near surfaces across which the solution is discontinuous. Note that near discontinuities, one would want $\delta>h$, a regime in which discontinuous finite element spaces perform optimally. For a more detailed discussion, see, e.g., Refs. \refcite{xchen,feifei2,feifei1}.  

\subsection{\textbf{Finite element grids and spaces}} \label{femgrid}

For the sake of simplicity of exposition, we assume that $\Omega$ is a polyhedral domain\footnote{Non-polyhedral domains can be handled by well-known methods documented in the finite element literature; see, e.g., Refs. \refcite{brenner,ciarlet}.}. Let ${\mcT}_{h,\Omega}$ denote a regular triangulation\footnote{We use the terminology ``triangulation'' to refer to general subdivisions of a domain, even if the domain is a subset of $\R$ or $\R^3$, and even if the subdomains are something other than triangles.} (see, e.g., Refs. \refcite{brenner,ciarlet}) of $\Omega$ into $K_\Omega$ finite elements $\{\mcE_k\}_{k=1}^{K_\Omega}$; we often refer to $\mcE_k$ as simply an element and in contexts for which the elements are indeed triangles, we will simply refer to them as triangles. Because $\Omega$ is a polyhedral domain, this triangulation is exact, i.e., $\cup_{k=1}^{K_\Omega} \mcE_k = \Omega$. As always, it is propitious to ensure that one ``triangulates into corners'', i.e., that {every vertex of $\Omega$ is also a vertex of the triangulation ${\mcT}_{h,\Omega}$.}

\begin{figure}[ht] 
\centerline{
\includegraphics[width=1in]{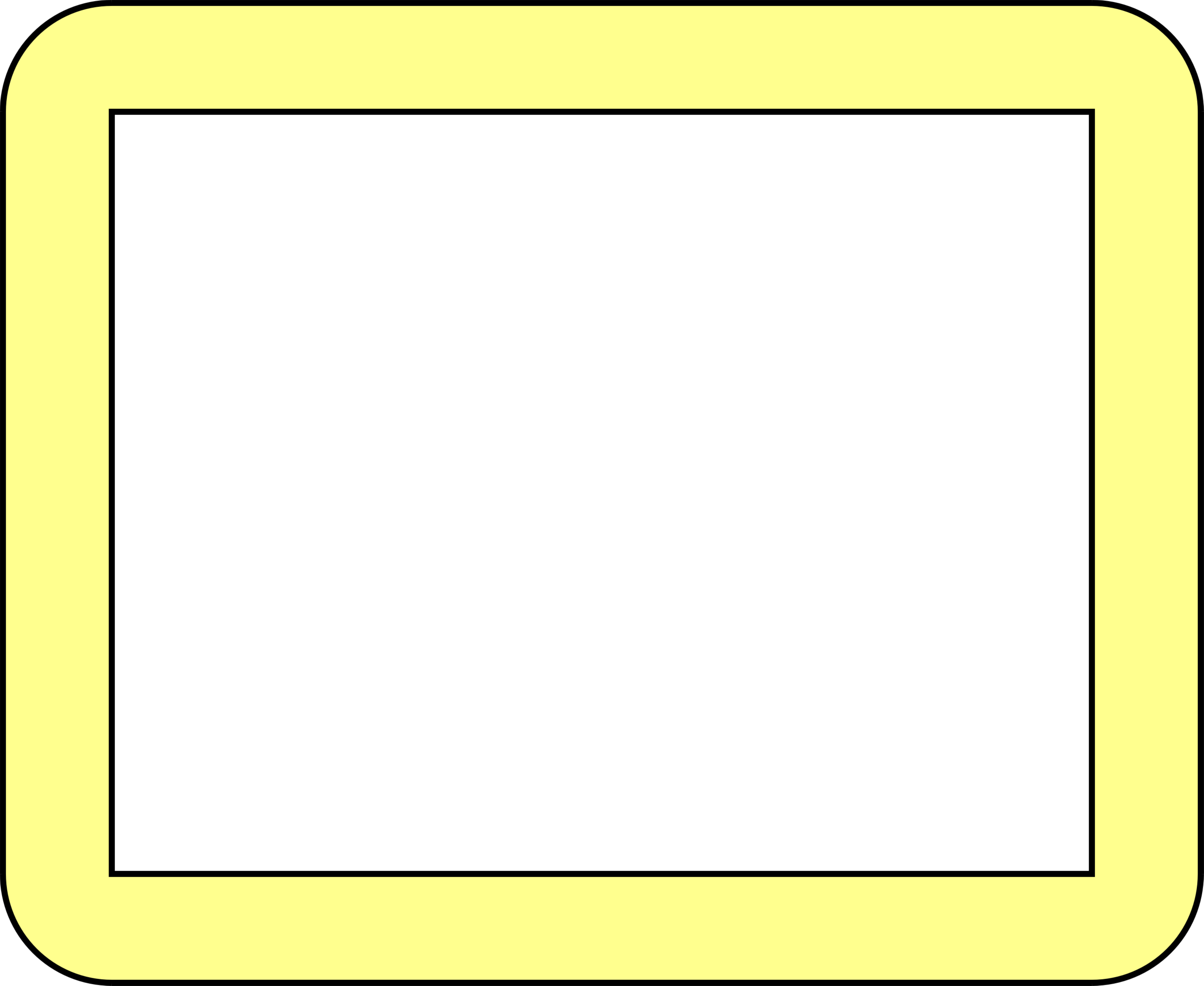}
\quad
\includegraphics[width=1in]{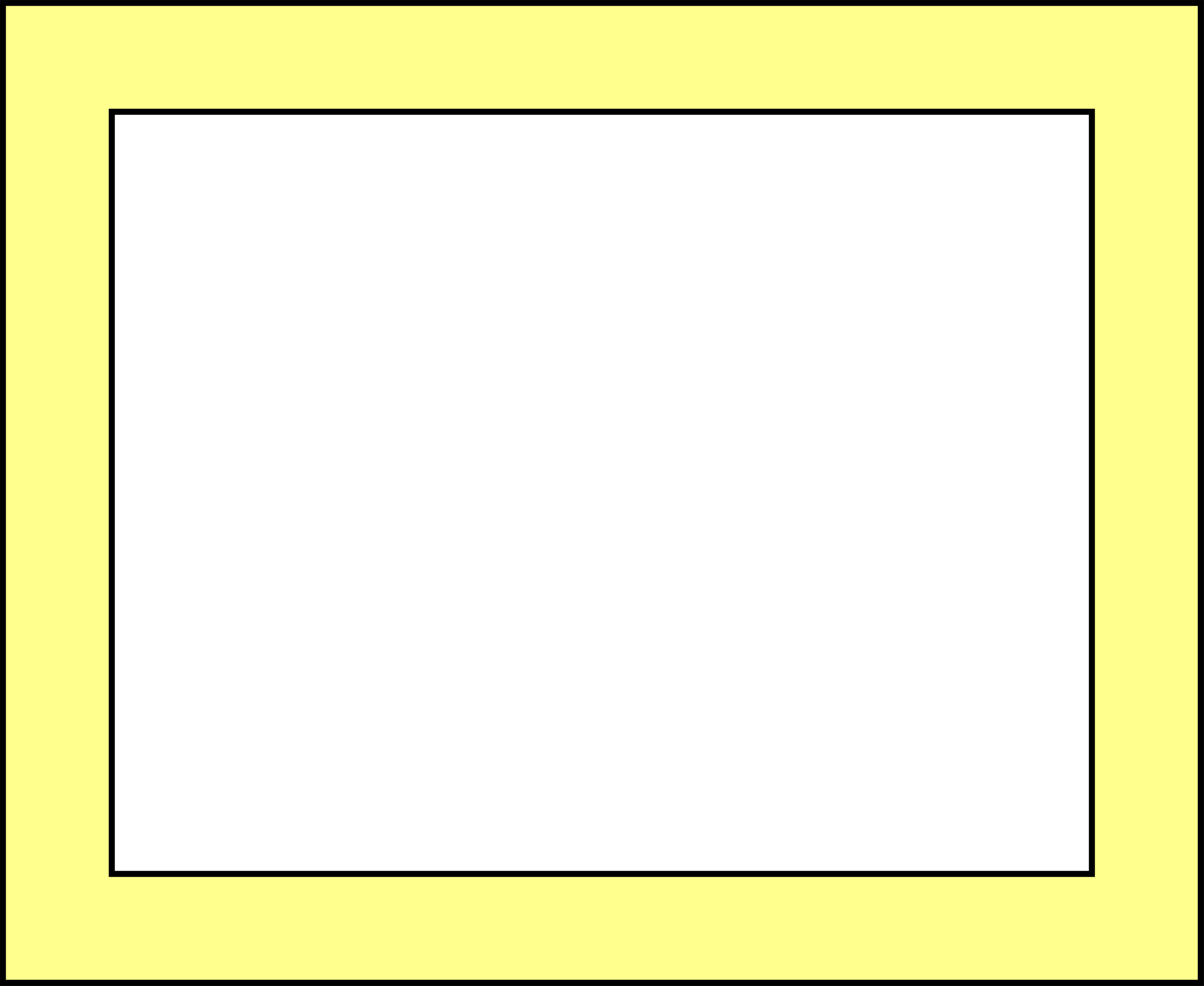}
}  
\caption{Left: a rectangular domain $\Omega$ (the white rectangle) and the corresponding interaction domain $\Omega_{\mathcal I}$ having rounded corners (in yellow). Right: the same rectangular domain and a polygonal approximate interaction domain, still referred to as $\Omega_{\mathcal I}$.}
\label{round-corners}
\end{figure}

For polyhedral $\Omega$, the corresponding interaction domain in case of Euclidean balls is in general not polyhedral, i.e., vertices of $\Omega$ cause rounded corners in $\Omega_{\mathcal I}$; see Fig. \ref{round-corners}-left for a simple illustration. As a result, $\Omega_{\mathcal I}$ cannot be exactly triangulated into elements with straight sides in two dimensions or with planar faces in three dimensions. Again, for the sake of simplicity of exposition, we approximate $\Omega_{\mathcal I}$ by a polyhedral domain by replacing rounded corners by vertices; see Fig. \ref{round-corners}-right for a simple illustration. We henceforth refer to that approximate domain also as $\Omega_{\mathcal I}$. No extension of the data $g(\xb)$ is needed because the added regions between the curved corners of the ``old'' $\Omega_{\mathcal I}$ and the polygonal corners of the ``new'' $\Omega_{\mathcal I}$ are never accessed during the finite element assembly process.

\begin{figure}[ht]
\centering
\begin{tabular}{cc} 
\includegraphics[width=1.25in]{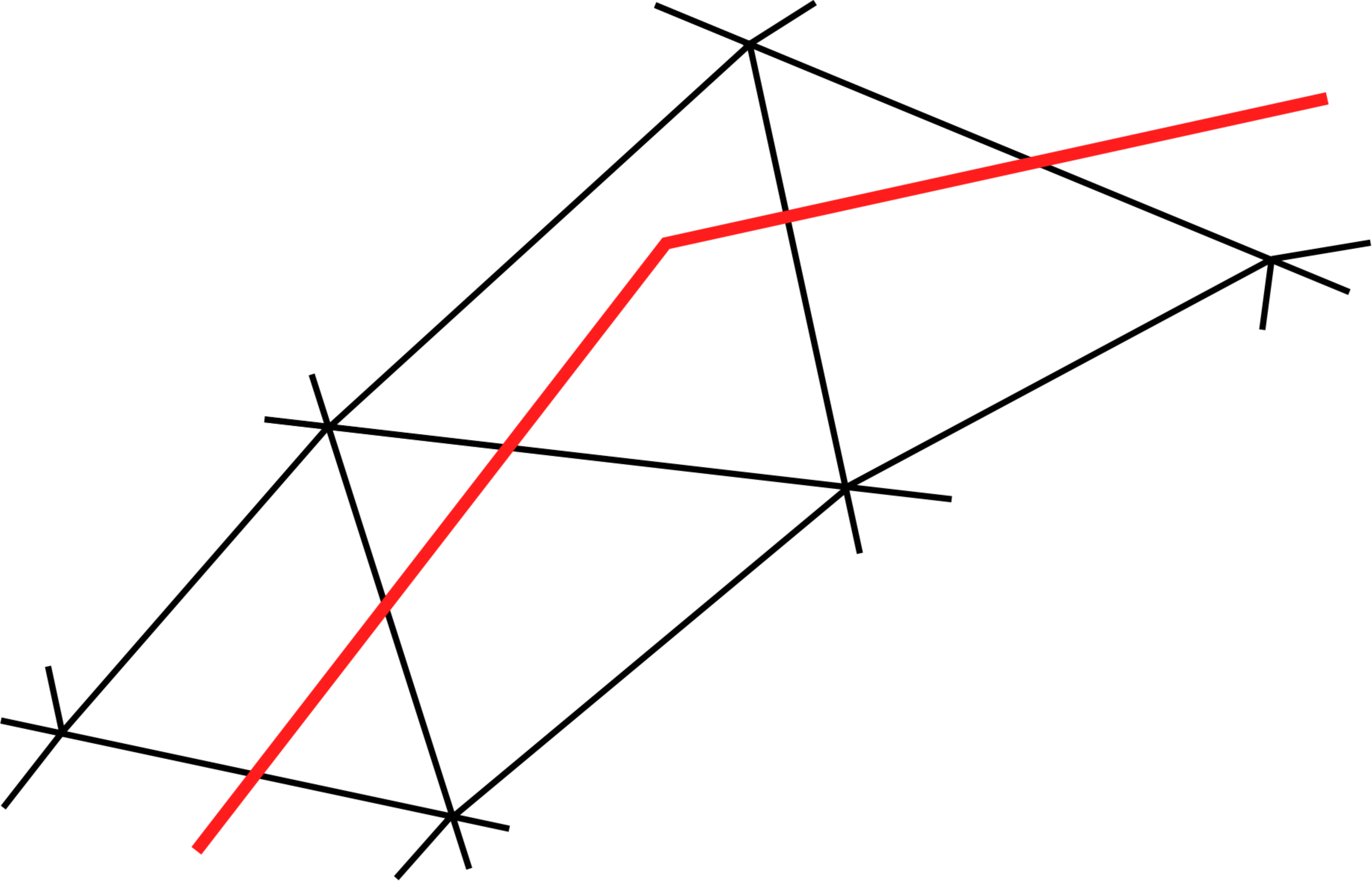}
&
\includegraphics[width=1.25in]{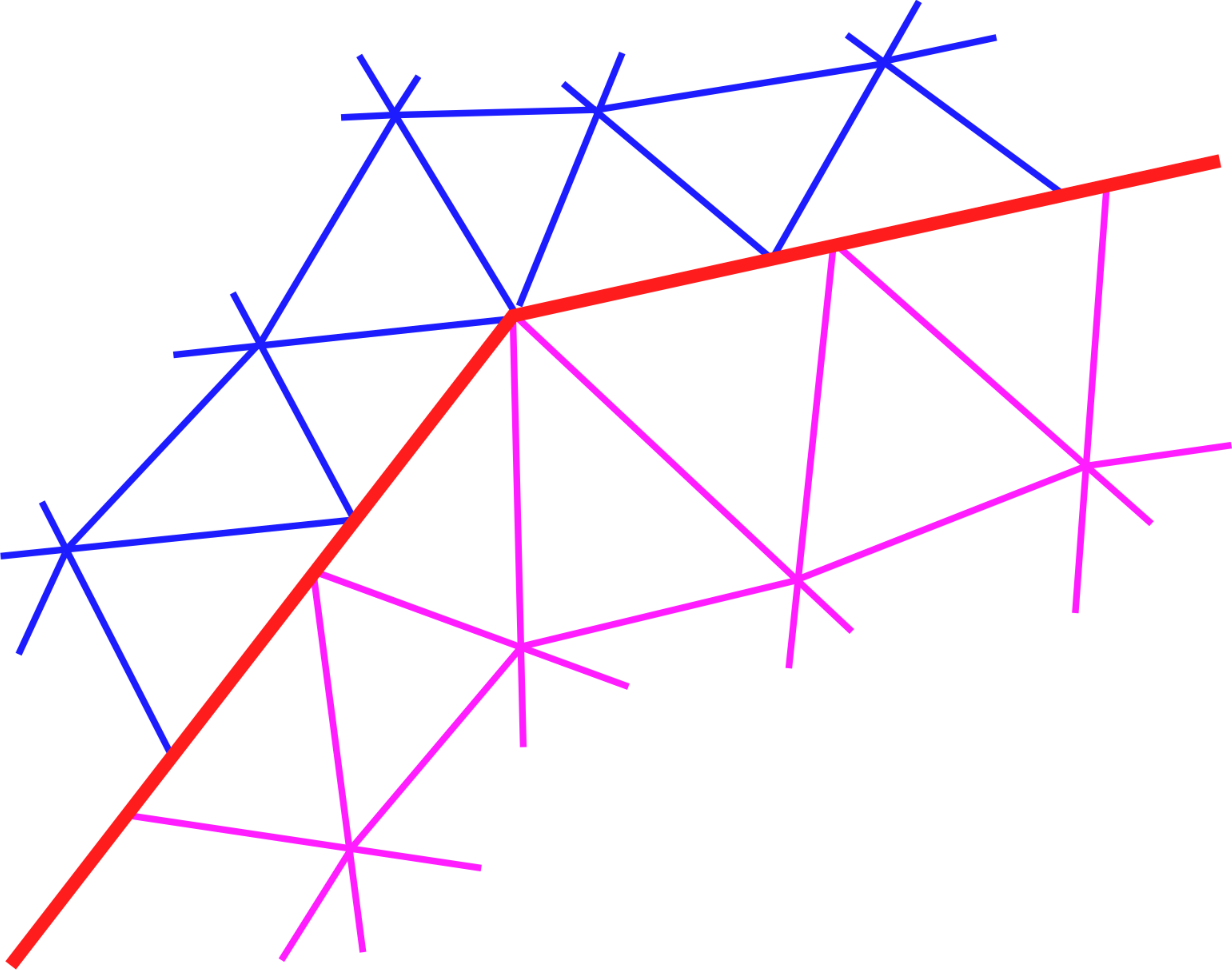}
\\
(a) & (b)
\\
\includegraphics[width=1.25in]{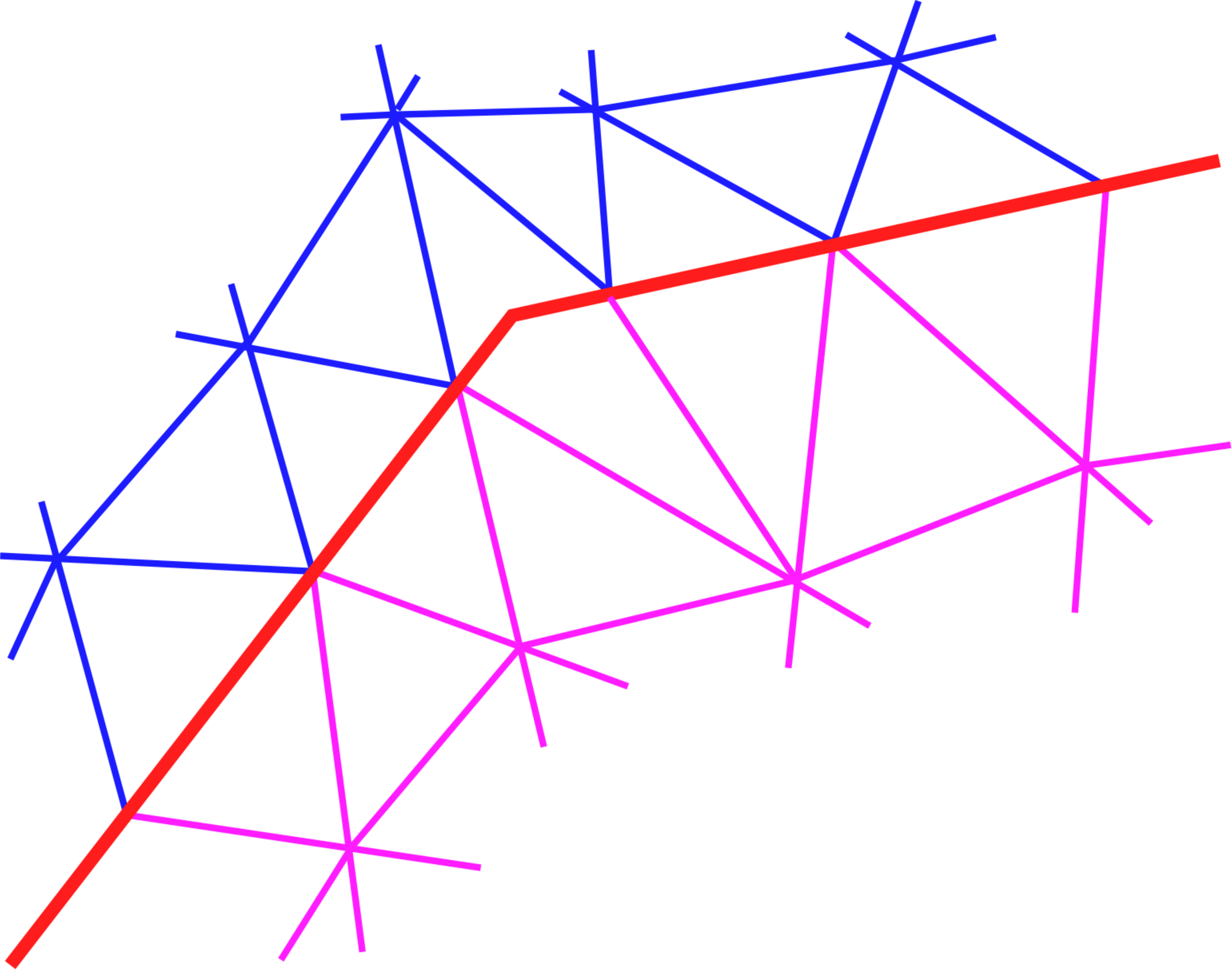}
&
\includegraphics[width=1.25in]{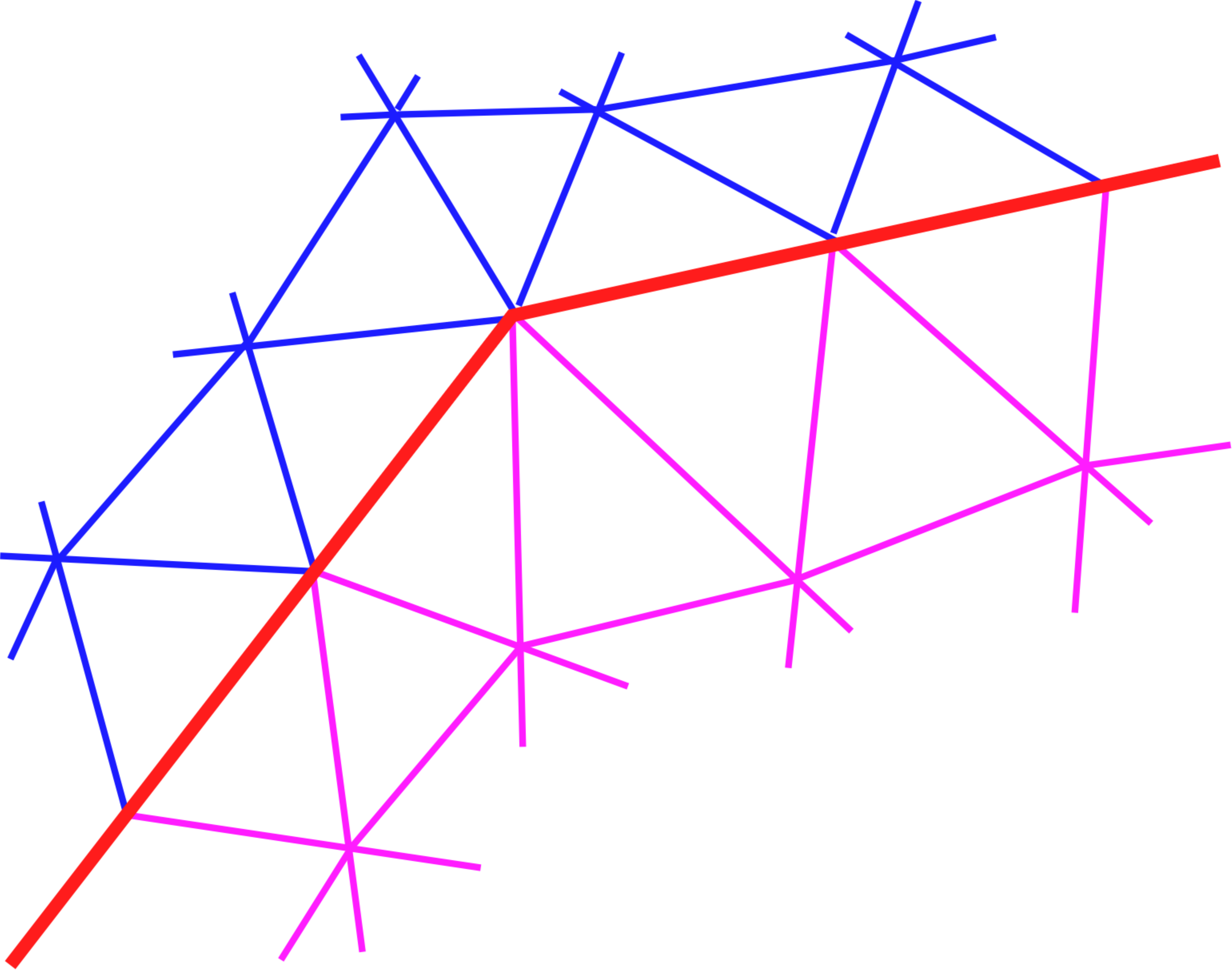}
\\
 (c) & (d)
\end{tabular} 
\caption{The red line segments are a part of the boundary of a polygonal domain $\Omega$. (a) The finite elements straddle across the boundary of $\Omega$. (b) The elements do not straddle across the boundary of $\Omega$ but the vertices that lie on that boundary corresponding to elements on opposite sides do not coincide. (c) The element vertices now coincide on the boundary of $\Omega$ but there is no triangle vertex located at a vertex of the boundary of $\Omega$. (d) A grid configuration that satisfies all requirements, namely, the elements do not straddle across the boundary of $\Omega$, a triangle vertex is placed at each vertex of that boundary, and the vertices on the boundary of elements on opposite side of that boundary coincide.}
\label{gridconstraints}
\end{figure}

Having assumed that $\Omega_{\mathcal I}$ is polyhedral, one can construct an exact regular triangulation $\mcT_{h,\Omega_{\mathcal I}}$ of $\Omega_{\mathcal I}$ into $K_{\Omega_{\mathcal I}}=K-K_{\Omega}$ finite elements $\{\mcE_k\}_{k=K_{\Omega}+1}^{K}$. Triangulating $\Omega$ and $\Omega_{\mathcal I}$ separately assures that {\em elements do not straddle across the common boundary of $\Omega$ and $\Omega_{\mathcal I}$}, i.e., across $\partial\Omega=\overline\Omega\cap\Omega_{\mathcal I}$, which is likely to occur if one directly triangulates $\Omega\cup\Omega_{\mathcal I}$. We require that the triangulations $\mcT_{h,\Omega_{\mathcal I}}$ and $\mcT_{h,\Omega}$ ``match'', i.e., that {\em along the boundary of $\Omega$, the vertices of the triangulations $\mcT_{h,\Omega}$ and $\mcT_{h,\Omega_{\mathcal I}}$ coincide.} In this case, the triangulation $\mcT_{h}=\mcT_{h,\Omega}\cup\mcT_{h,\Omega_{\mathcal I}}$ is itself a regular triangulation of $\OuO$ into $K$ elements $\{\mcE_k\}_{k=1}^{K}$. The constraints imposed on the triangulation $\mcT_{h}$ and the violations of those constraints are illustrated in Fig. \ref{gridconstraints}.

We restrict ourselves to continuous finite element spaces; discontinuous finite element spaces are also in use for discretizing nonlocal problems. However, the choice between the two types of spaces is, once again, not germane to the main goals of the paper; furthermore, arguments similar to those used in the following sections lead to the same conclusions for discontinuous finite element methods. We also restrict ourselves to Lagrange-type compactly supported piecewise-polynomial finite element bases that are defined with respect to a set of nodes associated with the triangulation $\mcT_h=\mcT_{h,\Omega}\cup\mcT_{h,\Omega_{\mathcal I}}$ of $\OuO$. For piecewise-linear and piecewise-bilinear bases, the associated nodes are merely the vertices of the elements, whereas for higher-degree polynomial bases, nodes placed on the edges or faces or even in the interior of the elements also come into play. Specifically, let $\{\widetilde\xb_j\}_{j=1}^J$ denote the set of nodes, with the nodes $\{\widetilde\xb_j\}_{j=1}^{J_\Omega}$ located in the open domain $\Omega$ and the nodes $\{\widetilde\xb_j\}_{j=J_\Omega+1}^{J}$ located in the closed domain ${\Omega}_{\mathcal I}$ so that the nodes located on $\partial\Omega={\overline\Omega}\cap{\Omega}_{\mathcal I}$ are assigned to ${\Omega}_{\mathcal I}$. Then, for $j=1,\ldots,J$, let $\phi_j(\xb)$ denote a piecewise-polynomial function such that $\phi_j(\widetilde\xb_{j'})=\delta_{jj'}$ for $j'=1,\ldots,J$, where $\delta_{jj'}$ denotes the Kroenecker delta function. We then define the finite element spaces 
$$
V^h = \mbox{span} \{ \phi_j(\xb) \}_{j=1}^J \subset V(\OuO)
\quad\mbox{and}\quad
V^h_c = \mbox{span} \{ \phi_j(\xb) \}_{j=1}^{J_\Omega}  \subset V_c(\OuO)
$$ 
of dimension $J$ and $J_\Omega$, respectively. By construction, functions belonging to $V^h$ and $V^h_c$ are continuous.

\subsection{\textbf{Finite element discretization of the weak formulations}}\label{femweak}

The finite element approximation $u_h\in V^h$ of the solution $u(\xb)$ of the nonlocal problem \eqref{weak} is determined as the solution of the discrete weak formulation
\beq{weakha}
    A(u_h,v_h)= F(v_h) \quad\forall\, v_h\in V_c^h.
\eeq
Here, the finite element approximation $u_h(\xb)$ has the form
\beq{appr}
  u_h(\xb) = 
  \sum_{j=1}^{J} U_j \phi_j(\xb) = \sum_{j=1}^{J_\Omega} U_j \phi_j(\xb) + \sum_{j=J_\Omega + 1}^J g(\widetilde\xb_j) \phi_j(\xb) \in V^h
\eeq
for a set of constants $\{U_j\}_{j=1}^J$, where the volume constraint in \eqref{prob} has been applied to set
\beq{intcon}
U_j = g(\widetilde\xb_j)\quad\mbox{for $j=J_\Omega+1,\ldots,J$.}
\eeq
Note that the volume constraint is applied at the nodes in ${\Omega}_{\mathcal I}$ that include the nodes located on the boundary $\partial\Omega$ between $\Omega$ and $\Omega_{\mathcal I}$.

The last term in \eqref{appr} is merely the interpolant of $g(\xb)$ in the space $V^h\setminus V^h_c$ so that it requires $g(\xb)$ to be continuous on $\Omega_{\mathcal I}$. On the other hand, well posedness of the weak formulations \eqref{weak0} or \eqref{weak} only requires that $g(\xb)\in L^2(\Omega_{\mathcal I})$. If $g(\xb)$ is not sufficiently smooth to posses a well-defined interpolant in $V^h\setminus V^h_c$, one can instead use, in \eqref{appr}, the $L^2(V^h\setminus V^h_c)$ projection of $g(\xb)$.

{Substituting \eqref{appr} into \eqref{weakha} and choosing $v_h(\xb)$ from the set of basis functions $\{\phi_{j'}(\xb)\}_{j'=1}^{J_\Omega}$ results in the linear system 
\beq{weakja}
\sum_{j=1}^{J_\Omega} A(\phi_{j},\phi_{j'})U_j = F(\phi_{j'}) \quad  \mbox{for $j'=1,\ldots,J_\Omega$}
\eeq}
from which the coefficients $U_j$, $j=1,\ldots,J_\Omega$, in \eqref{appr} are determined, where we have that the entries of the $J_\Omega\times J_\Omega$ stiffness matrix are given by 
\beq{bform3}
\bal
A(\phi_j,\phi_{j'}) &= \sum_{k=1}^{K_\Omega}\int_{\mcE_k}\int_{\Omega\cap B_\del(\xb)} 
\big(\phi_j(\yb)-\phi_j(\xb)\big)\big(\phi_{j'}(\yb)-\phi_{j'}(\xb)\big)\psi(\xb,\yb)d\yb d\xb 
\\&\qquad\qquad+ 
2\sum_{k=1}^{K_\Omega}\int_{\mcE_k}\phi_j(\xb)\phi_{j'}(\xb)\Big(\int_{\Omega_{\mathcal I} \cap B_\del(\xb)} \psi(\xb,\yb)d\yb\Big) \,d\xb 
\eal
\eeq
for $j,j'=1,\ldots,J_\Omega$, and the components of the $J_\Omega$-dimensional right-hand side vector are given by
\beq{lfunc3}
F(\phi_{j'}) = \sum_{k=1}^{K_\Omega}\int_{\mcE_k}\phi_{j'}(\xb) \Big(f(\xb) + 
2\int_{\Omega_{\mathcal I}\cap B_\del(\xb)} g(\yb)\psi(\xb,\yb)d\yb  
  \Big) \,d\xb 
\eeq
for $j'=1,\ldots,J_\Omega$. In \eqref{bform3} and \eqref{lfunc3} we have expressed the integrals over $\Omega$ as the sum of integrals over the sets of finite elements ${\mcT}_{h,\Omega}$ that cover $\Omega$. Also, because in \eqref{inter} we assumed that a point $\xb\in\Omega$ interacts only with the points $\yb\in B_\del(\xb)$, we  restricted the domain of integration of the inner integrals in \eqref{bform3} and \eqref{lfunc3} to the ball $B_\del(\xb)$. Also note that even for singular kernel functions $\psi(\xb,\yb)$, i.e., for kernel functions such that $\psi(\xb,\yb)\to\infty$ as $\yb\to\xb$, the inner integrals in the second term in \eqref{bform3} and in \eqref{lfunc3} are bounded because $\xb\in\Omega$ and $\yb\in\Omega_{\mathcal I}$, although some care must be exercised whenever $\xb\in\Omega$ and $\yb\in\Omega_{\mathcal I}$ are both close to the same point on the boundary of $\Omega$.

\subsection{\textbf{Estimate for the \em approximation \em error incurred by finite element discretization}}\label{feme}

Let the finite element space $V^h$ be the space of functions in $V(\OuO)$ that are piecewise polynomials of degree no more than $m$ defined with respect to the 
shape-regular triangulation $\mcT_{h}=\mcT_{h,\Omega}\cup\mcT_{h,\Omega_{\mathcal I}}$. If the exact solution is sufficiently smooth, we have the following result; see Ref. \refcite{dglz1}.

\begin{theo}[Approximation error due to finite element discretization]\label{thm:FEMaccuracy}
{Assume that the kernel $\gamma(\xb,\yb)$ in \eqref{prob} is square integrable or translationally invariant and integrable so that the energy space $V(\OuO)$ is equivalent to $L^2(\OuO)$.} Let $m$ denote a nonnegative integer and suppose that the domain $\Omega$ and the data $f(\xb)$ and $g(\xb)$ are such that $u(\xb)|_{\Omega}$ belongs to the Sobolev space {$H^{m+1}(\Omega)$}. Then, there exists a constant $C$ whose value is independent of $h$, $\delta$, and $u$ such that, for sufficiently small $h$,
\beq{eq:FEMaccuracy}
\|u-u^h\|_{L^2(\OuO)} \leq C {h^{m+1}} \|u\|_{H^{m+1}(\Omega)}.
\eeq
\end{theo}
In the case of piecewise-linear polynomials, i.e., $m=1$, \eqref{eq:FEMaccuracy} implies that the expected optimal convergence rate is quadratic, i.e., $\|u-u^h\|_{L^2(\OuO)}=\mathcal O(h^2)$. This result plays a fundamental role in the choice of quadrature rules for the outer and inner integrals and of approximations of the standard Euclidean balls. In the following three sections, we examine how the convergence rates for the errors introduced by such approximations compare with that of \eqref{eq:FEMaccuracy}. Of course, the ideal situation is the one in which those choices result in convergence rates that are commensurate with that of \eqref{eq:FEMaccuracy} so that the overall convergence rate remains optimal.

\vskip5pt
{\textbf{\em Finite element error estimates for non-integrable kernels.}}
{As has already been stated, we note that in this paper we limit ourselves to the case of square integrable kernel functions or translationally invariant integrable kernels so that we can refer to \eqref{eq:FEMaccuracy} whenever discussing convergence rates. For non-integrable kernels, $L^2$-norm error estimates are generally not available. Instead, if the energy space is a strict subspace of $L^2(\OuO)$, error estimates are only available with respect to the corresponding energy norm; see Ref. \refcite{dglz1}.}

\section{Approximate balls}\label{approxballs}

{As mentioned in Sec. \ref{sec:intro}, there are difficulties encountered in the finite element assembly process, difficulties that result from the use of Euclidean balls as interaction domains.} To alleviate these difficulties and thus simplify the assembly process and make it more efficient, in this section we define approximations $B_{\del,h}(\xb)$ of the Euclidean ball $B_{\del}(\xb)$ that appears in the domain of integration of the inner integrals in \eqref{bform3} and \eqref{lfunc3}.

To keep the exposition relatively simple, in this section we only consider the two-dimensional case and triangular meshes. Quadrilateral meshes can be handled in the same manner as triangular meshes; in fact, their treatment is, in many of the situations discussed in this section, simpler than it is for triangular meshes.

In Sec. \ref{approxballe} an estimate is given for the {\em geometric} error incurred as a result of using approximate balls $B_{\del,h}(\xb)$ instead of the true ball $B_{\del}(\xb)$. Then, in Sec. \ref{polyballs} we provide four specific examples of {\em polytopial} approximate balls $B_{\del,h}(\xb)$ and in Sec. \ref{subsec:shifted} we consider an approximation of the ball $B_{\del,h}(\xb)$ constructed by shifting the center of the ball. We apply the estimate of Sec. \ref{approxballe} to each of the five approximate balls discussed in Sections \ref{polyballs} and \ref{subsec:shifted}. 

\vskip5pt
{\em\textbf{Finite element discretization using approximate balls.}}
If an approximate ball $B_{\delta,h}(\xb)$ is used instead of the exact ball $B_{\del}(\xb)$, a finite element approximation $\wuh(\xb)$ is obtained from the system
\beq{D-Dh-weak}
A_h(\wuh,v_h) = F_h(v_h) \quad \forall\, v_h\in V^h_c,
\eeq
where, instead of \eqref{bform2} and \eqref{lfunc2}, we have the approximate bilinear form
\beq{bform3hh}
\bal
A_h(u,v) =& \int_{\Omega}\int_{\Omega\cap B_{\delta,h}(\xb)} 
            (u(\yb)-u(\xb))(v(\yb)-v(\xb))\psi(\xb,\yb)\,d\yb\,d\xb
           \\&\qquad + 2\int_{\Omega}u(\xb)v(\xb)\bigg(\int_{\Omega_{\mathcal I}\cap B_{\delta,h}}\psi(\xb,\yb)d\yb\bigg) \,d\xb
            \quad\forall\, u,v\in V_c^h
\eal
\eeq
and approximate linear functional
\beq{lfunc3hh}
F_h(v) = \int_{\Omega}v(\xb) \bigg(f(\xb) + 
2\int_{\Omega_{\mathcal I}\cap B_{\delta,h}} g(\yb)\psi(\xb,\yb)d\yb  
  \bigg) \,d\xb \quad\forall\, v\in V_c^h.
\eeq
The corresponding stiffness matrix entries, instead of \eqref{bform3} and \eqref{lfunc3}, are given by
\beq{bform3h}
\bal
&A_h(\phi_j,\phi_{j'}) 
\\&\qquad= \sum_{k=1}^{K_\Omega}\int_{\mcE_k}\int_{\Omega\cap B_{\del,h}(\xb)} 
\big(\phi_j(\yb)-\phi_j(\xb)\big)\big(\phi_{j'}(\yb)-\phi_{j'}(\xb)\big)\psi(\xb,\yb)d\yb d\xb 
\\&\qquad\qquad\qquad+ 
2\sum_{k=1}^{K_\Omega}\int_{\mcE_k}\phi_j(\xb)\phi_{j'}(\xb)\Big(\int_{\Omega_{\mathcal I} \cap B_{\del,h}(\xb)} \psi(\xb,\yb)d\yb\Big) \,d\xb 
\eal
\eeq 
for $j,j'=1,\ldots,J_\Omega$,
and the components of the right-hand side vector are given by
\beq{lfunc3h}
F_h(\phi_{j'}) = \sum_{k=1}^{K_\Omega}\int_{\mcE_k}\phi_{j'}(\xb) \Big(f(\xb) + 
2\int_{\Omega_{\mathcal I}\cap B_{\del,h}(\xb)} g(\yb)\psi(\xb,\yb)d\yb  
  \Big) \,d\xb 
\eeq
for $j'=1,\ldots,J_\Omega$.

\subsection{\textbf{Estimates for the \em geometric \em error incurred by using approximate balls}}\label{approxballe}

In this section we provide general results about the error incurred as a result of the use of approximate balls. In the following proposition we show that the energy norm of $(u_h-\wuh)$ can be bounded by the volume of the symmetric difference between $B_\delta$ and $B_{\delta,h}$, i.e., by the volume $|\Delta B_{\delta,h}|$ of the set  of $\Delta B_{\delta,h}=(B_\delta\setminus(B_\delta\cap B_{\delta,h}))\cup(B_{\delta,h}\setminus(B_\delta\cap B_{\delta,h}))$. We refer $\Delta B_{\delta,h}(\xb)$ as the {\em ball difference}. We assume that for all $\xb\in\OuO$, the kernel function $\psi(\xb,\yb)$ is bounded for all $\yb\in\Delta B_{\delta,h}(\xb)$. This is generally true because, e.g., for singular kernels, the singular point is at the center of the ball $B_\delta(\xb)$ and, in general, is also in $B_{\delta,h}(\xb)$ so that it is not in $\Delta B_{\delta,h}(\xb)$.

The following proposition provides an error estimate for the energy norm of $(u_h-\wuh)$; the proof is given in \ref{appendixa}. The convergence rate with respect to $h$ of the energy norm of $(u_h-\wuh)$ determines whether or not the approximate balls introduced in Sections \ref{polyballs} and \ref{subsec:shifted} compromise the overall accuracy of the finite element approximations. 

\begin{prop}[Geometric error due to the use of approximate balls]\label{prop:error1}
Let $B_\delta(\xb)$ denote the $\ell^2$-ball and $B_{\delta,h}(\xb)$ be an approximation of that ball, and let $u_h$ and $\wuh$ denote the corresponding finite element solutions obtained from \eqref{weakha} and \eqref{D-Dh-weak}, respectively. Assume that for all $\xb\in\OuO$, the kernel function $\psi(\xb,\yb)$ is bounded for all $\yb\in\Delta B_{\delta,h}(\xb)$ and also that all inner and outer integrals in \eqref{weakha} and \eqref{D-Dh-weak} are exactly evaluated. Then,
\begin{equation}\label{eq:energy-diff1}
|||u_h-\wuh |||\leq K \;\sup_{\xb\in\Omega}|\Delta B_{\delta,h}(\xb)|,
\end{equation}
where $K$ is a positive constant that depends on the data $f$ and $g$ but is independent of $\delta$ and $h$.
\end{prop}
{The following corollary is immediate because of the equivalence between the norms $|||\cdot|||$ and $\|\cdot \|_{L^2(\OuO)}$ in the case of square integrable kernel functions or translationally invariant integrable kernels.}
\begin{cor}\label{prop:error}
Assume the hypotheses of Proposition \ref{prop:error1}. Also, assume that the kernel function $\psi(\xb,\yb)$ is square integrable or integrable and translationally invariant.
Then,
\begin{equation}\label{eq:energy-diff}
\|u_h-\wuh \|_{L^2(\OuO)}\leq C_e K \;\sup_{\xb\in\Omega
}|\Delta B_{\delta,h}(\xb)|,
\end{equation}
where $C_e$ denotes a norm-equivalence constant.
\end{cor}

As a consequence of Corollary \ref{prop:error}, for piecewise-linear finite element approximations, the (optimal) quadratic convergence rate is preserved as long as the ball difference has volume $\sup_{\xb\in\Omega
}|\Delta B_{\delta,h}(\xb)|\sim \mcO(h^r)$ with $r\geq 2$, provided the outer and inner integrals are sufficiently accurately approximated.

As already noted, the ball difference $\Delta B_{\delta,h}(\xb)$ does not contain the centers of the balls $B_\delta(\xb)$ or $B_{\delta,h}(\xb)$ so that even in the case of singular kernels, the proof of Proposition \ref{prop:error1} holds. Hence, the result \eqref{eq:energy-diff1} applies to singular kernels as well.

\subsection{\textbf{Polytopial approximate balls}}
\label{polyballs}

In this section, we consider four polytopial approximations of the ball $B_\delta(\xb)$. The construction process is based on the finite element grid $\mcT_h$ in the sense that in the two cases considered in Sections \ref{poly} and \ref{acaps}, the approximate balls consist of a subset of the finite element triangles and additional triangles each of which is itself a subset of a finite element triangle whereas in Sections \ref{bary} and \ref{overlapball}, the approximate balls consist of a subset of only the finite elements triangles.

The construction of the approximate polytopial balls we consider requires that at least some of the following tasks be executed, based on a given finite element mesh.

\hangthree 1.~Determination of the location of the barycenter of an element.

\hangthree 2.~Identification of elements that intersect the ball.

\hangthree 3.~Identification of those elements identified in 2. that are wholly contained within a ball.

\hangthree 4.~Identification of those elements identified in 2. that partially overlap with a ball.

\hangthree 5.~Identification of the points at which the boundary of the ball intersects the boundary of the elements.

\hangthree 6.~Determining a subdivision of a polygon into triangles.

\begin{figure}[ht]
\centerline{
\includegraphics[width=1.5in]{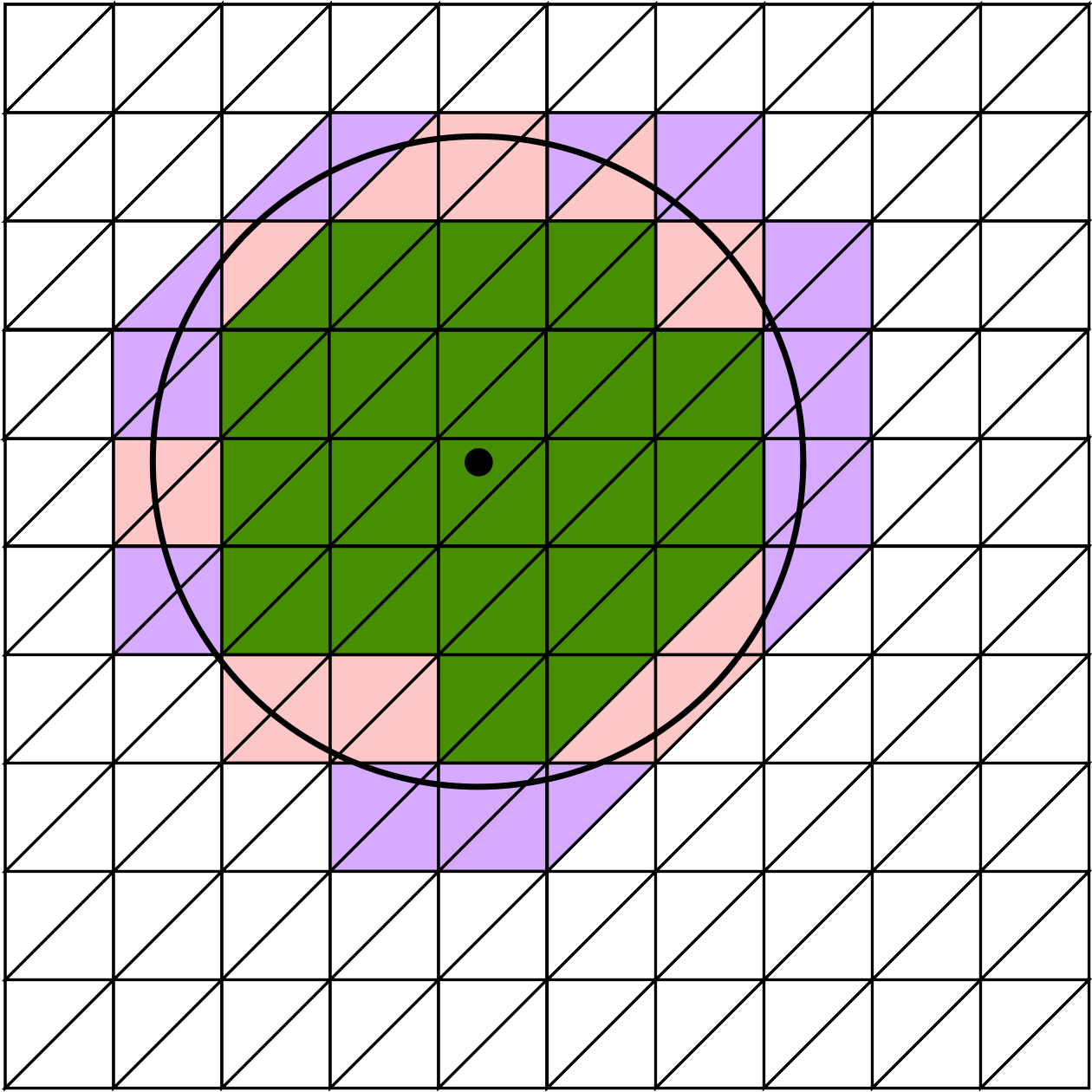}
}
\vskip10pt
\centering
\begin{tabular}{ll}
color of triangle $\mcE_{k}$ & type of triangle $\mcE_{k}$
\\
\hline
green & whole triangles intersecting the ball $B_\del(\xb)$, 
\\& i.e., $\mcE_k\cap B_\del(\xb)=\mcE_k$
\\
pink + magenta & partial triangles intersecting the ball $B_\del(\xb)$, 
\\& i.e., $\emptyset\ne\mcE_k\cap B_\del(\xb)\subsetneq\mcE_k$
\\
pink & partial triangles whose barycenters are
\\&  inside the ball $B_\del(\xb)$
\\
magenta & partial triangles whose barycenters are
\\&  outside the ball $B_\del(\xb)$
\\
white & whole triangles outside the ball $B_\del(\xb)$,
\\& i.e., $\mcE_k\cap B_\del(\xb)=\emptyset$
\\
\hline
\end{tabular} 
\caption{
The circle depicts the boundary of the ball $B_\del(\xb)$ of radius $\del$ centered at the black dot $\xb$. The colored triangles denote the elements $\mcE_{k}\in\mcT_h$ that overlap with the ball $B_\del(\xb)$. The color coding of the triangles depict the nature of the overlap, as listed in the table.
} 
\label{ftriangletypes}
\end{figure}

\noindent Efficient means for accomplishing these tasks are considered in Sec. \ref{ballconstr}. These tasks help to classify the finite elements into several categories, as illustrated in Fig. \ref{ftriangletypes}; this classification is used in the construction of the polytopial approximate balls. 

Suppose the black dot in Fig. \ref{ftriangletypes} is the center of the ball $B_\del(\xb)$. The colored triangles highlight all the triangles that overlap with the ball. Those triangles can be further categorized according to their geometric characteristics. Thus, we see both whole triangles and partial triangles intersecting the ball and differentiate between partial triangles whose barycenters are inside and outside the ball. 

\begin{figure}[ht]
\centering
\begin{tabular}{cc}
\includegraphics[height=1.45in]{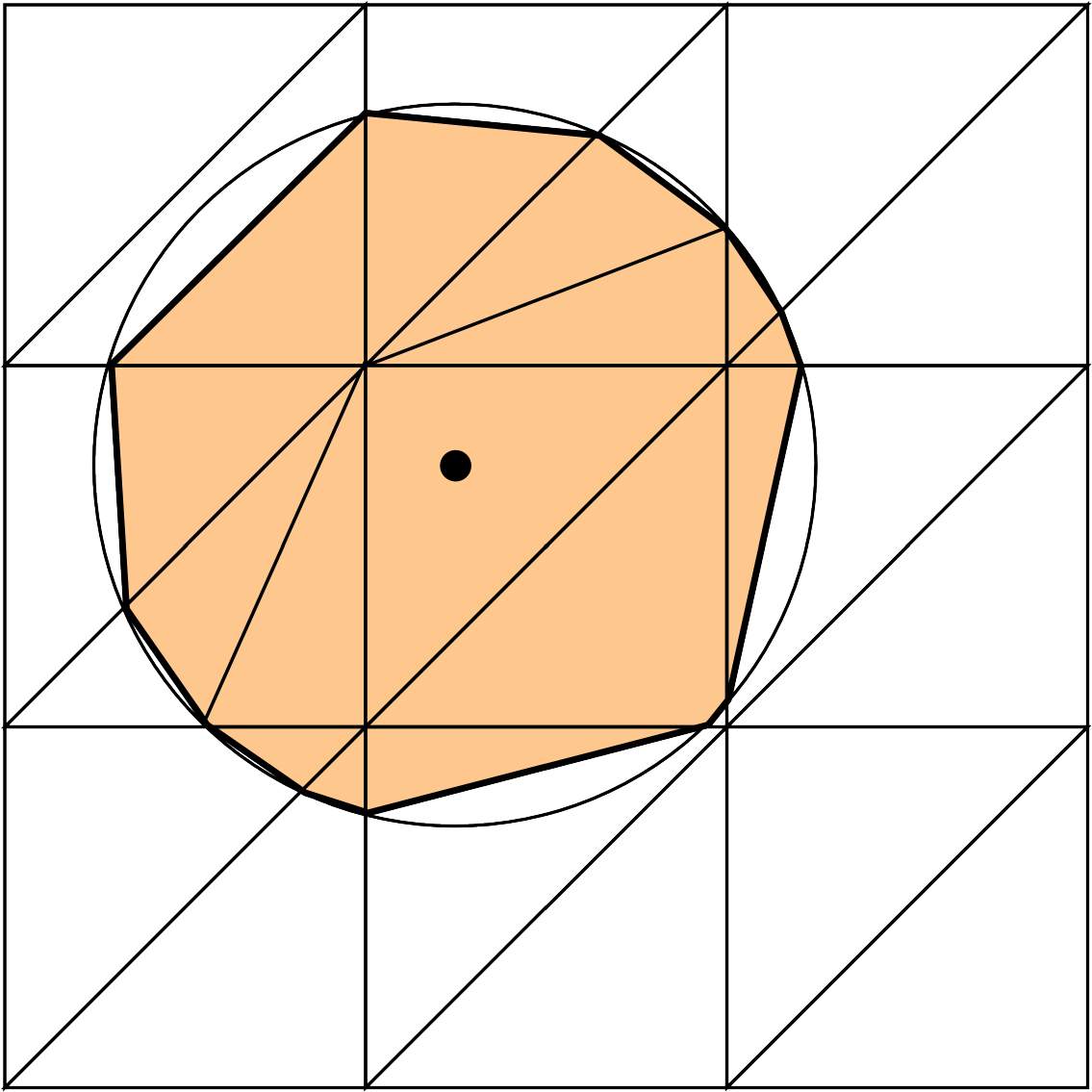} &
\includegraphics[height=1.45in]{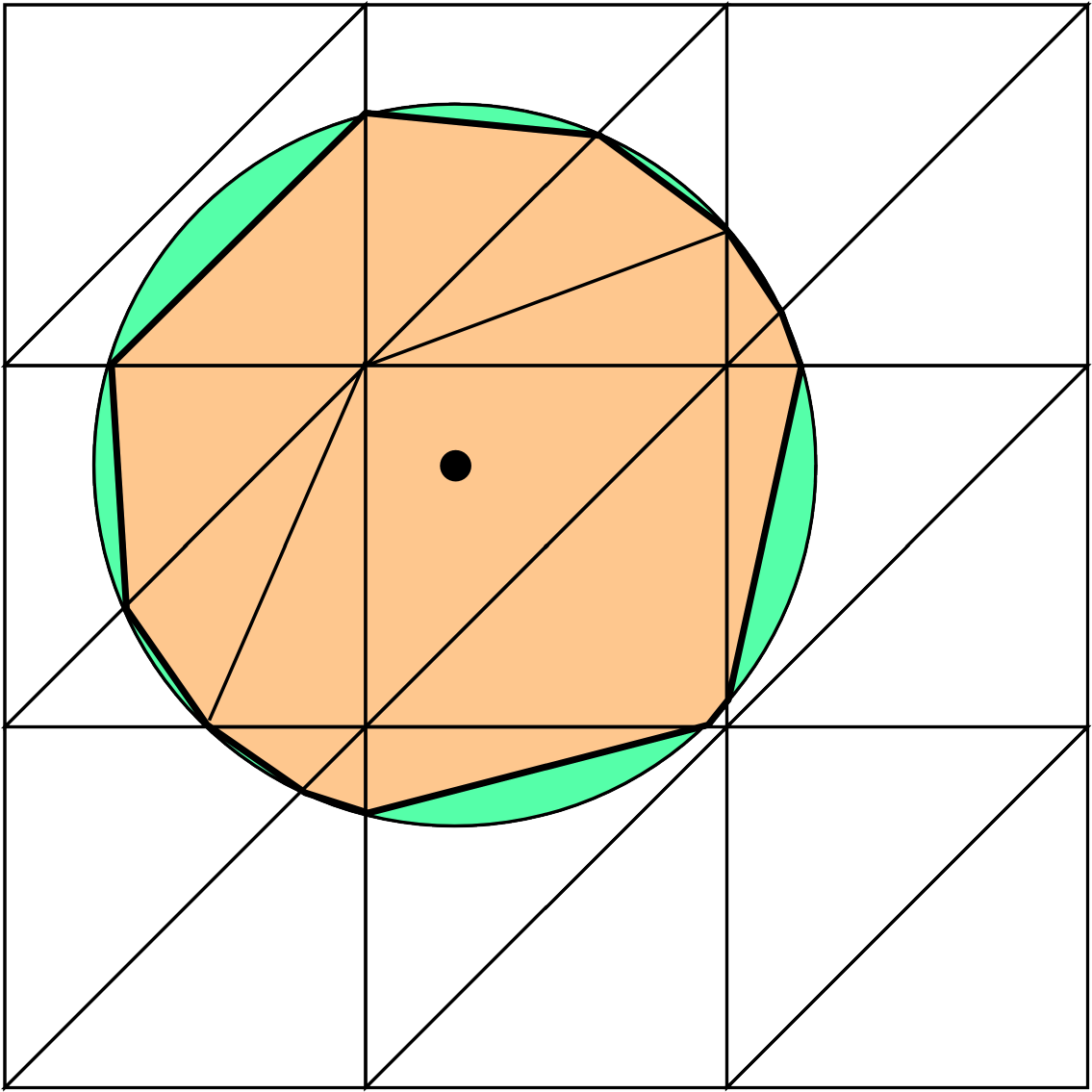} 
\\
(a) & (b) 
\\[1ex]
\includegraphics[height=1.45in]{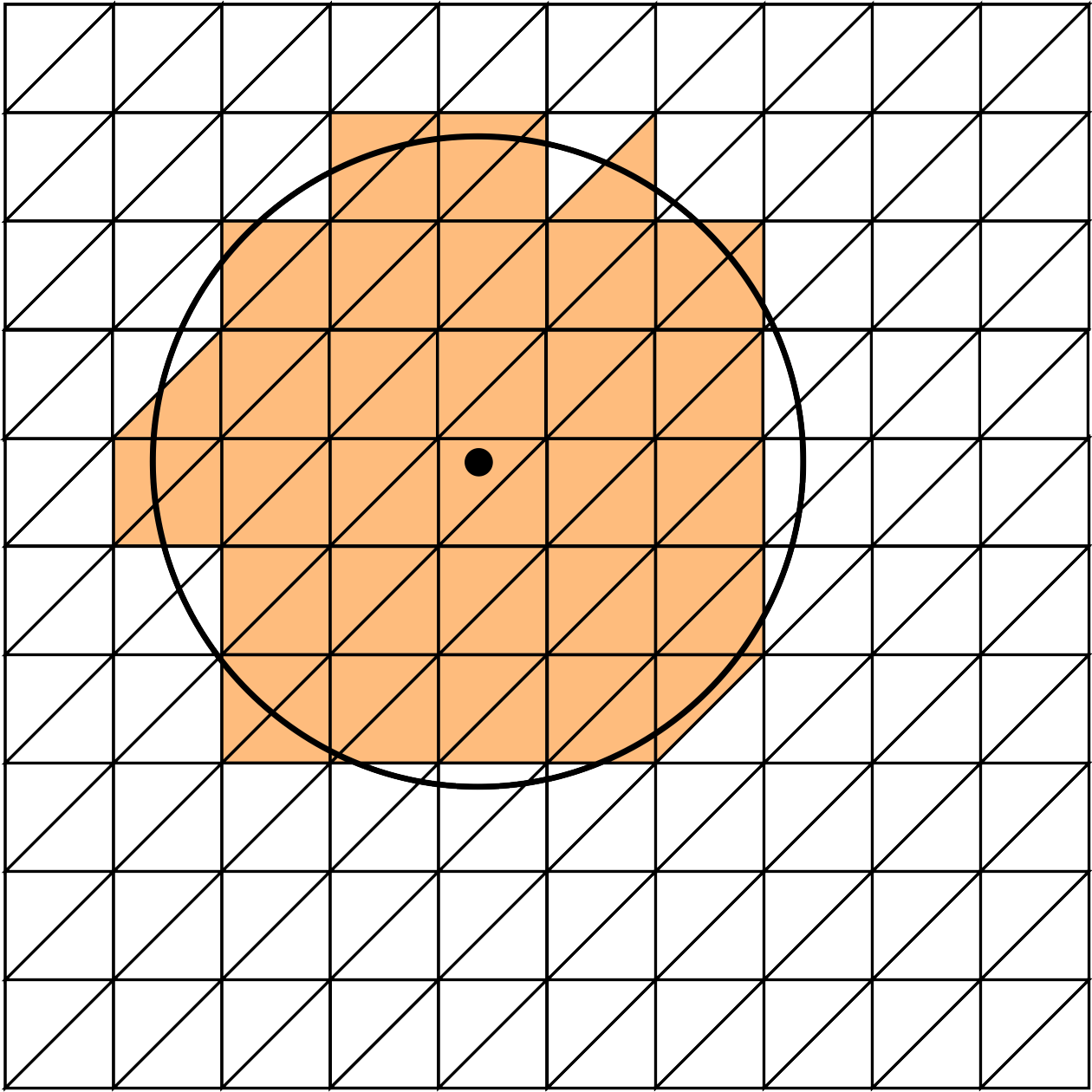} &
\includegraphics[height=1.45in]{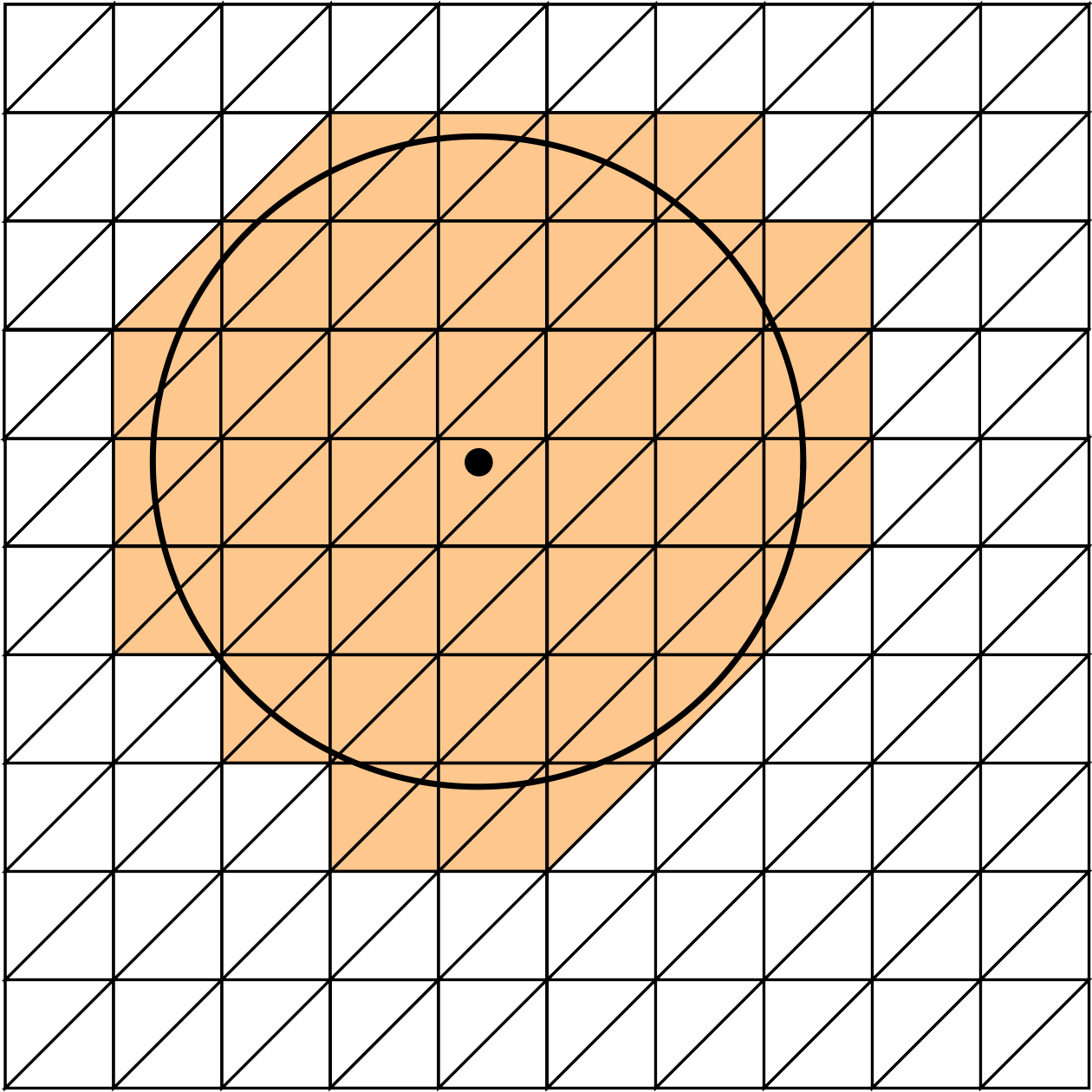}
\\
 (c) & (d)
\end{tabular} 
\caption{(a) An inscribed polygonal approximation of the ball. (b) Subdivision of the ball into the polygon of (a) and circular caps. An inscribed polygonal approximation of the ball is defined by approximating the green caps by triangles. (c) Approximation of the ball by whole finite element triangles that intersect the ball and for which the barycenter lies within the ball. (d) Approximation of the ball by all whole finite element triangles that intersect the ball.}
\label{fappoxball}
\end{figure}

In Sections \ref{poly} to \ref{overlapball} we provide specific examples of polytopial approximate balls and discuss how they are constructed and the geometric and solution errors incurred by replacing the exact ball by an approximate ball. The discussion makes use of the four geometric configurations depicted in Fig. \ref{fappoxball}. 
 
\subsubsection{\textbf{Inscribed triangle-based polygonal approximations of balls - Fig. \ref{fappoxball}a}}\label{poly}

The ball $B_\del(\xb)$ is approximated by an inscribed polygon $B^{nocaps}_{\del,h}(\xb)$ according to the following recipe.

\hangthree 1.~Determine the triangles ${\mathcal E}_k$ that are wholly contained within the ball, i.e., the triangles for which $\ {\mathcal E}_k\cap B_\del(\xb) = {\mathcal E}_k$.

\hangthree 2.~Determine the triangles ${\mathcal E}_k$ that are only partially contained within the ball, i.e., the triangles for which $\emptyset\ne  {\mathcal E}_k\cap B_\del(\xb) \ne {\mathcal E}_k$.

\hangthree 3.~For each triangle selected in step 2, determine the points at which the boundary of the ball intersects the sides of the triangle.
 
\hangthree 4.~Construct the polygon having vertices at the intersection points found in step 3. 

\noindent As a result of these steps, we have an inscribed polygon that is subdivided into triangles and polygons having more than three sides. For the latter we add one more step.

\hangthree 5.~Subdivide all polygons having more than three sides into triangles.

\noindent Fig. \ref{fappoxball}a illustrates the result of the five-step recipe. Note the two orange polygonal subregions that are divided into triangles. The sides of the polygon $B^{nocaps}_{\del,h}(\xb)$ so constructed are cords of the circular ball $B_\del(\xb)$ and, because they are necessarily shorter than the longest side of the triangle, the cords have lengths of ${\mathcal O}(h)$.

As a result of the five-step recipe, the approximate ball $B^{nocaps}_{\del,h}(\xb)$ is exactly subdivided into a set of nonoverlapping triangles ${\mcT}^{nocaps}_{\del,h,\xb}$ which consists of a subset of the finite element triangles in ${\mcT}_{\del,h}$ and also the triangles created by steps 2 to 5. For example, in Fig. \ref{fappoxball}a, $B^{nocaps}_{\del,h}(\xb)$ is subdivided into 14 triangles, only two of which are whole finite element triangles. Note that the membership of ${\mcT}^{nocaps}_{\del,h,\xb}$ depends on the horizon $\del$, the grid size $h$, and the position of the center of the exact ball $B_{\del}(\xb)$.

\begin{figure}[ht]
\begin{center}
\includegraphics[height=1.5in]{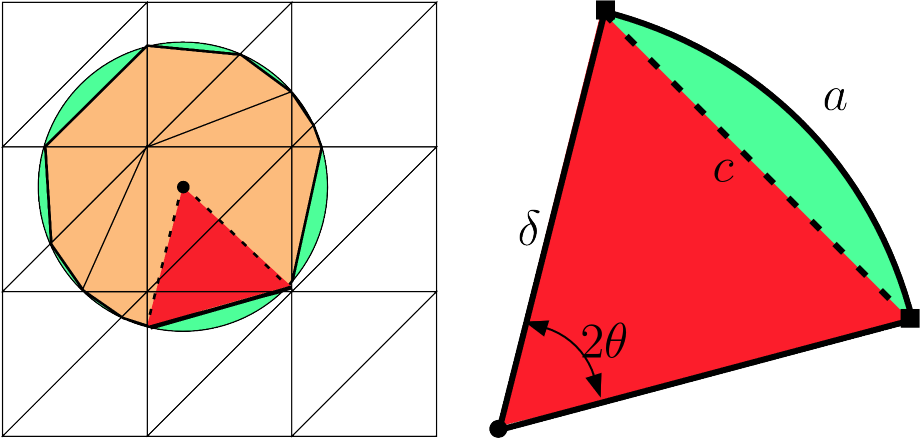}
\end{center}
\caption{Left: The red triangle and its abutting green circular cap depict one of the sectors defined by the center of the ball (the black dot) and a cord of the ball (the thick line segment). Right: A typical sector and the geometrical quantities used to estimate the area of the cap and the length of the arc.}
\label{caparea}
\end{figure}

\vskip5pt
{\em{\textbf{Geometric error.}}} A geometric error is incurred by replacing the ball $B_\del(\xb)$ by the polygon $B^{nocaps}_{\del,h}(\xb)$. Fig. \ref{caparea}-right highlights a typical sector of the ball $B_\del(\xb)$; such sectors are used to estimate the areas of a circular cap depicted in green.\footnote{What we refer to a ``circular caps'' or just ``caps'' are often referred to as ``circular segments.''} Circular caps are formed whenever the circular boundary of the ball intersects the sides of a triangle. The line segment joining the two intersection points is a cord of the circle and also a side of the polygon $B^{nocaps}_{\del,h}(\xb)$. In Fig. \ref{caparea}-left, we have 11 such triangles, hence there are 11 caps (highlighted in green) and 11 cords. 

The difference between the ball $B_\del(\xb)$ and its polygonal approximation $B_{\del,h}^{nocaps}(\xb)$ are the circular caps depicted in green in Fig. \ref{caparea}-left. To estimate the error associated with this approximation, according to Corollary \ref{prop:error}, we need to estimate the area of the ball difference $\Delta B_{\delta,h}^{nocaps}(\xb)$ for $B_{\delta,h}(\xb)=B_{\del,h}^{nocaps}(\xb)$; thus we need to estimate the areas of the caps and the number of caps. To do so, we consider a sector of the ball such as the one illustrated in Fig. \ref{caparea}-left by the red triangle and its abutting green cap. A typical sector is depicted in Fig. \ref{caparea}-right. The black dot denotes the center $\xb$ of the ball $B_\del(\xb)$ having radius $\delta$. The black squares denote the intersection points of the ball and the boundary of a triangle ${\mathcal E}_k$ of the grid. The dashed line connecting those two points is the cord $c$ that, along with the radius $\del$, defines the sector angle $2\theta$ and the circular arc $a$.

We first consider the case $h\ll\delta$ for which we have that

\hangfour~--~the length of the cord $c$, which we also denote by $c$, is smaller than the length of the longest side of the triangle ${\mathcal E}_k$ so that $c={\mathcal O}(h)$ and $\frac{c}{2\del}={\mathcal O}(\frac{h}\del)$

\hangfour~--~in terms of the radius $\delta$ and the cord length $c$, the 
\beq{areacap}
\mbox{area of a circular cap} =
\delta^2\bigg(\arcsin \Big(\frac{c}{2\delta}\Big) -  
\Big(\frac{c}{2\delta}\Big)\sqrt{ 1- \Big(\frac{c}{2\delta}\Big)^2 }\,\,\bigg) 
\eeq

\hangfour~--~if $h\ll\delta$, we easily see that area of a circular cap $ = {\mathcal O} (\frac{h^3}{\del})$.

\noindent We next estimate the number of sides of the polygon $B^{nocaps}_{\del,h}(\xb)$. We have that

\hangfour~--~$\sin\theta = \frac{c}{2\del}\le \frac{h}{2\del}$ so that for $h\ll\del$ we have  $\theta = {\mathcal O}(\frac{h}\del)$

\hangfour~--~the length of the circular arc $= \theta \del = {\mathcal O}(h)$

\hangfour~--~the perimeter of the circle is $2\pi\delta$;

\hangfour~--~therefore the number of circular arcs (= number of cords = the number of caps) is of ${\mathcal O}(\frac{\del}{h} )$.

\noindent Therefore, the total area of the circular caps $= {\mathcal O} (\frac{h^3}{\del})\, {\mathcal O}(\frac{\del}{h}) = {\mathcal O}(h^2)$. Clearly, we then have that the difference between the areas of the Euclidean ball and the inscribed polygon is estimated, for all $\xb \in \Omega$, by 
\beq{pdifare}
|\Delta B_{\delta,h}^{nocaps}(\xb)| = \mcO(h^2)\qquad \mbox{if $h\ll\del$}.
\eeq

In the mechanics setting, several authors set $\del=$ constant$\,\times\, h$; for example, in Refs. {\refcite{Bobaru12,Parks08}}, the choice $\del=3h$ is advocated. In such cases we have that 

--~the area of the ball is of ${\mathcal O}(h^2)$

--~the cord length $c={\mathcal O}(h)$ so that $\frac{c}\del = {\mathcal O}(1)$ 

--~the area of the cap is of ${\mathcal O}(\del^2)={\mathcal O}(h^2)$

--~the length of the circular arc is of ${\mathcal O}(h)$

--~the number of the circular arcs is of ${\mathcal O}(1)$

--~the total area of all of the circular caps is of ${\mathcal O}(\del^2)={\mathcal O}(h^2)$.

\noindent Thus, \eqref{pdifare} also holds for the case of $\del=$ constant$\,\times\, h$.

\vskip5pt
{\em{\textbf{Solution errors.}}}
Because of \eqref{pdifare}, according to Corollary \ref{prop:error} and, if the kernel is integrable and translationally invariant or just square integrable, we respectively have that, for piecewise-linear finite element approximations and for sufficiently smooth solutions,
\beq{polyerror}
  \|u_h-\wuh\|_{L^2(\OuO)}\leq C |||u_h-\wuh ||| \sim \mcO(h^2).
\eeq

\subsubsection{\textbf{Inscribed cap-based polygonal approximations of balls -- Figure \ref{fappoxball}b}}\label{acaps}

Given the results of Sec. \ref{poly}, it seems unnecessary to try to obtain a better polygonal approximation of a ball $B_\del(\xb)$. However, having such an approximation might be valuable. Although the ${\mathcal O}(h^2)$ accuracy in \eqref{polyerror} is good enough to preserve the second-order accuracy of the approximate solution, having a better approximation of the ball reduces the constant in the order relation.  

\begin{figure}[ht]
\centerline{
\includegraphics[height=1.4in]{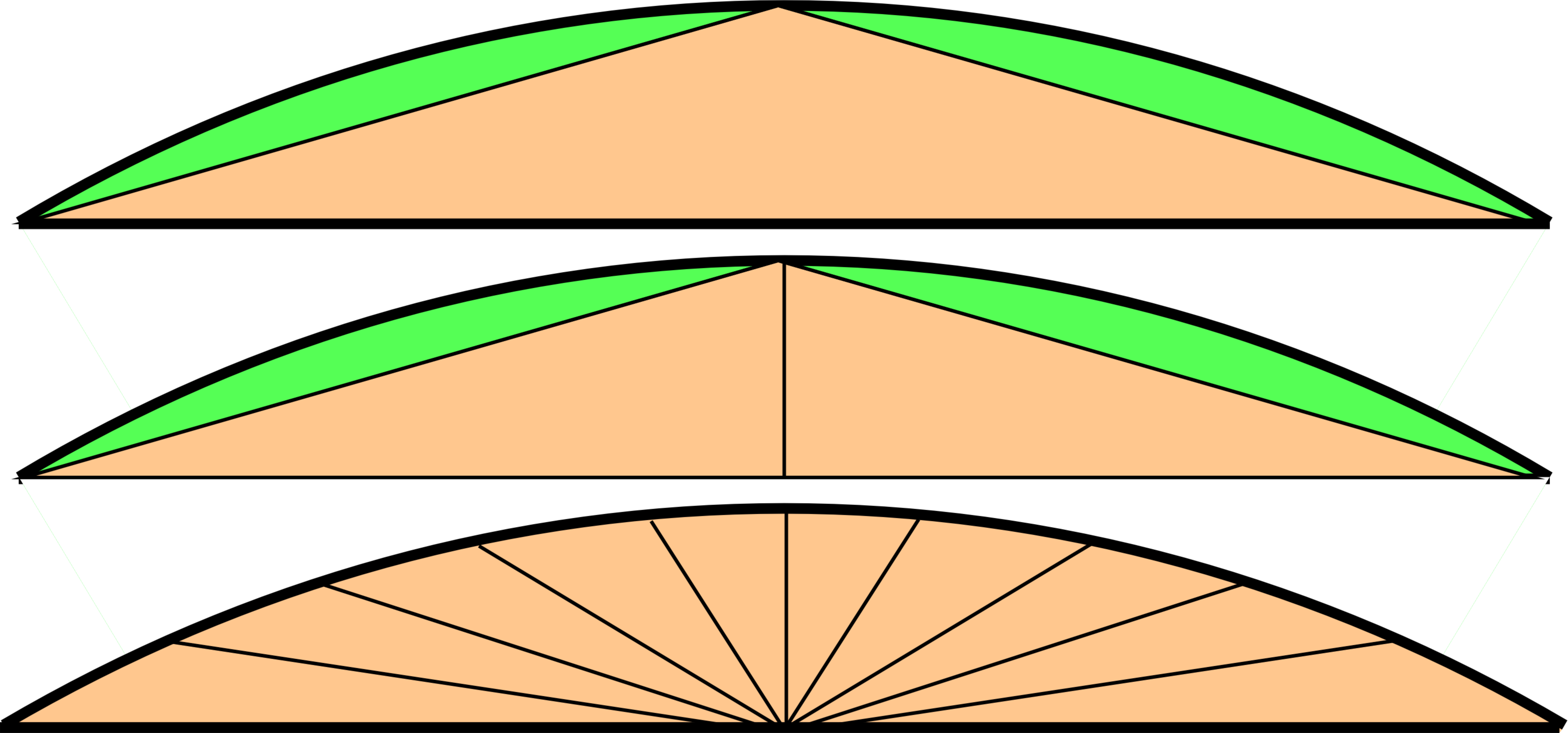}
\quad
\includegraphics[height=1.4in]{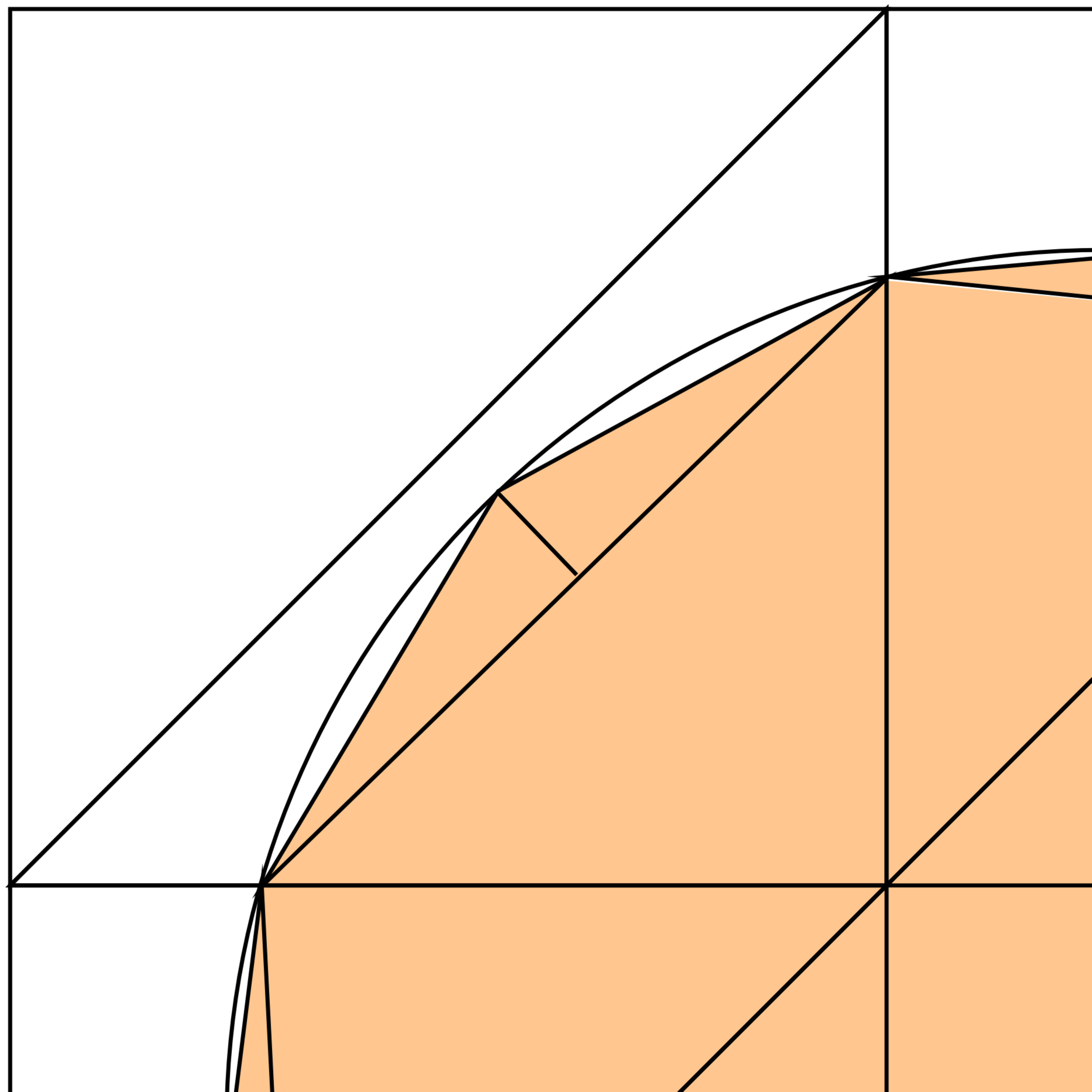}
}
\caption{Left: approximation of a cap by one, two, and ten triangles, where the latter is defined by dividing the circular arc into ten smaller arcs of equal arc length. Right: a zoom in illustrating that adding approximate caps to the ball $B^{nocaps}_{\del,h}(\xb)$ results in a better approximation of the exact ball $B_\del(\xb)$.}
\label{discap}
\end{figure}

In this section, we consider approximating the circular caps by triangles so that, together with the inscribed polygon $B^{nocaps}_{\del,h}(\xb)$ of Sec. \ref{poly}, there results in a different inscribed polygonal approximation $B^{approxcaps}_{\del,h}(\xb)$ of the ball. As is the case for $B^{nocaps}_{\del,h}(\xb)$, $B^{approxcaps}_{\del,h}(\xb)$ is subdivided into triangles. Fig. \ref{discap}-left illustrates a cap approximated by one, two, and ten triangles. With ten triangles one cannot, with the image resolution and image size used, see the part of the cap that lies outside of the triangles. Fig. \ref{discap}-right is a zoom-in illustrating how adding an approximate cap to the approximate ball $B^{nocaps}_{\del,h}(\xb)$ results in a better geometric approximation of the exact ball. In that figure, the large orange triangles (some of which are only partially depicted) are part of the approximate ball $B^{nocaps}_{\del,h}(\xb)$ whereas the two small orange triangles are what is added when forming the approximate ball $B^{approxcaps}_{\del,h}(\xb)$. 

Clearly the approximate ball $B^{approxcaps}_{\del,h}(\xb)$ is subdivided into a set ${\mcT}^{approxcaps}_{\del,h,\xb}$ of non-overlapping triangles consisting of the triangles in ${\mcT}^{nocaps}_{\del,h,\xb}$ plus the triangles added by approximating the caps. The membership of ${\mcT}^{approxcaps}_{\del,h,\xb}$ depends on the horizon $\del$, the grid size $h$, and the position of the center of the exact ball $B_{\del}(\xb)$.

\vskip5pt
{\em{\textbf{Geometric error.}}} 
Approximating each cap by one or a few triangles would not change the second-order convergence rate of the difference in the area $|\Delta B_{\delta,h}^{approxcaps}(\xb)|$ between $B_\del(\xb)$ and $B^{approxcaps}_{\del,h}(\xb)$, i.e., \eqref{pdifare} would hold for $B^{approxcaps}_{\del,h}(\xb)$ as well. However, the constant in the order relation is reduced. For example, consider the one or two triangle cases of Fig. \ref{discap}-left. We see that an omitted cap in the construction of $B_{\del,h}^{nocaps}(\xb)$ is replaced by triangles and two omitted smaller caps. The total areas omitted in the two cases are $\frac{\del^2}2 (2\theta - \sin 2\theta)$ and $\del^2 (\theta - \sin \theta)$, respectively, so that if $\theta\ll1$, i.e., if $h\ll\del$, it is easily seen that the constant in the order relation is reduced by a factor of four. Using more than two triangles to approximate a cap would reduce the constant even further, but would also incur additional costs.

\vskip5pt
{\em{\textbf{Solution errors.}}}
Because $|\Delta B_{\delta,h}^{approxcaps}(\xb)|=\mcO(h^2)$, the error estimates in \eqref{polyerror} also hold for $B^{approxcaps}_{\del,h}(\xb)$ with possibly smaller constants.

\vskip5pt
{\em{\textbf{Thin obtuse triangles and hanging nodes.}}}
{ 
In the one and ten triangle cases of Fig. \ref{discap}-left, we see that thin obtuse triangles are used to approximate the cap. This can also occur for the approximate ball $B^{nocaps}_{\del,h}(\xb)$; see Fig. \ref{fappoxball}a. In Fig. \ref{discap}-right, we see that the two triangle case results in a ``hanging node'' as would also occur for the ten triangle case, where by ``hanging node'' we mean that a vertex of a triangle is not also a vertex of an abutting triangle. Both thin obtuse triangles and hanging nodes are considered to be anathemas for finite element discretizations. However, {\em here, we use the triangulation of approximate balls only to define composite quadrature rules for the inner integrals; they are not used to define finite element discretizations.} The latter are always effected using only finite element triangles, i.e., the triangles in the set $\mcT_h$.
}

\subsubsection{\textbf{Whole-triangle ball approximation based on barycenter location - Figure \ref{fappoxball}c}}\label{bary}

In this section we consider an approximate ball $B^{barycenter}_{\del,h}(\xb)$ that, for any point $\xb$, can be constructed without having to deal with caps nor with intersections of the ball boundary and element edges. In fact, the recipe for constructing this type of approximate ball is simply 
\beq{baryball}
\bal
    B^{barycenter}_{\del,h}(\xb) &= \big\{ \cup_{k=1}^K \,{\mathcal E}_k \quad\mbox{such that the barycenter of ${\mathcal E}_k\in B_{\del}(\xb)$} \big\}
    \\
    &= \big\{ \cup_{k=1}^K \,{\mathcal E}_k \quad\mbox{such that $|\xb-\xb^{barycenter}_k|\le\del$} \big\},
\eal
\eeq
where $\xb^{barycenter}_k$ denotes the barycenter of the finite element ${\mathcal E}_k$. Thus all elements whose barycenters are in the ball $B_{\del}(\xb)$ are part of the approximate ball but those whose barycenters are outside the ball are not. An illustration of the approximate ball $B^{barycenter}_{\del,h}(\xb)$ is given in Fig. \ref{fappoxball}c.  Unlike the approximate balls discussed in Sections \ref{poly} and \ref{acaps}, the approximate ball of \eqref{baryball} includes areas outside the ball $B_{\del}(\xb)$ and leaves out areas inside that ball.

The approximate ball $B^{barycenter}_{\del,h}(\xb)$ is subdivided into a set ${\mcT}^{barycenter}_{\del,h,\xb}$ of {\em whole} finite element triangles, i.e., ${\mcT}^{barycenter}_{\del,h,\xb}\subset{\mcT}_h$. The membership of ${\mcT}^{barycenter}_{\del,h,\xb}$ depends on the horizon $\del$, the grid size $h$, and the position of the center of the exact ball $B_{\del}(\xb)$.

\vskip5pt
{\em{\textbf{Geometric error.}}}
It is obvious that as $h\to0$ the approximate ball $B^{barycenter}_{\del,h}(\xb)$ reduces to the ball  $B_{\del}(\xb)$ and certainly the area of the former converges to the area of the latter. It is also easy to prove that the convergence is at least linear in $h$ because each partial triangle included or left out has an area of ${\mathcal O}(h^2)$ and, similarly to what we saw in Sec. \ref{poly}, the number of such partial triangles is of ${\mathcal O}(\frac1h)$. Thus, we have that 
\beq{bdifare}
|\Delta B^{barycenter}_{\del,h}(\xb)|={\mathcal O}(h).
\eeq
This estimate also holds for the case $\del=$ constant $\,\times\,h$. 

\vskip5pt
{\em{\textbf{Lack of sharpness of the estimate \eqref{bdifare}.}}}
The estimate \eqref{bdifare} may not be sharp because it does not take into account the ``cancellation'' of areas, i.e., that some of the whole triangles in $B^{barycenter}_{\del,h}(\xb)$ add area to the ball $B_\del(\xb)$ (see the pink triangles in Fig. \ref{ftriangletypes}) whereas some of the triangles that intersect $B_\del(\xb)$ are left out of $B^{barycenter}_{\del,h}(\xb)$ and thus subtract area (see the magenta triangles in Fig. \ref{ftriangletypes}). Thus, we conjecture that the cancellation due to areas added and areas subtracted might result in
\beq{bdifaree}
|\Delta B^{barycenter}_{\del,h}(\xb)|={\mathcal O}(h^\acancel)
\qquad\mbox{with $\acancel>1$}
\eeq
and possibly $\alpha\approx2$. This would occur if the difference in the area inside of the ball that is not included and that of area outside the ball that is included is of ${\mathcal O}(h^\acancel)$. This second conjecture seems to be reasonable, at least for locally quasi-uniform grids. Support for the veracity of these conjectures is provided by numerical results given in Sec. \ref{numerics} in which further discussions about the conjectures are also given. 

\vskip5pt
{\em{\textbf{Solution error.}}}
According to Corollary \ref{prop:error}, and if the kernel is integrable and translationally invariant or just square integrable, we have, at least conjecturally, that
\beq{errorbary}
\|u_h-\wuh\|_{L^2(\OuO)}\leq C |||u_h-\wuh ||| \sim O(h^\acancel)
\eeq
with $\acancel>1$ and possibly $\alpha\approx2$. 

\subsubsection{\textbf{Whole-triangle ball approximation based on overlap with ball - Figure \ref{fappoxball}d}}\label{overlapball}

In this section, we consider another approximate ball $B^{overlap}_{\del,h}(\xb)$ that, for any point $\xb$, can be constructed without having to deal with caps nor with intersections of the ball boundary and element edges nor with the location of triangle barycenters. The recipe for constructing this type of approximate ball is even simpler than that for $B^{barycenter}_{\del,h}(\xb)$; it is given by 
\beq{overlap}
    B^{overlap}_{\del,h}(\xb) = \big\{ \cup_{k=1}^K \,{\mathcal E}_k \quad\mbox{such that ${\mathcal E}_k\cap B_{\del}(\xb)\ne\emptyset\big\}$},
\eeq
i.e., all elements that overlap with the ball $B_{\del}(\xb)$ are part of the approximate ball, and those that do not overlap are not. An illustration of the approximate ball $B^{overlap}_{\del,h}(\xb)$ is given in Fig. \ref{fappoxball}d. Unlike the approximate balls discussed in Sections \ref{poly} and \ref{acaps}, the approximate ball of \eqref{overlap} includes areas outside the ball $B_{\del}(\xb)$ but unlike the ball discussed in Sec. \ref{bary}, the ball of \eqref{overlap} covers the ball $B_{\del}(\xb)$. 

The approximate ball $B^{overlap}_{\del,h}(\xb)$ is subdivided into a set ${\mcT}^{overlap}_{\del,h,\xb}$ of {\it whole} finite element triangles triangles, i.e., ${\mcT}^{overlap}_{\del,h,\xb}\subset{\mcT}_h$. The membership of ${\mcT}^{overlap}_{\del,h,\xb}$ depends on the horizon $\del$, the grid size $h$, and the position of the center of the exact ball $B_{\del}(\xb)$.

\vskip5pt
{\em{\textbf{Geometric error.}}}
It is obvious that as $h\to0$ the approximate ball $B^{overlap}_{\del,h}(\xb)$ reduces to the ball  $B_{\del}(\xb)$ and certainly the area of the former converges to the area of the latter. It is also easy to prove, as it is for the ball $B^{barycenter}_{\del,h}(\xb)$, that the convergence is linear in $h$. Thus, we have that
\beq{odifare}
|\Delta B^{overlap}_{\del,h}(\xb)|={{\mathcal O}(h)}.
\eeq
However, unlike the case of Sec. \ref{bary}, for $B^{overlap}_{\del,h}(\xb)$, there is no possibility of the convergence rate of $|\Delta B_{\del,h}^{overlap}(\xb)|$ being better than one because there is no opportunity for the cancellation of areas.
\vskip5pt
{\em{\textbf{Solution error.}}}
According to Corollary \ref{prop:error}, and, if the kernel is integrable and translationally invariant or just square integrable, we have that 
\beq{erroroverlap}
\|u_h-\wuh\|_{L^2(\OuO)}\leq C |||u_h-\wuh ||| \sim O(h)
\eeq
This estimate is sharp, as is illustrated by the numerical results in Sec. \ref{numerics}.

\subsection{\textbf{Shifted center approximate ball}}\label{subsec:shifted}

\begin{figure}[ht]
\begin{center}
\includegraphics[height=1.5in]{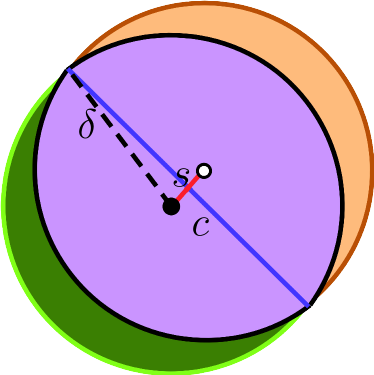}
\end{center}
\caption{The green-violet ball $B_\delta(\xb)$ is centered at a point $\xb$ (the filled dot). The orange-violet ball $B_\delta(\xb^{shifted})$ is shifted so that it is centered at another point $\xb^{shifted}$ (the open dot). Each half of the violet area is a circular cap having cord length $c$ and radius $\delta$ of the green or orange balls; $s$ is the separation distance between the centers $\xb$ and $\xb^{shifted}$.}
\label{fshifted}
\end{figure}

The polygonal approximate balls constructed in Sec. \ref{polyballs} share the same center as that of the exact ball $B_\del(\xb)$ but differ in their shape. Here, we consider an approximate ball $B^{shifted}_{\del,h}(\xb)$ that differs from the exact ball only in the position of their centers. For example, in Fig. \ref{fshifted}, the exact ball is centered at the filled dot and is depicted by the green and violet areas, whereas the shifted ball is centered at the open dot and is depicted by the orange and violet areas.  Specifically, when we use shifted balls, we shift the center $\xb$ of the ball to a new point $\xb^{shifted}$ in such a way that $s=|\xb - \xb^{shifted}|=\mcO(h)$. In particular, in our experiments we choose the barycenter of the triangle for the center of the shifted ball.

\vskip5pt
{\em{\textbf{Geometric error.}}}
It is obvious that as $h\to0$ the approximate ball $B^{shifted}_{\del,h}(\xb)$ reduces to the ball  $B_{\del}(\xb)$ and, of course, the area of the former is the same as the area of the latter. Thus, here, the geometric error is solely due to the shift of the center.

Referring to Fig. \ref{fshifted}, we estimate the areas of the two lunes (the green and orange areas) by subtracting the area of violet region from the area of the ball. Note that each half of the violet region is a circular cap for one of the balls; those caps are defined by the radius $\delta$ of the ball (the dashed line segment), the cord length $c$ (the blue line segment), and $\frac{s}2$, where $s$ denotes the separation distance between the two centers of the balls (the red line segment). We have that $s={\mathcal O}(h)$ so that
$$
\bal
\delta^2 = \big(\frac{c}{2}\big)^2 + \big(\frac{s}{2}\big)^2
\quad&\Rightarrow\quad
\sqrt{1 - \big(\frac{c}{2\del}\big)^2} = \frac{s}{2\del}
\\
\frac{c}{2\del} = \sqrt{1 - \big(\frac{s}{2\del}\big)^2} \approx 1 - \frac12 \big(\frac{s}{2\del}\big)^2
\quad&\Rightarrow\quad
\frac{c}{2\del}\sqrt{1 - \big(\frac{c}{2\del}\big)^2} \approx \frac{s}{2\del}
\\
\arcsin \Big(\frac{c}{2\delta}\Big) \approx \arcsin \Big( 1 - \frac12 \big(\frac{s}{2\del}\big)^2  \Big)
&\approx \frac{\pi}2 - \frac{s}{2\delta},
\eal
$$
where here the symbol $\approx$ means that terms of ${\mathcal O}(\frac{s^3}{\del^3})$ have been neglected. Then, from \eqref{areacap}, we have that
$$
|\mbox{violet regions}| = 2|\mbox{circular cap}| \approx \pi\del^2 - 2s\del
$$
so that
$$
\bal
|\mbox{green lune}| = |\mbox{orange lune}| &= |B_\del^{shifted}(\xb)| - |\mbox{violet region}|
\\&\approx \pi\del^2 - (\pi\del^2 - 2s\del) =  2s\del = \del{\mathcal O}(h).
\eal
$$
This implies that, for the shifted ball approximation, 
$$
|\Delta B_{\del,h}^{shifted}(\xb)|=\del\mcO(h).
$$
However, as was the case for the approximate ball of Sec. \ref{bary}, numerical evidence given in Sec. \ref{numerics} indicates that this estimate may not be sharp. A possible explanation for the better observed rate of convergence is that again a cancellation effect comes into play due to the symmetric placement of quadrature points with respect to the barycenter.

\vskip5pt
{\em{\textbf{Solution error.}}} 
According to Corollary \ref{prop:error}, and, if the kernel is integrable and translationally invariant or just square integrable, we have that 
\begin{align*}
  \|u_h-\wuh\|_{L^2(\OuO)}\leq C |||u_h-\wuh |||  &\sim \mcO(h) \quad~\text{if\,\,\, $h\ll\del$} \quad \mbox{and} \\
  \|u_h-\wuh\|_{L^2(\OuO)}\leq C |||u_h-\wuh |||  &\sim \mcO(h^2) \quad \text{if\,\,\, $\del=\,{\rm constant}\,\times\,h$}.
\end{align*}
Again, this estimate may not be sharp.

\vskip5pt
{\em{\textbf{Pairing with other approximate balls.}}}
This shifted-center approximation can be paired with any of the four approximate balls considered in Sec. \ref{polyballs} in which case one is approximating both the position of the center of the ball and the ball shape. 

\section{Approximating inner integrals}\label{apinner}

We consider three approaches for the approximation of the inner integrals appearing in  \eqref{bform3} and \eqref{lfunc3}  or \eqref{bform3h} and \eqref{lfunc3h}. In  Sec. \ref{gqr} we consider {\em global} quadrature rules for the exact ball $B_\del(\xb)$. We then consider {\em composite} quadrature rules, in Sec. \ref{ecaps} for the exact ball and then in Sec. \ref{cqr} for approximate balls $B_{\del,h}(\xb)$.

\subsection{Global quadrature rules for balls}\label{gqr}

We consider {\em global} quadrature rules over the exact ball $B_\del(\xb)$ so that the approximate balls of Sec. \ref{approxballs} do not come into play. Thus, considering the inner integrals in \eqref{bform3} and \eqref{lfunc3}, the task at hand is to effect the approximation 

\beq{gqr1}
\bal
    &   \int_{B_\del(\xb)\cap\Omega} \big(\phi_j(\yb)-\phi_j(\xb)\big)\big(\phi_{j'}(\yb)-\phi_{j'}(\xb)\big) \psi(\xb,\yb) d\yb 
    \\&  \quad \approx
       \sum_{q=1}^{Q_{global}} w^{global}_q \big(\phi_j(\yb^{global}_{q})-\phi_j(\xb)\big)\big(\phi_{j'}(\yb^{global}_{q})-\phi_{j'}(\xb)\big)\psi(\xb,\yb^{global}_{q})  
\eal
\eeq
and similar terms appearing in \eqref{bform3} and \eqref{lfunc3}, where $\{w^{global}_q,\yb^{global}_q\}_{q=1}^{Q_{global}}$ denotes a set of quadrature weights and points. Such rules are given in, e.g., Ref. \refcite{AandS}.

\begin{figure}[ht]
\centering
\begin{tabular}{ccc}
\includegraphics[height=1.35in]{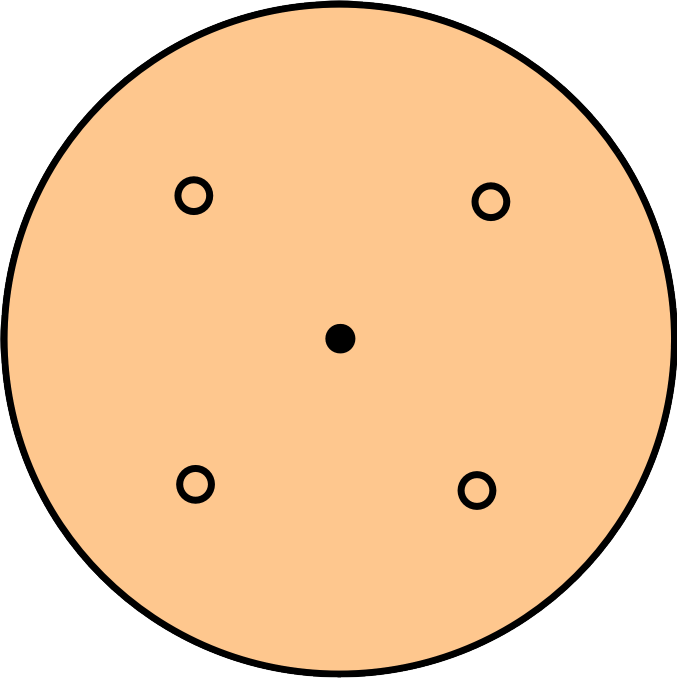} &
\includegraphics[height=1.35in]{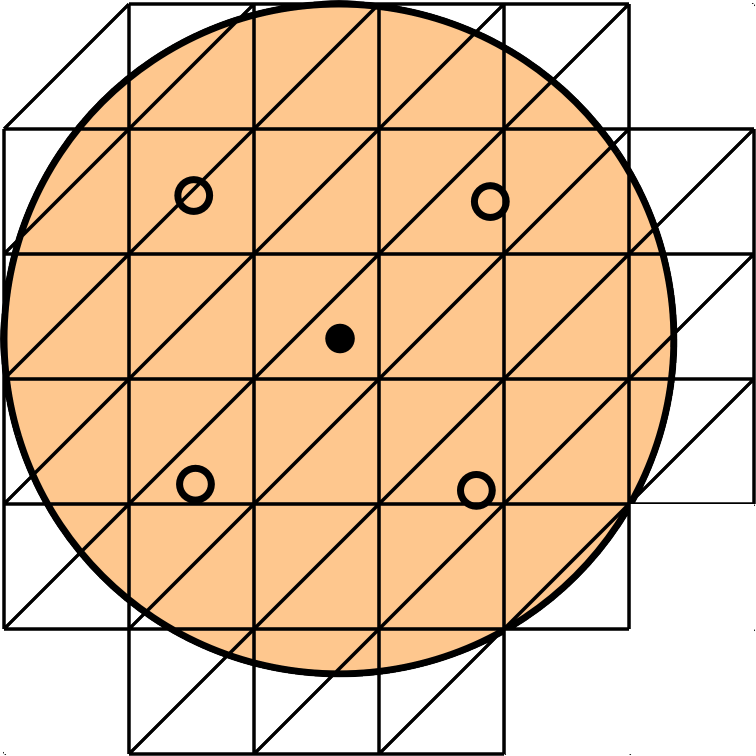} &
\includegraphics[height=1.35in]{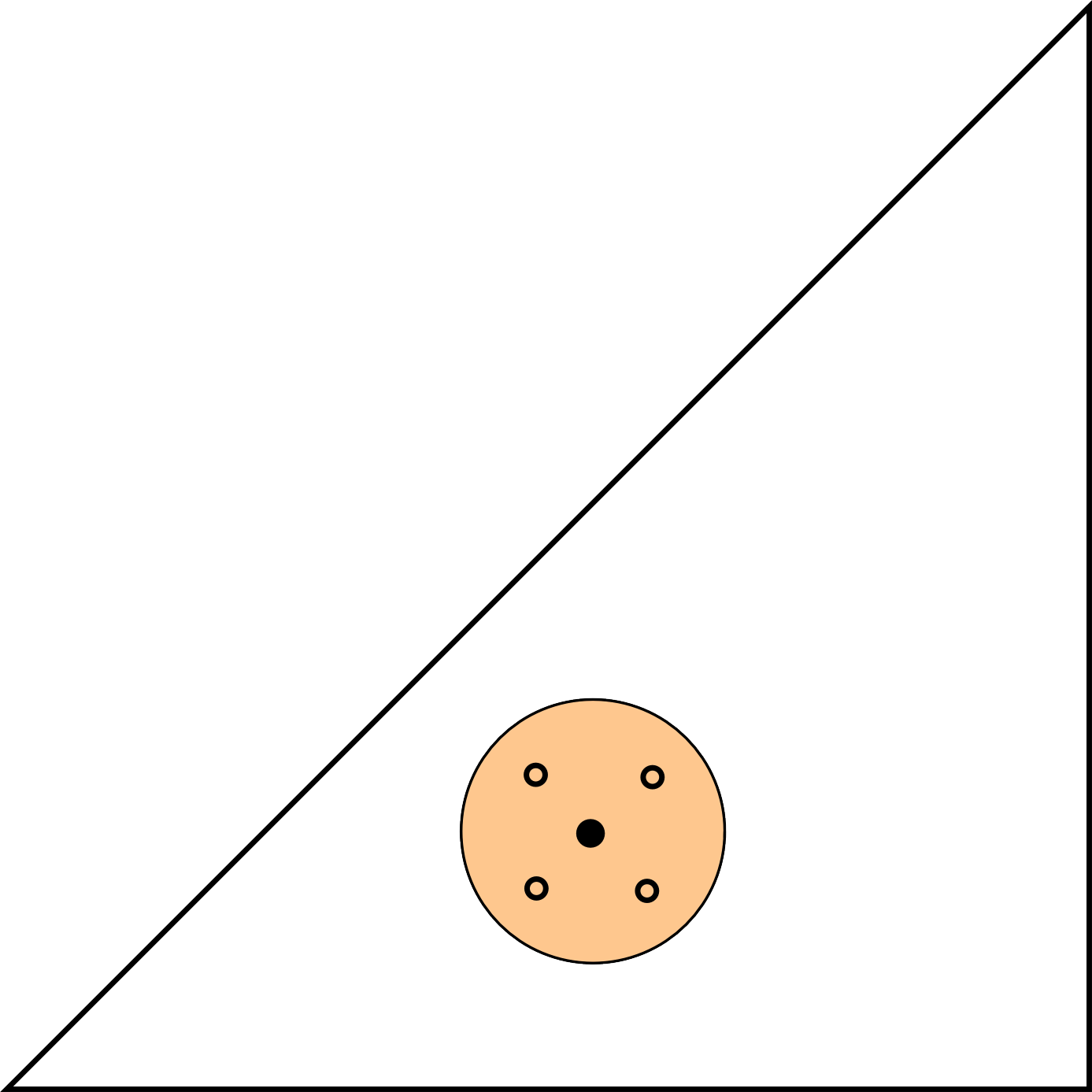} 
\\
(a) & (b) & (c) 
\end{tabular}
\caption{(a) The orange disc depicts the ball $B_\del(\xb)$ centered at $\xb$ and having radius $\del>h$. The open circles depict the quadrature points of a global quadrature rule that can be used to approximate the integral in \eqref{gqr1}. 
(b) The triangular grid is a portion of the finite element triangulation $\mcT_h$. 
(c) The situation in which the ball $B_\del(\xb)$ centered at a point $\xb$ in the interior of the finite element $\mcE_k$ has radius $\del$ that is sufficiently small relative to the grid size $h$ so that the whole ball is contained within $\mcE_k$. 
}
\label{globalquad}
\end{figure}

A main advantage accruing from using a global quadrature rule is that one does not have to deal with triangles when one approximates the inner integral, one simply integrates over the ball, as is implied by Fig. \ref{globalquad}a. As already mentioned, a second advantage is that there is no need to approximate the ball so that no geometric error is incurred. However, in the setting in which $\del$ is fixed and $\del>h$ (that is of most interest to us) there are two serious disadvantage stemming from using a global quadrature rule that outweighs these advantages, so that we do not pursue the use of such rules beyond what is written in this section. First, the integrand in \eqref{gqr1} involves piecewise-polynomial functions defined with respect to the finite element grid; see Fig. \ref{globalquad}b. Most commonly, these functions are continuous but are not continuously differentiable. Such functions are not sufficiently smooth to take advantage of the accuracy potential of even low-precision global rules.
The second disadvantage is that the error incurred by the use of a global quadrature rule depends on $\del$, so that if $h\ll\del<1$, one would need a very high-order quadrature rule to balance the quadrature error with the other errors incurred which, if piecewise-linear finite element spaces are used, are of $\mcO(h^2)$. However, as we just commented, the use of high-order quadrature rules is compromised due to the lack of smoothness of the integrand so that, in the end, one cannot balance the $\del$ with the $h$ errors. We just mention that there is a third disadvantage in that for the term involving the data $g(\xb)$ in \eqref{lfunc3}, the domain of integration is a partial ball.

There is the situation illustrated in Fig. \ref{globalquad}c for which the use of a global quadrature rule on balls may be applicable, namely $\del$ being sufficiently small compared to $h$. We note that the setting of $\del$ small compared to $h$ arises relatively rarely in applications, but is useful for illustrating that a nonlocal model reduces to a local one as the horizon $\del\to0$. Here, the ball center at $\xb$ would have to lie in the interior of an element $\mcE_k$. Furthermore, the radius $\del$ of the ball $B_\del(\xb)$ would have to be sufficiently small (relative to the diameter of the element and the position of the point $\xb$) so that the whole ball is contained within the element $\mcE_k$. In this situation, the domain of integration in \eqref{gqr1} does not straddle across triangle boundaries so that the integrand is smooth. Note that in this case, the error in the quadrature rule depends on $h$ and not on $\del$. As a result, a relatively low-precision quadrature rule can be chosen in \eqref{gqr1} so that the quadrature error is commensurate with other $h$-dependent errors incurred, e.g., due to finite element approximation. However, there is a complication in handling the inner integral in \eqref{bform3} and \eqref{lfunc3} over the domain $\Omega_{\mathcal I}\cap B_\del(\xb)$. Necessarily, that domain is always a partial ball so that one would need to use a global integration rule that can handle arbitrary partial balls that are created by cutting off part of a ball by a cord. {Such rules do exist;  see Ref. \refcite{sectorrule}. Because their integrands only involve given data, this complication may not add significantly to the cost of the assembly process.}

\subsection{\textbf{Composite quadrature rules for exact balls}}\label{ecaps}

In this section we consider composite quadrature rules for the whole ball so that there is no error incurred due to geometric approximation; errors are due only to the use of quadrature rules. As illustrated in Fig. \ref{fappoxball}b, we subdivide the ball into the polygon of Fig. \ref{fappoxball}a (the orange region) and the circular caps (the green regions). Specifically, let $\mcT_{\del,h,\xb}^{exactcaps}$ denote the set of caps and recall that $\mcT_{\del,h,\xb}^{nocaps}$ denotes the set of triangles in the approximate ball of Sec. \ref{poly}. Letting $\mcT_{\del,h,\xb}^{exact}=\mcT_{\del,h,\xb}^{nocaps}\cup\mcT_{\del,h,\xb}^{exactcaps}$, we have that 
$$
B_\del(\xb) = \mcT_{\del,h,\xb}^{exact}=\mcT_{\del,h,\xb}^{nocaps}\cup\mcT_{\del,h,\xb}^{exactcaps}=
\big(\cup_{\widetilde\mcE_{k'}\in\mcT_{\del,h,\xb}^{nocaps}} {\widetilde\mcE}_{k'}\big) \cup \big(\cup_{{\widehat\mcE}_{k'}\in\mcT_{\del,h,\xb}^{exactcaps}} {\widehat\mcE}_{k'}\big),
$$
where $\widetilde\mcE_{k'}$ denotes a typical triangle in $\mcT_{\del,h,\xb}^{nocaps}$ and ${\widehat\mcE}_{k'}$ denotes a typical cap in $\mcT_{\del,h,\xb}^{exactcaps}$. Then, considering \eqref{bform3} and \eqref{lfunc3}, the task at hand is to effect the approximation

\beq{cqr3}
\bal
    &
     \int_{B_\del(\xb)\cap\Omega} \big(\phi_j(\yb)-\phi_j(\xb)\big)\big(\phi_{j'}(\yb)-\phi_{j'}(\xb)\big) \psi(\xb,\yb) d\yb
  \\&
   \quad =\sum_{\widetilde\mcE_{k'}\in \mcT_{\del,h,\xb}^{nocaps}} \int_{\widetilde\mcE_{k'} \cap \Omega} \big(\phi_j(\yb)-\phi_j(\xb)\big)\big(\phi_{j'}(\yb)-\phi_{j'}(\xb)\big) \psi(\xb,\yb) d\yb
 \\&   \qquad\qquad+
   \sum_{{\widehat\mcE}_{k'}\in \mcT_{\del,h,\xb}^{exactcaps}} \int_{{\widehat\mcE}_{k'}\cap \Omega} \big(\phi_j(\yb)-\phi_j(\xb)\big)\big(\phi_{j'}(\yb)-\phi_{j'}(\xb)\big) \psi(\xb,\yb) d\yb
    \\&  \quad\approx
           \sum_{\widetilde\mcE_{k'}\in \mcT_{\Omega,\del,h,\xb}^{nocaps}}
      \sum_{q=1}^{Q_{nocaps}} w^{nocaps}_{q,k'} \big(\phi_j(\yb^{nocaps}_{q,k'})-\phi_j(\xb)\big)
\\&
\qquad\underbrace{\phantom{\qquad\qquad\qquad}      
     \times \big(\phi_{j'}(\yb^{nocaps}_{q,k'})-\phi_{j'}(\xb)\big)\psi(\xb,\yb^{nocaps}_{q,k'}) }_{\mbox{\em composite quadrature rule over the triangles in $\mcT_{\del,h,\xb}^{nocaps}$}}
     \\&  \qquad\qquad +      
     \sum_{{\widehat\mcE}_{k'}\in \mcT_{\Omega,\del,h,\xb}^{exactcaps}}
     \sum_{q=1}^{Q_{caps}} w^{exactcaps}_{q,k'} \big(\phi_j(\yb^{exactcaps}_{q,k'})-\phi_j(\xb)\big)
     \\&
 \qquad\qquad    \underbrace{
 \phantom{\qquad\qquad\qquad}
 \times\big(\phi_{j'}(\yb^{exactcaps}_{q,k'})-\phi_{j'}(\xb)\big)\psi(\xb,\yb^{exactcaps}_{q,k'})\}}_{\mbox{\em composite quadrature rule over the caps in $\mcT_{\del,h,\xb}^{exactcaps}$}}
\eal
\eeq
and similar terms appearing in \eqref{bform3h} and \eqref{lfunc3h}. Here, $\{w^{nocaps}_{q,k'},\yb^{nocaps}_{q,k'}\}_{q=1}^{Q_{nocaps}}$ denotes a set of quadrature weights and points for the composite rule for the polygon of  Sec. \ref{poly}. 

Here, no geometric error is incurred because we are using a whole ball. We suppose that the outer integral is integrated exactly. Then,  for  piecewise-linear  basis  functions  and assuming that the kernel function $\psi(\xb,\yb)$ is constant, we have that the integrand is a polynomial of degree two in the components of $\yb$. Thus, we would use a precision two quadrature rule so that the quadrature error is commensurate with the $\mcO(h^2)$ error incurred by the finite element discretization. For this purpose, we can use a three-point symmetric Gaussian quadrature rule for triangles. We expect these rules to also work equally well for smooth non-constant kernel functions.

In \eqref{cqr3}, $\{w^{exactcaps}_{q,k'},\yb^{exactcaps}_{q,k'}\}_{q=1}^{Q_{caps}}$ denotes a set of quadrature weights and points for the composite quadrature rule for the set of caps $\mcT_{\del,h,\xb}^{exactcaps}$. The error incurred when using piecewise linear finite element basis functions is of $\mcO(h^2)$. To render the error incurred by the quadrature rule for caps to also be of $\mcO(h^2)$, a one-point centroid rule would more than suffice. Referring to Fig. \ref{caparea}, that point is located along the bisector of the circular sector at a distance $4\del\sin^3\theta/3(2\theta-\sin 2\theta)$  from the center of the ball. The quadrature weight is the area of the cap which is given by $\frac14 \del^2  (2  \theta - \sin(2  \theta))$. If a higher-order finite element approximation is used, then the quadrature rule used for the caps has to be commensurately higher-order as well. A family of such rules is given in Ref. \refcite{sectorrule}. 

\subsection{\textbf{Composite quadrature rules for polytopial approximations of balls}}\label{cqr}

For $\sharp\in\{nocaps,\,approxcaps,\,barycenter,\,overlap\}$, in Sections \ref{poly}--\ref{overlapball} we have the  approximate balls $B_{\del,h}^{\sharp}(\xb)$, each of which is covered by a set $\mcT_{\del,h,\xb}^{\sharp}$ of disjoint triangles. We consider composite quadrature rules over those approximate balls. Thus, considering \eqref{bform3h} and \eqref{lfunc3h}, the task at hand is to effect the approximation
\beq{cqr1} 
\bal
&    \int_{B_{\del,h}^{\sharp}(\xb)\cap\Omega} \big(\phi_j(\yb)-\phi_j(\xb)\big)\big(\phi_{j'}(\yb)-\phi_{j'}(\xb)\big) \psi(\xb,\yb) d\yb
\\&\qquad    =\sum_{\widetilde\mcE_{k'}\in\mcT_{\del,h,\xb}^{\sharp}}
    \int_{\widetilde\mcE_{k'} \cap \Omega} \big(\phi_j(\yb)-\phi_j(\xb)\big)\big(\phi_{j'}(\yb)-\phi_{j'}(\xb)\big) \psi(\xb,\yb) d\yb
    \\&\qquad   \approx
      \sum_{\widetilde\mcE_{k'}\in \mcT_{\Omega,\del,h,\xb}^{\sharp}}
     \sum_{q=1}^{Q_{\sharp}} w^{\sharp}_{q,k'} \big(\phi_j(\yb^{\sharp}_{q,k'})-\phi_j(\xb)\big)\big(\phi_{j'}(\yb^{\sharp}_{q,k'})-\phi_{j'}(\xb)\big)\psi(\xb,\yb^{\sharp}_{q,k'}) ,
\eal
\eeq
where, for each member $\widetilde\mcE_{k'}$ of the set of triangles in $\mcT_{\del,h,\xb}^{\sharp}$, we use a quadrature rule with weights and points $\{w^{\sharp}_{q,k'},\yb^{\sharp}_{q,k'}\}_{q=1}^{Q_{\sharp}}$. Because the subdomains within each of the four approximate balls are all triangles, one can use the same quadrature rule for all triangles within the approximate ball. 

\subsubsection{\textbf{Error-commensurate and heuristics choices of quadrature rules}}\label{quade_inner}

We discuss two ``philosophies'' for choosing quadrature rules for inner integrals. Because the geometric error incurred by the use of approximate balls is of $\mcO(h^2)$ at best, we restrict our discussion to piecewise-linear finite element approximations, for which the rate of convergence is also $\mcO(h^2)$ at best.

\vskip5pt
{
{\em{\textbf{Error-commensurate choices of quadrature rules.}}}
In Section \ref{sec:weak}, four sources of errors were listed, including one due to the use of quadrature-rule approximations of inner integrals. The choice of what rule to use is, in principle, governed by the minimum precision needed to render the inner integral quadrature error commensurate with other errors incurred while at the same time using the fewest number of quadrature points needed to achieve that precision. Because the finite element approximation error is at best of $\mcO(h^2)$, it seems that one should avoid rules that have higher accuracy than that. Even lower-accuracy rules seem appropriate if the geometric error is of $\mcO(h)$.

This philosophy results in the following choices of quadrature rules, where, for simplicity, we restrict our discussion to constant kernel functions $\psi(\xb,\yb)$.
\hangseven--~For $\sharp\in\{nocaps,approxcaps\}$, the geometric error is of $\mcO(h^2)$ so that any precision one rule can be used, i.e., any rule that integrates quadratic polynomials exactly can be used.

\hangseven--~For $\sharp= overlap$, the geometric error is of $\mcO(h)$ so that even though the finite element error is of $\mcO(h^2)$, the overall error cannot be better than $\mcO(h)$. Thus, in principle, a precision zero rule, i.e., one that integrates constants can be used.

\hangseven--~For $\sharp\in \{barycenter,shifted\}$, the geometric error is provably of $\mcO(h)$ so that a precision zero rule is seemingly called for. However, numerical results given in Section  \ref{numerics} indicate that the geometric errors for these two balls may be better than that, so that a precision one rule may be a better choice.

\noindent Another approach for choosing quadrature rules is discussed below. In that context, the precision of the quadrature rules suggested above should be viewed as what is minimally required to not ruin the accuracy achieved by finite element and geometric approximations. 

}

\vskip5pt
{\em{\textbf{Heuristic-based choices of quadrature rules.}}}
When using finite element methods for second-order elliptic PDE problems with smooth coefficients, one chooses a quadrature rule such that $\nabla\phi_j(\xb) \cdot \nabla\phi_{j'}(\xb)$ is integrated exactly (see Refs. \refcite{brenner,ciarlet}), where here $\phi_j(\xb)$ denotes a finite element basis function. Thus, letting $\mcE$ denote a generic finite element triangle and letting $\{\xb_q,w_q\}_{q=1}^Q$ denote the points and weights of a quadrature rule over $\mcE$, it is required that
	\beq{outloc}
	\int_{\mcE} \nabla\phi_j(\xb) \cdot \nabla\phi_{j'}(\xb) d\xb
	= \sum_{q=1}^Q w_q  \nabla\phi_j(\xb_q) \cdot \nabla\phi_{j'}(\xb_q).
	\eeq
For piecewise-linear finite element approximations, the integrand $\nabla\phi_j(\xb) \cdot \nabla\phi_{j'}(\xb)$ is constant, so that a rule that integrates piecewise constants should suffice.
 
We use the same reasoning to {\em heuristically} decide about what precision is needed for quadrature rules in the nonlocal case. Following that reasoning, and assuming that the kernel function $\psi(\xb,\yb)$ is a constant, we then seek a quadrature rule that is exact for the inner integrals appearing in \eqref{bform3} and \eqref{lfunc3}. For piecewise-linear finite element approximations, the integrand is quadratic so that a precision two quadrature rule is needed for exact integration. In our computations, we choose to use the heuristic philosophy so that we use a three-point symmetric Gaussian quadrature rule for triangles; see Ref. \refcite{AandS}. We expect these rules to also work equally well for smooth non-constant kernel functions. The precision of the heuristic choice for the quadrature rule is higher than that of the commensurate rules discussed above. We choose to use the heuristic rule because we have empirically found that the additional cost of using the three-point rule instead of a one-point rule is dominated by other costs incurred during the assembly process and, in addition, the error due to quadrature is dominated by the other errors incurred so that the overall error is smaller than when using a one-point rule.

We have tacitly glossed over an important difference between finite element methods for local and nonlocal problems. {\em Because there are no derivatives involved in nonlocal models, for the same polynomial finite element space, the integrands for nonlocal models involve higher-degree polynomials and thus require higher-precision quadrature rules compared to local models.}

\section{Approximating outer integrals}\label{apouter}

Superficially, it would seem that making a good choice of a quadrature rule to approximate the outer integrals in  \eqref{bform3} and \eqref{lfunc3}  or \eqref{bform3h} and \eqref{lfunc3h} is one of the simpler decisions one has to make in the assembly process. After all, the outer integrals seem to be the same as the single integrals encountered in the PDE setting, i.e., both involve a sum of integrals over the finite elements. However, as we explain in this section, there are subtle issues that render the approximation of the outer integral in nonlocal models not as straightforward as it first seems. For simplicity, we again assume were are dealing with triangular finite elements and with piecewise-linear finite element approximations.

To investigate the approximation of outer integrals, we fix an outer integral triangle $\mcE_k$, $k=1,\ldots,K_\Omega$, and and inner integral triangle\footnote{For simplicity, we refer to $\mcE_{k'}$ as a ``triangle'' for all cases, even though for $\sharp=exactcaps$, some $\mcE_{k'}$ are exact caps.}  $\mcE_{k'}$, $k'=1,\ldots,K$, and consider the double integral
\begin{equation} \label{def:outer_task_approx}
\bal
&\int_{\mcE_k }   \mcK^\sharp_{k';j,j'}(\xb) d\xb  \quad\mbox{with integrand}
\\
&\qquad  \mcK^\sharp_{k';j,j'}(\xb) = 
\int_{\mcE_{k'} \cap  B^\sharp_{\del,h}(\xb)} \big(\phi_{j}(\xb)-\phi_{j}(\yb)\big)\big(\phi_{j'}(\xb)-\phi_{j'}(\yb)\big)\psi(\xb,\yb)d\yb,
\eal
\eeq
where $\sharp\in\{exactcaps,\,nocaps,\,approxcaps,\,barycenter,\,overlap,\,shifted\}$. We note that the evaluation of $\mcK^\sharp_{k';j,j'}(\xb)$ at a point $\xb \in \mcE_k$ requires the approximation of inner integrals as discussed in Sec. \ref{apinner}. In Fig. \ref{fig:outerq}a, we depict the three types of interactions between an outer integral triangle $\mcE_k$ (in blue) and an inner integral triangle $\mcE_{k'}$. In that figure, the orange regions depict the interaction region for $\mcE_k$. If $\mcE_{k'}$ is the yellow triangle, then there is no interaction and therefore the integrand $\mcK^\sharp_{k';j,j'}(\xb)=0$ for all $\xb\in\mcE_k$. Thus, we focus on the other two types of interactions illustrated by the green triangle, all of which overlaps with the orange interaction region for $\mcE_k$, and the violet triangles for which the overlap is only partial.

\begin{figure}[ht]
\begin{center}
\begin{tabular}{ccc}
\includegraphics[height=1.3in]{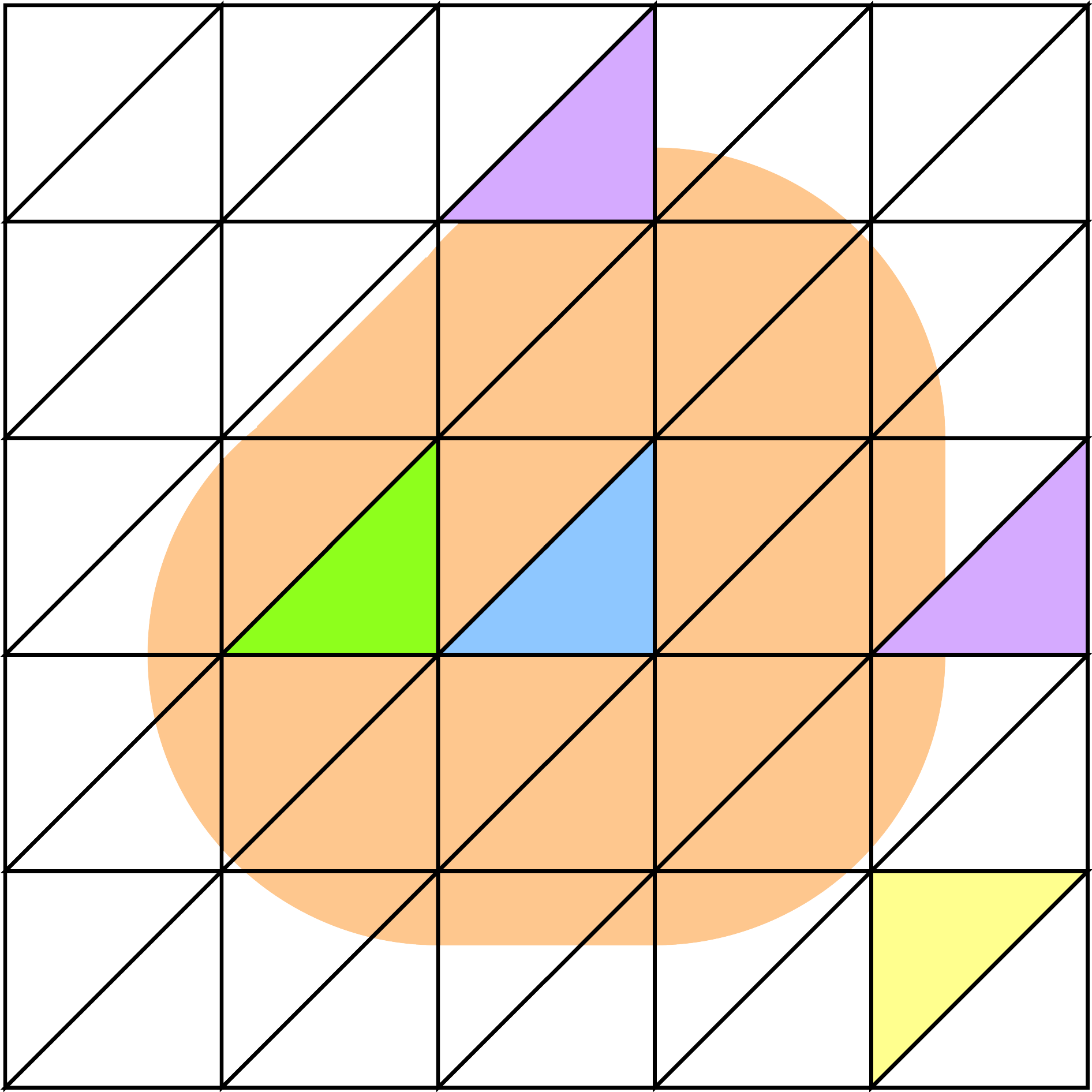}
&
\includegraphics[height=1.3in]{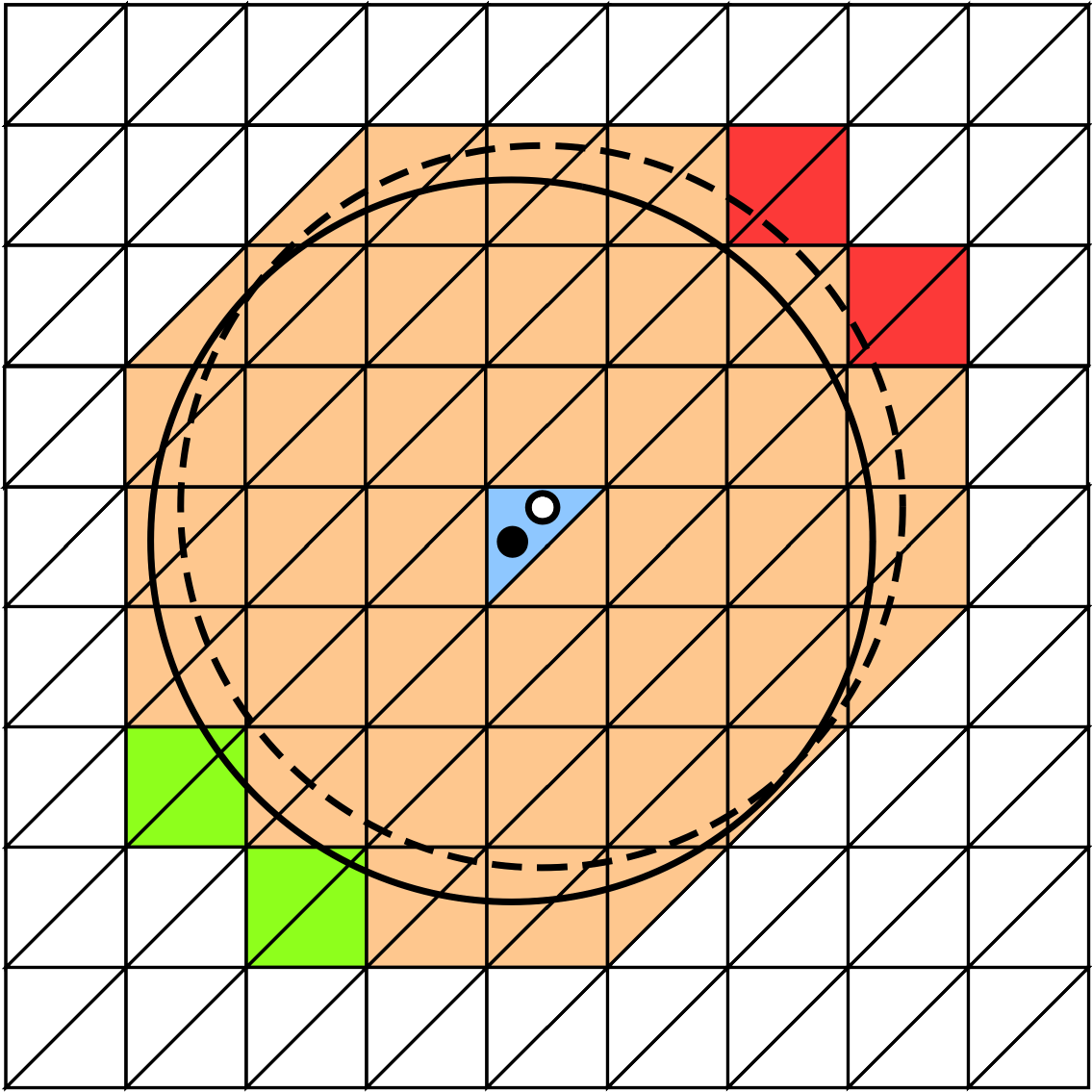}
&
\includegraphics[height=1.3in]{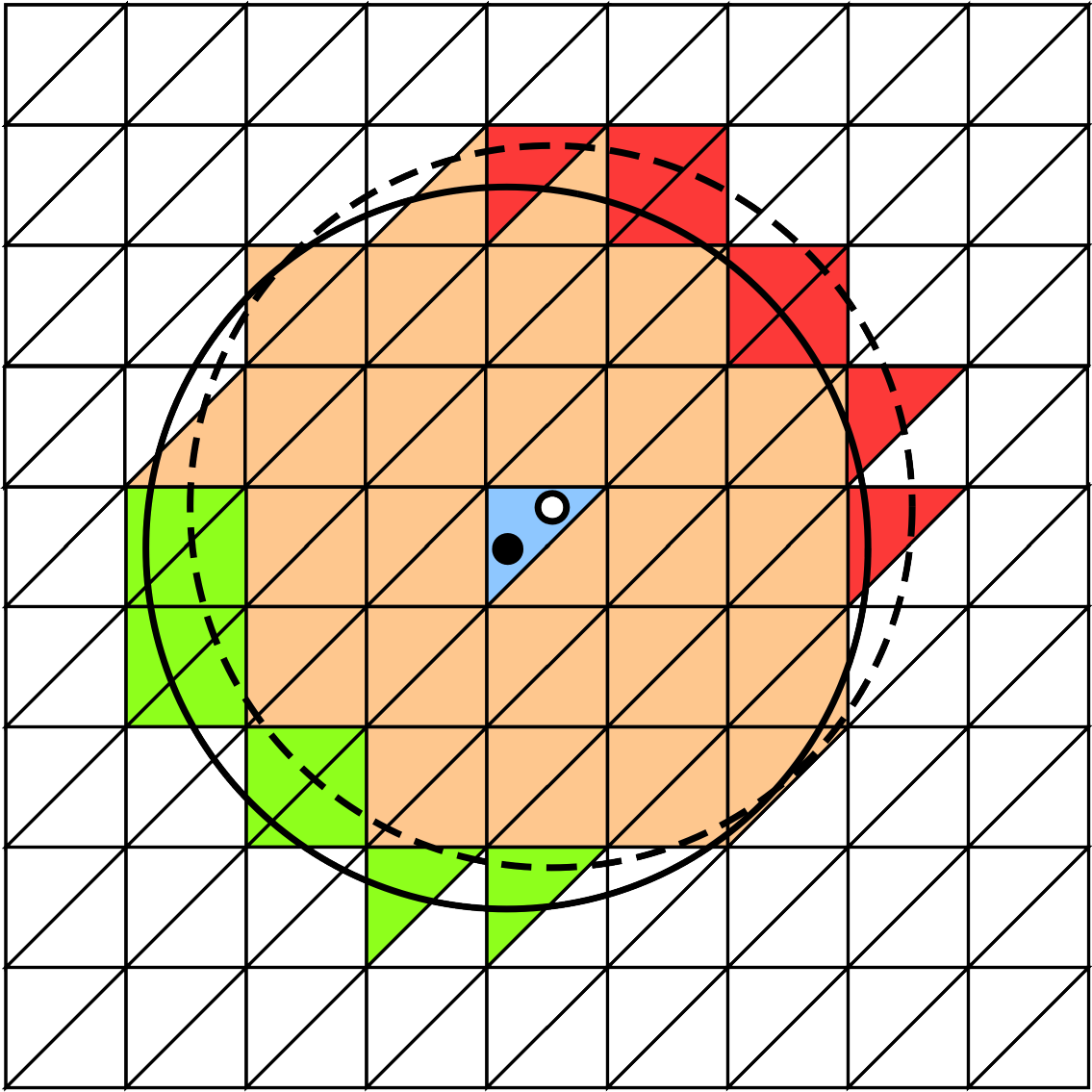}
\\
(a) & (b) & (c)
\\[1ex]
\includegraphics[height=1.3in]{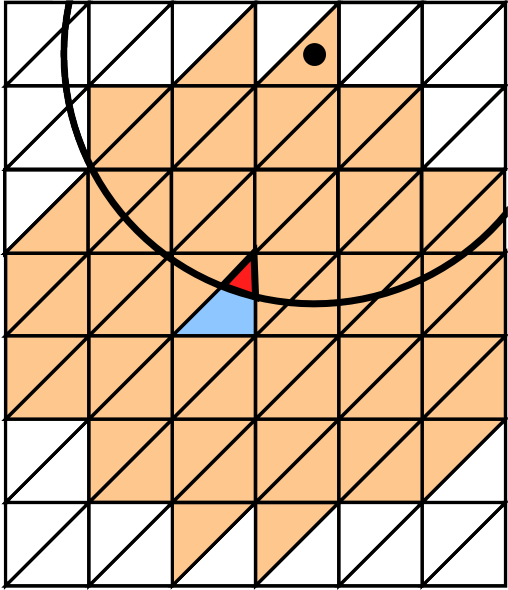}
&
\;\includegraphics[height=1.3in]{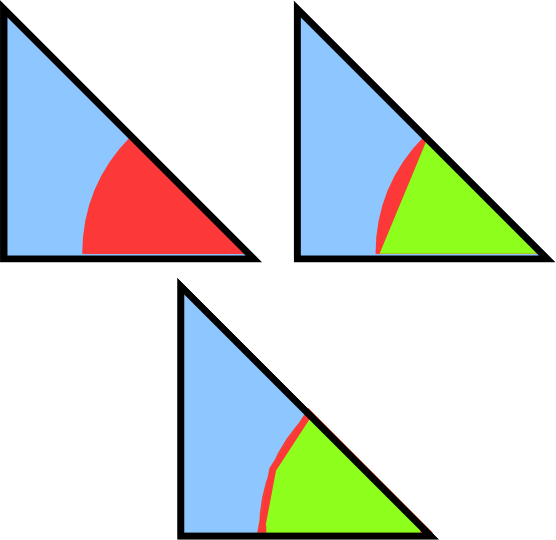}
&
\includegraphics[height=1.2in]{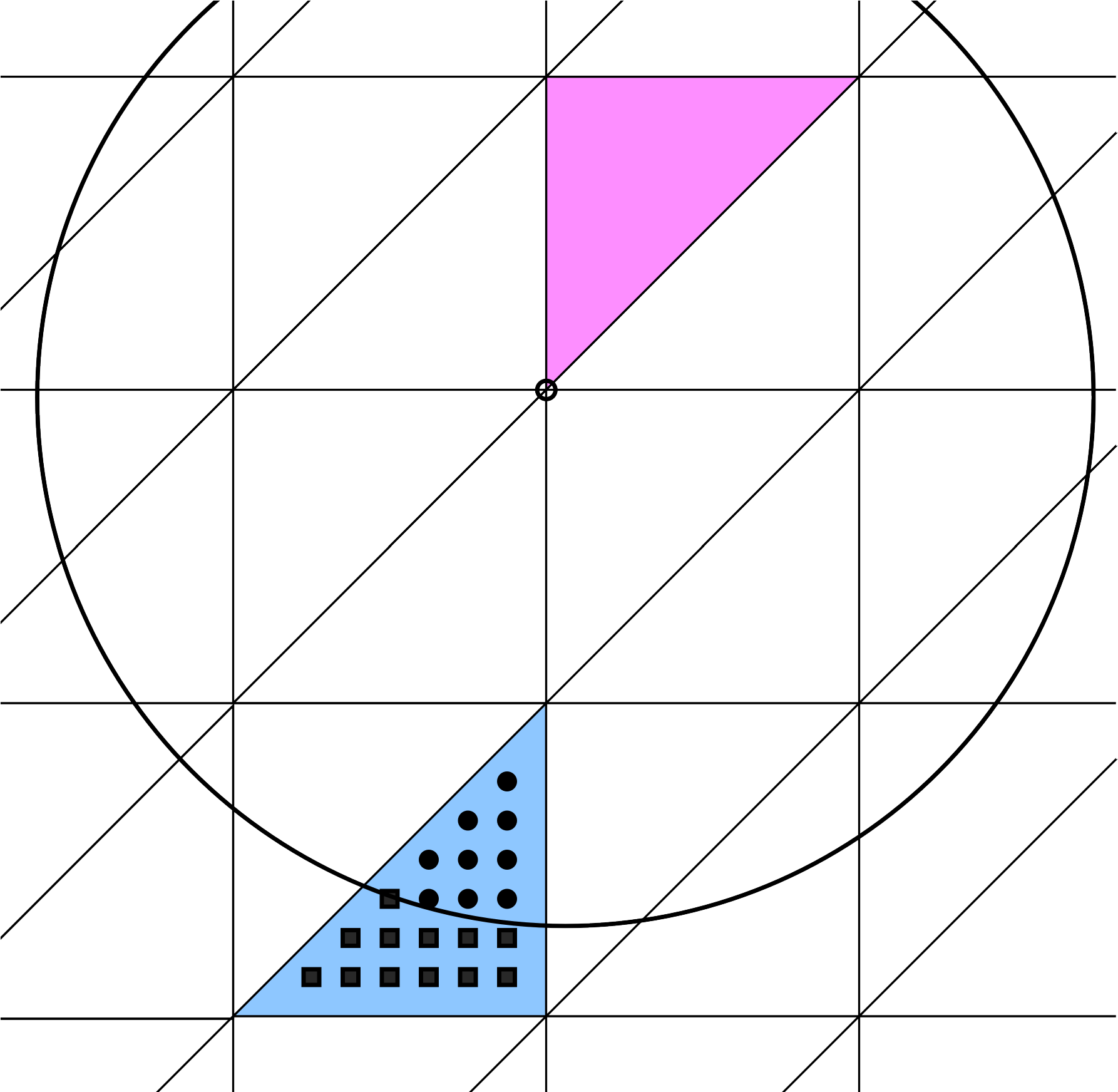}
\\
(d) & (e) & (f)
\end{tabular}
\end{center}
\caption{The blue triangles depict an element $\mcE_k$ for the outer integral. 
(a) The interaction region $\mcI_{\del,\mcE_k}$ for $\mcE_k$ for the ball $B_\del^{exact}(\xb)$ is depicted in orange. The yellow, green, and violet triangles $\mcE_{k'}$ do not, wholly, and partially interact with $\mcE_k$, respectively, with one of the violet triangles having a relatively large interaction area compared to the other.
(b) The orange + red and orange + green triangles are the approximate balls $B_{\del,h}^{overlap}(\xb)$ corresponding to two points in $\mcE_k$.
(c) The same as (b) but for the approximate ball $B_{\del,h}^{barycenter}(\xb)$.
(d) {For the barycenter case, the} support region (in red) of the outer integral over $\mcE_k$ is determined by the intersection of $\mcE_k$ and the ball centered at the barycenter of $\mcE_{k'}$. 
(e) For the barycenter case, the exact support region (in red) and two approximate support regions (in green) using the approximate balls $\{barycenter + nocaps\}$ and $\{barycenter + approxcaps\}$. 
(f) For the barycenter case, an illustration of quadrature points in $\mcE_k$ that are located within the support region (circles) and outside the support region (squares).}
\label{fig:outerq}
\end{figure}

To reveal the difficulties that arise when choosing a quadrature rule for the outer integral triangle $\mcE_k$, we examine the support of the integrand $\mcK^\sharp_{k';j,j'}$ given as 
$$
\bal
\mcS^\sharp_{k,k'} := supp(\mcK^\sharp_{k';j,j'}) 
&= \{\xb\in \mcE_k
\,\,\colon\,\, B^\sharp_\delta(\xb) \cap \mcE_{k'} \neq \emptyset\}
\\&=\mcI_{\mcE_{k'}} \cap \mcE_k =\{\xb\in \mcE_k
\,\,\colon\,\, \exists \,\,\yb\in \mcE_{k'} \,\,\mbox{s.t.}\,\, B^\sharp_\delta(\yb)\cap 
\mcE_{k}\neq \emptyset\},
\eal
$$
where $\mcI_{\mcE_{k'}}$ denotes the interaction domain of $\mcE_{k'}$. 
We have the relations
{ 
\begin{equation*}
 \mcS^{barycenter}_{k,k'}\subset \mcS^{nocaps}_{k,k'}\subset \mcS^{approxcaps}_{k,k'}\subset \mcS^{exactcaps}_{k,k'} =\mcS^{overlap}_{k,k'} \subset \mcE_k.
\end{equation*}
Note that $\mcS^{shifted}_{k,k'} \in \{\mcE_k,\emptyset\}$, depending on whether the shifted ball intersects the inner integral triangle $\mcE_{k'}$ or not. 
}

In what follows, we {\em distinguish between the cases for which $supp(\mcK^\sharp_{k';j,j'}) = \mcE_k$ and the more delicate situation $supp(\mcK_{k';j,j'}) \subsetneqq  \mcE_k $}. {Note that for local PDEs, the second case does not occur because the support of the integrand is always the whole triangle $\mcE_k$.}

\subsection{\textbf{Case 1 -- support of the integrand of the outer integral is the whole outer integral triangle}} \label{subsec:case1}

Consider the case $supp(\mcK^\sharp_{k';j,j'}) = \mcE_k$ almost surely (so that $\mcK^\sharp_{k';j,j'}$ is {almost surely} nonzero for all $\xb\in\mcE_k$) that occurs whenever $\mcE_{k'}$ is wholly contained within the interaction region of $\mcE_k$ as is illustrated by the green triangle in Fig. \ref{fig:outerq}a. This is the simple case that does not require a special treatment of the outer integral. In fact, we can approximate the outer integral using a standard ${Q^{outer}_k}$-point quadrature rule $\{{\xb^{outer}_{k,q}},{w^{outer}_{k,q}}\}_{q=1}^{Q^{outer}_k}$, $k=1,\ldots,K_\Omega$, to obtain, e.g.,
\beq{a_oqr_2}
\bal
&
\int_{\mcE_k }  \mcK^\sharp_{k';j,j'}(\xb) d\xb ~
\approx \sum_{q=1}^{Q^{outer}_k} w^{outer}_{k,q} \mcK^\sharp_{k';j,j'}(\xb^{outer}_{k,q}) .
\eal
\eeq
As discussed in Sec. \ref{quadrules_outer}, a good choice is a four-point symmetric Gaussian quadrature rule of precision three; see Ref. \refcite{AandS}.
 

\subsection{\textbf{Case 2 -- support of the integrand of the outer integral is not the whole outer integral triangle}} \label{subsec:case2}

We now consider the case 
$\mcE^\sharp_k\backslash\mcS^\sharp_{k,k'}\ne\emptyset$ so that 
$supp(\mcK^\sharp_{k';j,j'}) \subsetneqq  \mcE_k $ 
and $\mcK^\sharp_{k';j,j'}$ vanishes on a strict subset of $\mcE_k$ that has positive $d$-dimensional volume. Note that, in this case, triangles $\mcE_{k'}$ are not fully contained within the (approximate) interaction domain of the outer integral triangle $\mcE_{k}$ and thereby are located on the periphery of that interaction domain as is illustrated by the violet triangles in Fig. \ref{fig:outerq}a and the triangle in Fig. \ref{fig:outerq}d having its barycenter depicted by the black dot. As a consequence, there are two issues that arise when choosing a quadrature rule for the outer integral, the first related to precision and the other being a geometric one so that not only the precision of the rule but also the \textit{location of the quadrature points} within the triangle $\mcE_k$ play important roles. Specifically, we discuss the following issues.

\hangseven--~\underline{Lack of smoothness of the integrand} -- For $\sharp\in\{barycenter,\,overlap\}$ the integrand $\mcK^\sharp_{k';j,j'}$ is {\em discontinuous} on $\mcE_k$ (more details below). For the exact ball and the other approximate balls, the integrand is continuous but may {\em not be differentiable} on $\mcE_k$. As a result, the accuracy of any quadrature rule on $\mcE_k$ that requires greater smoothness may be compromised.

\hangseven--~\underline{Missing triangles} -- If \emph{all} quadrature points $\xb^{outer}_{k,q}$ are located in the complement {$\mcE_k \setminus\mcS^\sharp_{k,k'}\ne\emptyset$} of the support of the integrand, then the double integral \eqref{def:outer_task_approx} is approximated by zero despite the fact that $\mcE_k$ and $\mcE_{k'}$ are a pair of interacting elements.

\noindent In the next two subsections we provide details about how these two issues arise and how they influence the choice of the quadrature rule for an outer integral triangle. 

\subsubsection{\textbf{Lack of smoothness in the integrand}}\label{subsubsec:discontinuity_outer}

We divide the discussion into four sub-cases because the issue ensuing from a lack of smoothness differs between them, as are the mitigating approaches for addressing the issue.

\underline{The $\sharp \in \{exactcaps, approxcaps, nocaps\}$ cases}.
For these cases, the integrand $\mcK_{k';j,j'}$ is continuous on $\mcE_k$ but may not be smoother than that. Because the support region $supp(\mcK^\sharp_{k';j,j'})$ is a strict subset of $  \mcE_k $, for a chosen quadrature rule on $\mcE_k$, some of the quadrature points may be located in the complement domain $\mcE_k\setminus\mcS^\sharp_{k,k'}$ on which the integrand vanishes; see Fig. \ref{fig:outerq}f for an illustration. Thus, the accuracy of a quadrature rule defined over all of $\mcE_k$ may be corrupted, i.e., it does not achieve its full potential accuracy, because {\em the integrand is not sufficiently smooth over $\mcE_k$}. The resulting approximations of the outer integrals then take the form of \eqref{a_oqr_2}. The numerical results presented in Sec. \ref{numerics} give rise to the conjecture that the seven-point rule (i.e., $Q^{outer}_k = 7$) in Fig. \ref{fig:7points} does not only fulfill an important placement feature (as illuminated in Sec. \ref{missing}) but also consists of sufficiently many quadrature points to produce stable second-order convergence rates for the exact ball as well as the ball approximations $\sharp\in\lk approxcaps,~nocaps\rk$.

\underline{The $\sharp = overlap$ case}.
For the $overlap$ ball approximation $B^{overlap}_{\delta,h}(\xb)$, the situation is even worse because, in this case, the integrand
$$
\bal
&\mcK^{overlap}_{k';j,j'}(\xb) 
\\&\quad= 
\left\{
\bal
&\int_{\mcE_{k'}} \big(\phi_{j}(\xb)-\phi_{j}(\yb)\big)\big(\phi_{j'}(\xb)-\phi_{j'}(\yb)\big)\psi(\xb,\yb)d\yb \ne 0
\\
&\qquad\qquad\qquad\qquad\qquad\qquad\qquad\qquad\qquad
\mbox{if\,\,\,}
\mcE_{k'} \cap B_{\del}(\xb)  \neq \emptyset
\\
&0\qquad\qquad\qquad\qquad\qquad\qquad\qquad\qquad\quad\;\;
\mbox{if\,\,\,}
\mcE_{k'} \cap B_{\del}(\xb) = \emptyset
\eal\right.
\eal
$$
has a jump discontinuity within $\mcE_k$, i.e., for $\xb$ such that the overlap $\mcE_{k'} \cap B_{\del}(\xb)$ is tiny, the whole element $\mcE_{k'}$ interacts with that $\xb$ but a slight change in the position of $\xb$ can cause the overlap to vanish, in which case $\mcE_{k'}$ no longer interacts with $\xb$. Fig. \ref{fig:outerq}b illustrates the strong dependence of the approximate ball $B^{overlap}_{\delta,h}(\xb)$ on where in $\mcE_{k}$ the point $\xb$ is located.

This issue is particularly preponderant if either the interaction horizon $\delta$ is comparable to the grid size or if $\delta$ is small compared to the grid size. For both these situations, Case 2 dominates for pairs of interacting triangles $(\mcE_k,\mcE_{k'})$. One is also naturally confronted with this issue when aiming to numerically investigate the local limit as $\delta\to 0$ for a fixed finite element mesh.

In order to handle the difficulty caused by the discontinuity of the integrand, it is best to numerically identify the support region $S^{overlap}_{k,k'} = \mcI^{overlap}_{\mcE_{k'}} \cap \mcE_k= \{\xb\in \mcE_k\,\colon\, B_{\delta}(\xb) \cap \mcE_{k'} \neq \emptyset\}$ and then only place quadrature points inside this region. However, this approach is computationally expensive because $\mcI^{overlap}_{\mcE_{k'}}\cap \mcE_k $ is determined by infinitely many ball intersections. Another approach is to use adaptive quadrature rules that automatically take care of the determination of the support. However, because the evaluation of $\mcK^\sharp _{k';j,j'}$ at a point $\xb\in\mcE_k$ is expensive, one wants to avoid as many function evaluations as possible. An in between approach is to use a quadrature rule $\{{\xb^{outer}_{k,q}},{w^{outer}_{k,q}}\}_{q=1}^{Q^{outer}_k}$ that consists of more points than are used in Case 1.

\underline{The $\sharp = barycenter$ case}.
The integrand $\mcK^{barycenter}_{k';j,j'}(\xb)$ corresponding to the barycenter based ball approximation also has a jump discontinuity on $\mcE_k$ because even a slight change in the position of a point $\xb\in\mcE_k$ can cause the barycenter of the element $\mcE_{k'}$ to be inside or outside the ball $B_{\delta}(\xb)$. Fig. \ref{fig:outerq}c illustrates the strong dependence of the approximate ball $B^{barycenter}_{\delta,h}(\xb)$ on where in $\mcE_{k}$ the point $\xb$ is located. However, unlike the overlap ball case, for the barycenter ball case one can numerically determine the support region.

More precisely, by definition we have that
$$ 
\bal
&\mcK^{barycenter}_{k';j,j'}(\xb) 
\\&\quad= 
\left\{
\bal
&\int_{\mcE_{k'}} \big(\phi_{j}(\xb)-\phi_{j}(\yb)\big)\big(\phi_{j'}(\xb)-\phi_{j'}(\yb)\big)\psi(\xb,\yb)d\yb \ne 0
\\
&\qquad\qquad\qquad\qquad\qquad\qquad\qquad\qquad\qquad
\mbox{if\,\,\,}
|\xb-\xb^{barycenter}_{k'}|\le\del
\\
&0\qquad\qquad\qquad\qquad\qquad\qquad\qquad\qquad\quad\;\;
\mbox{if\,\,\,}
|\xb-\xb^{barycenter}_{k'}|>\del
\eal\right.
\eal
$$
so that the resulting support region can be characterized as
$$
S^{barycenter}_{k,k'} = \{\xb \in \mcE_k \colon  |\xb-\xb^{barycenter}_{k'}|\le\del \} = \mcE_k \cap B_\del(\xb^{barycenter}_{k'}),
$$
i.e., the support region is determined as the intersection of the outer element $\mcE_{k}$ with the ball of radius $\del$ centered at the barycenter of the element $\mcE_{k'}$.

In contrast to the cases $\sharp\in\{exactcaps,nocaps,approxcaps\}$, the support region $S^{barycenter}_{k,k'} = \mcE_k \cap B^{barycenter}_\del(\xb_{k'})$ is characterized by exactly one ball intersection; see Fig. \ref{fig:outerq}d for an illustration. As a result we can apply one of the ball approximations $B^{nocaps}_{\del,h}(\xb_{k'})$ or $B^{approxcaps}_{\del,h}(\xb_{k'})$ introduced in Sec. \ref{approxballs} to $B_\del(\xb^{barycenter}_{k'})$ or even the exact ball $B_\del(\xb)$ in order to define a composite quadrature rule for the outer integral triangle in the fashion of Sections \ref{ecaps} and \ref{cqr}; see the three examples in Fig. \ref{fig:outerq}e. In Fig. \ref{fig:outerq}f, we illustrate how a quadrature rule on an outer integral triangle $\mcE_k$ for which the support of its integrand is only a portion of $\mcE_k$ {may have} quadrature points that lie outside that support.

\underline{The $\sharp = shifted$ case}. In contrast to the exact ball and all other approximate balls, the shifted approximate ball $B_{\del,h}^{shifted}(\xb) = B_{\del}(\xb_k^{barycenter})$ is a special case in that it does not depend on $\xb \in \mcE_k$, i.e., all quadrature points in $\mcE_k$ use the same ball $B_{\del}(\xb_k^{barycenter})$ to determine which inner elements $\mcE_{k'}$ they interact with. Thus we have the integrand 
$$
\bal
&\mcK^{shifted}_{k';j,j'}(\xb) 
\\&\quad= 
\left\{
\bal
&\int_{\mcE_{k'} \cap  B_{\del,h}(\xb_k^{{barycenter}})} \big(\phi_{j}(\xb)-\phi_{j}(\yb)\big)\big(\phi_{j'}(\xb)-\phi_{j'}(\yb)\big)\psi(\xb,\yb)d\yb \ne 0
\\
&\qquad\qquad\qquad\qquad\qquad\qquad\qquad\qquad\qquad
\mbox{if\,\,\,}
\mcE_{k'} \cap B_{\del}(\xb_k^{barycenter})  \neq \emptyset
\\
&0\qquad\qquad\qquad\qquad\qquad\qquad\qquad\qquad\quad\;\;
\mbox{if\,\,\,}
\mcE_{k'} \cap B_{\del}(\xb_k^{barycenter}) = \emptyset
\eal\right.
\eal
$$
so that $\mcS_{k,k'}^{shifted} \in \{\mcE_k,\emptyset\}$. Thus, the discontinuity issue does not arise because $\mcK^{shifted}_{k';j,j'}(\xb)$ is either nonzero or zero for all $\xb\in \mcE_k$. Therefore, for the shifted ball approximation, we can use the same quadrature rule as that chosen for the outer integral in {Case 1} in Sec. \ref{subsec:case1}. 

\subsubsection{\textbf{Missing triangles -- affecting the location of quadrature points}}\label{missing}

By using a quadrature rule with quadrature points that are interior to $\mcE_k$ and that has the minimum number of quadrature points needed for exact integration of cubic polynomials on triangles (see also Sec. \ref{quadrules_outer}), one can miss interactions between the outer integral triangle $\mcE_k$ and an inner integral triangle $\mcE_{k'}$. This is precisely the case if \emph{all} quadrature points $\xb^{outer}_{k,q}$ are located in the complement of the support of the integrand $\mcK^\sharp_{k';j,j'}(\xb)$. 

This observation is illustrated in Fig. \ref{fig:missing-triangles}. In (a), the violet area indicates the interaction region of the blue outer integral triangle $\mcE_k$, i.e., $\mcI_{\mcE_k} = \{\yb : |\xb-\yb|\leq\delta \; {\rm for}\; \xb\in\mcE_k\}$. In (b) and (c), the orange area indicates the union of the balls centered at three quadrature points in $\mcE_k$ indicated by the black dots. For simplicity, we are using exact balls but similar pictures would hold for approximate balls with the exception of the shifted ball for which there is only a single ball for all quadrature points.
In (b), the points are interior to $\mcE_k$ whereas for (c) they are at the vertices. We see that the three vertices result in much better coverage of the true interaction region $\mcI_{\mcE_k}$ than do the three interior points. Still, a vertex rule may miss an inner integral triangle that interacts with $\mcE_k$ as depicted in (d), with a zoom-in in (e). More precisely, the black part of the inner integral triangle $\mcE_{k'}$ colored in red and black overlaps with the interaction domain of the outer integral triangle $\mcE_{k}$ so that those two triangles interact. However, {\em because that black region does not intersect the orange region, the contribution of the two interacting triangles $\mcE_{k}$ and $\mcE_{k'}$ is missed.} Looking at (f), we see that by adding the midpoints of the sides of the triangle $\mcE_{k}$ to the vertex points results in even better coverage of the true interaction domain and thus there is even less likelihood that a triangle will be missed compared to just having vertex points. In (g), the orange triangles are those that overlap with one or more of the three balls and in (h) the same is true for the orange and magenta triangles, with the magenta triangles are those that are missed in (g). In fact, in (h), no triangles are missed, i.e., the magenta and orange triangles account for all triangles that intersect the true interaction region for $\mcE_k$.

\begin{figure}[ht]
\begin{center}
\begin{tabular}{ccc}
\includegraphics[width=1.3in]{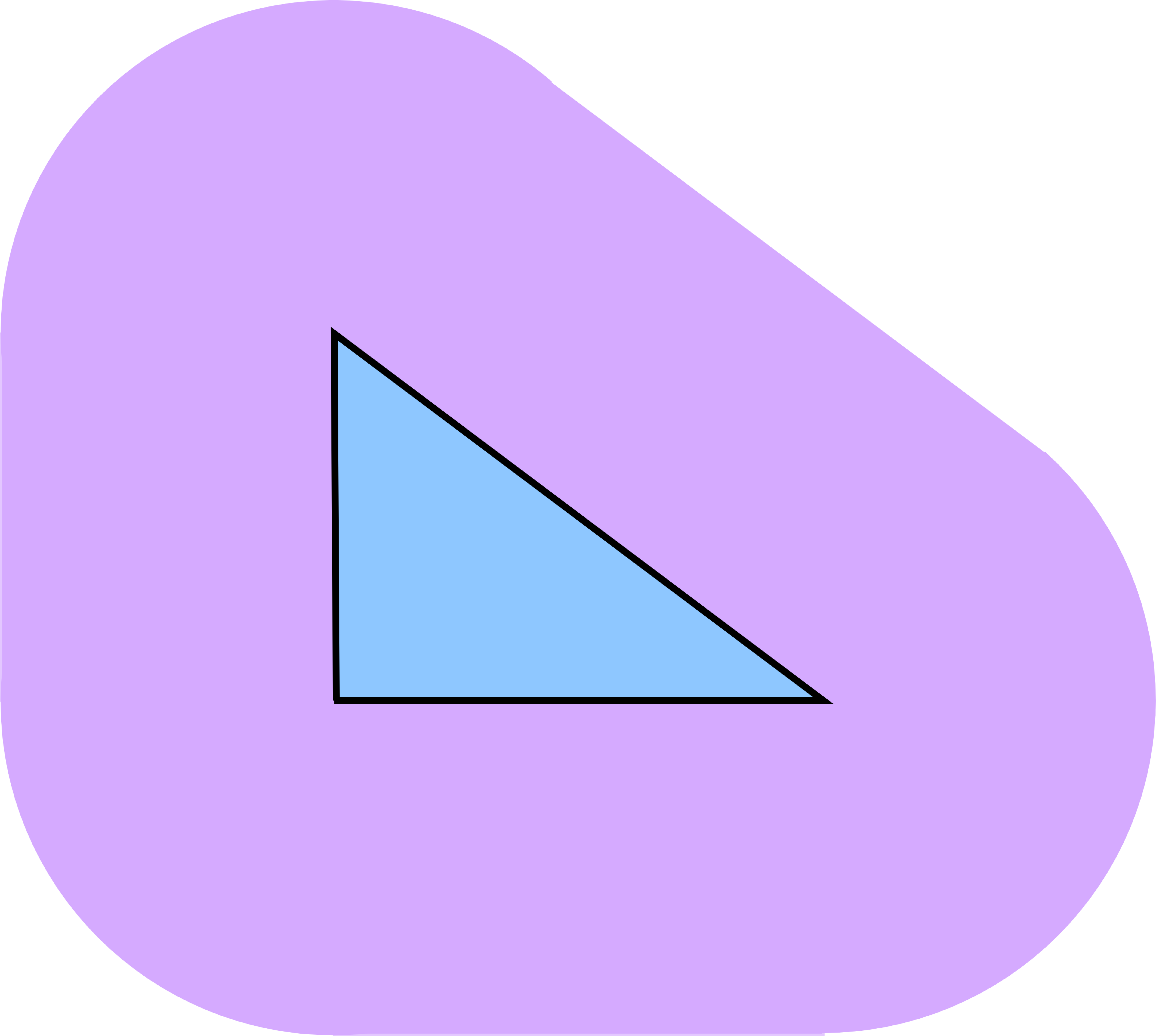}
&
\includegraphics[width=1.3in]{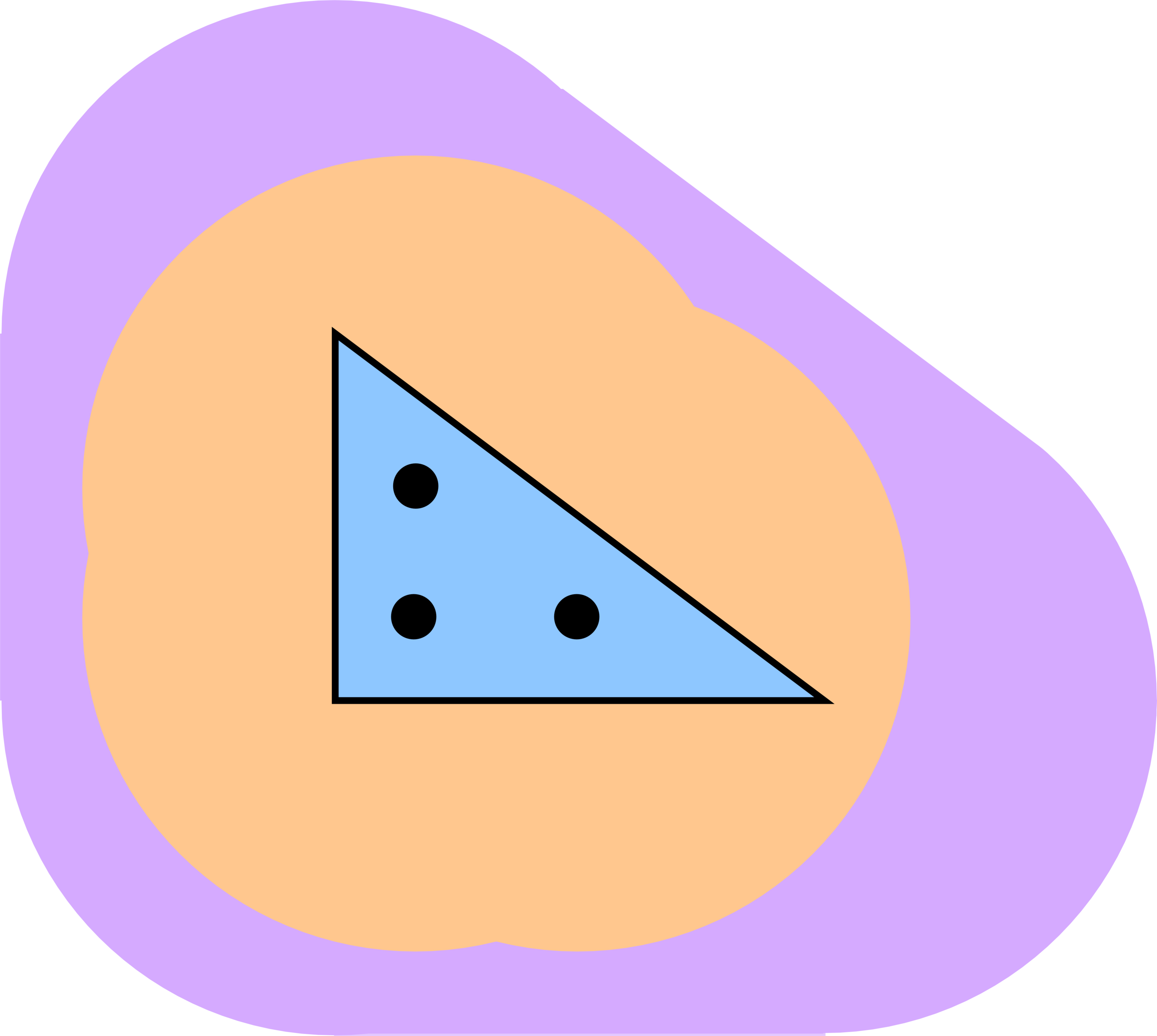}
&
\includegraphics[width=1.3in]{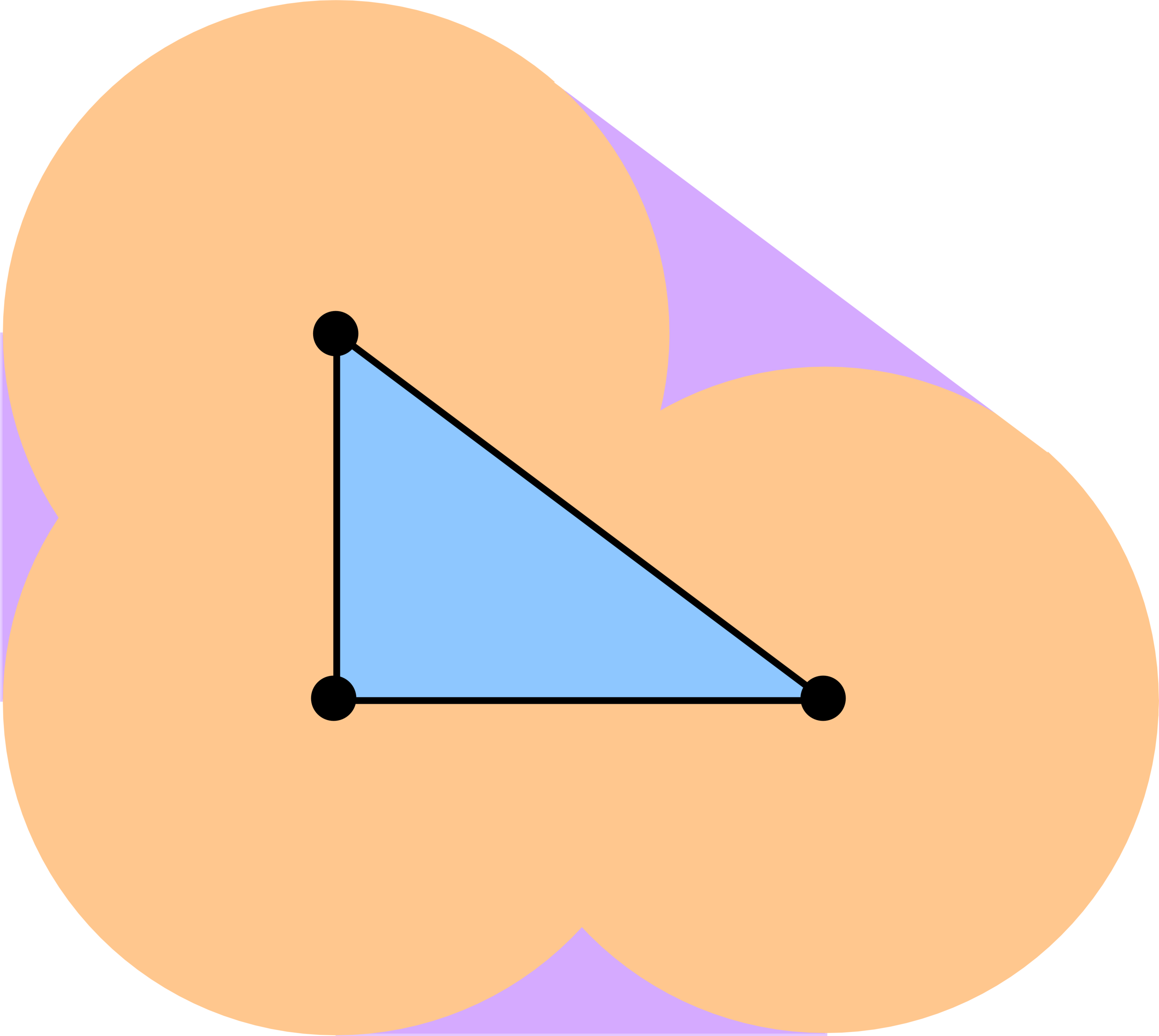}
\\
(a) & (b) & (c)
\\
\includegraphics[width=1.3in]{figures/missing-triangle4.pdf}
&
\includegraphics[width=1.2in]{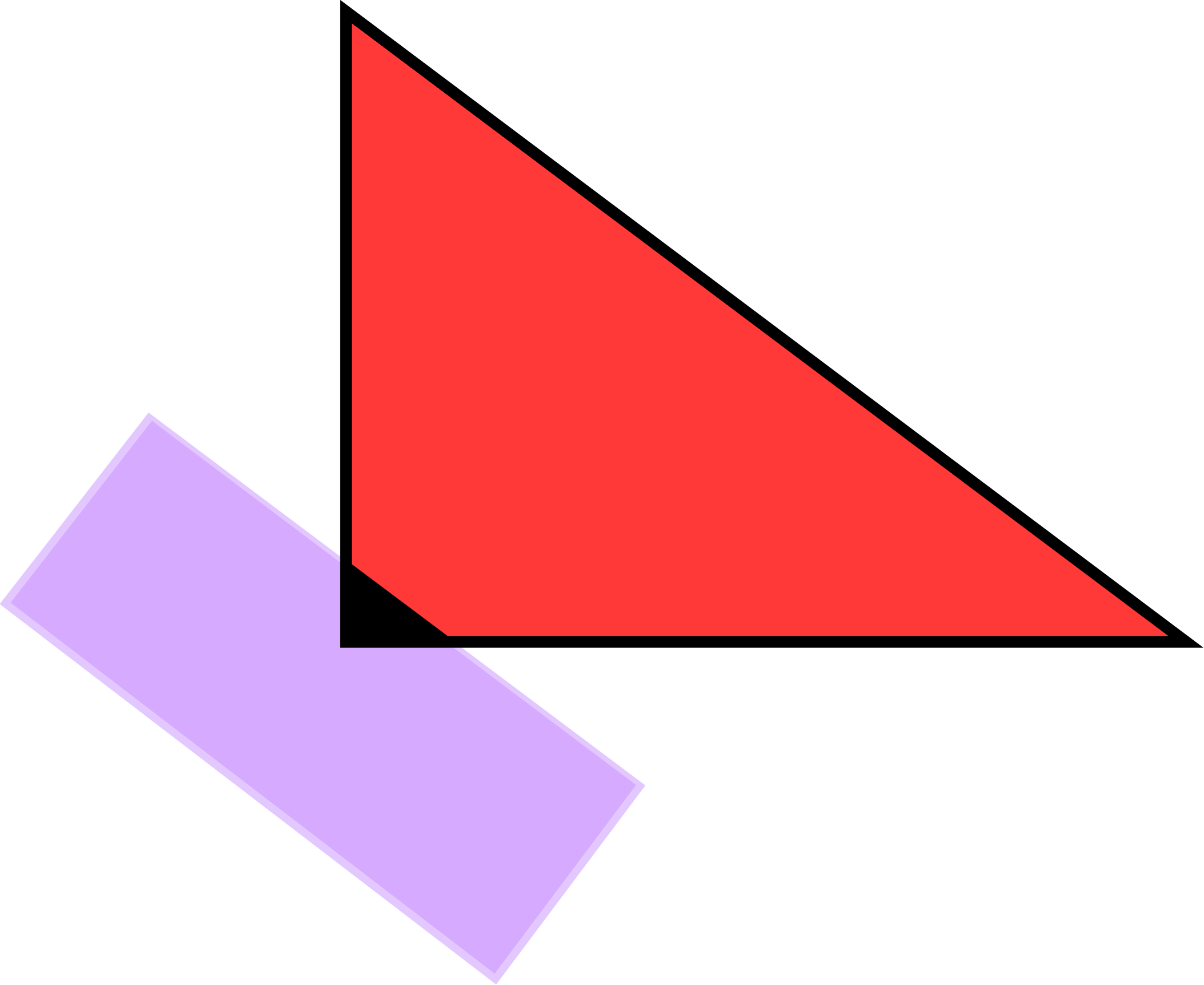}
&
\includegraphics[width=1.3in]{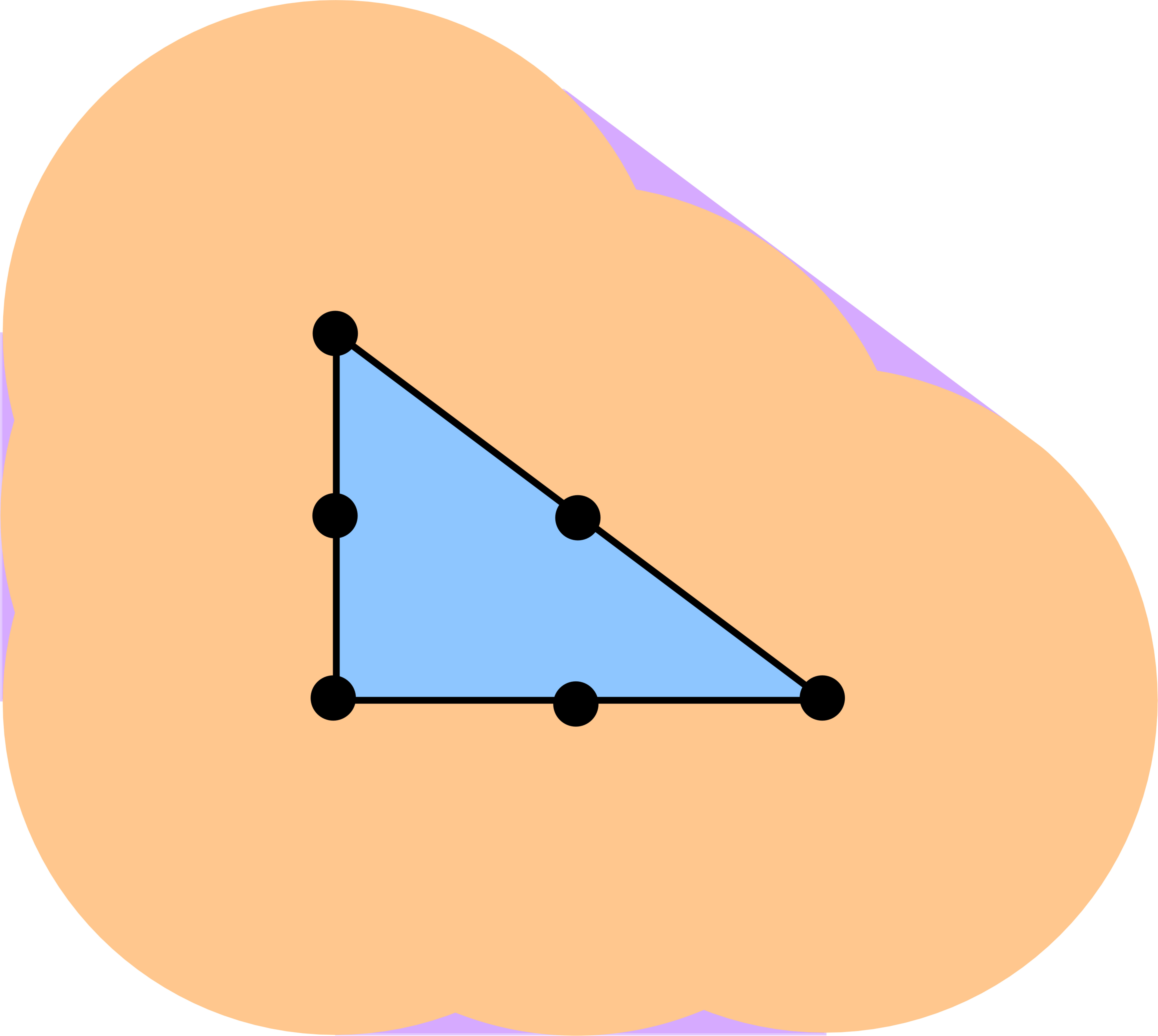}
\\
(d) & (e) & (f)
\\[1ex]
\includegraphics[width=1.3in]{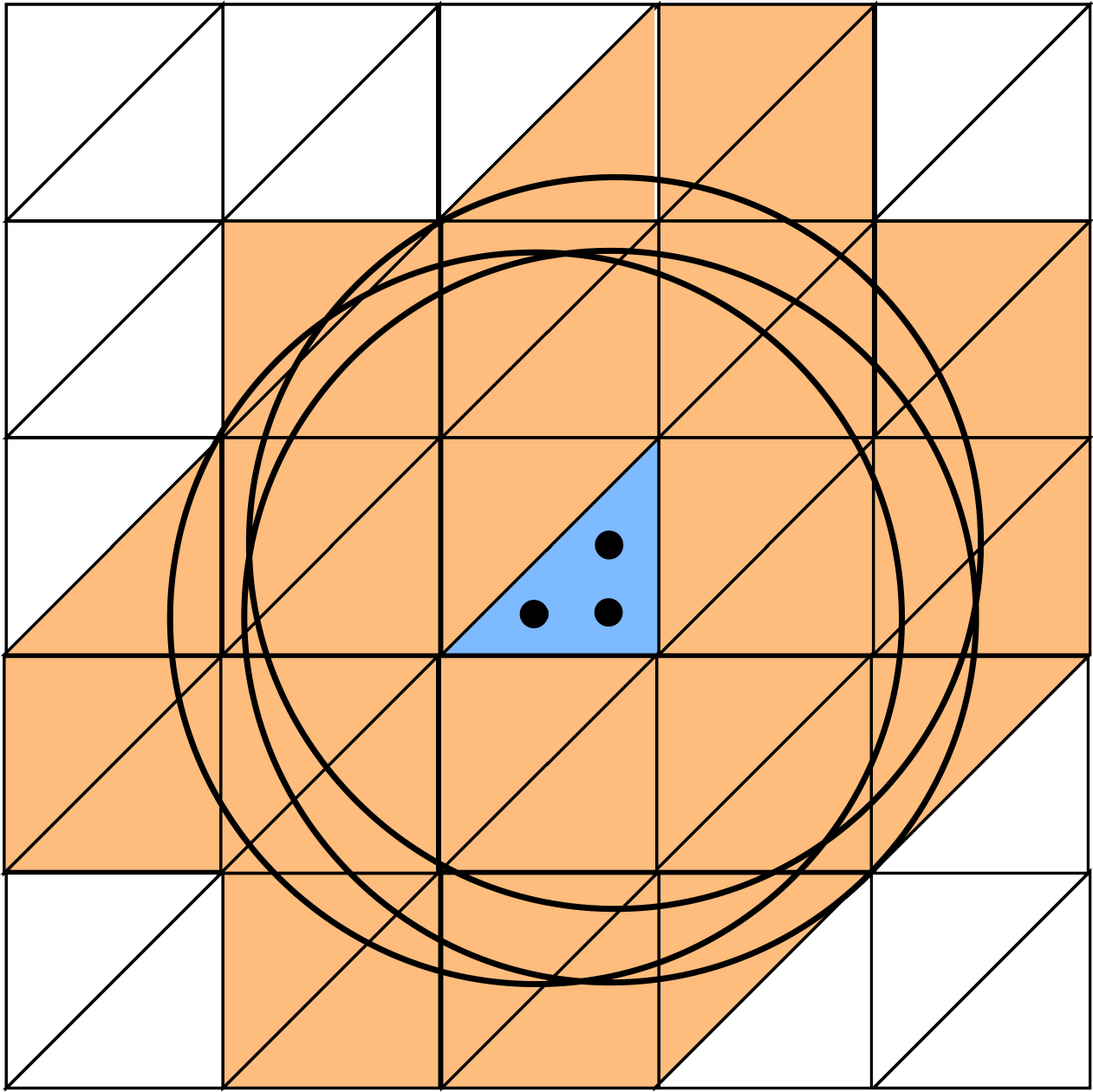}
&
\includegraphics[width=1.3in]{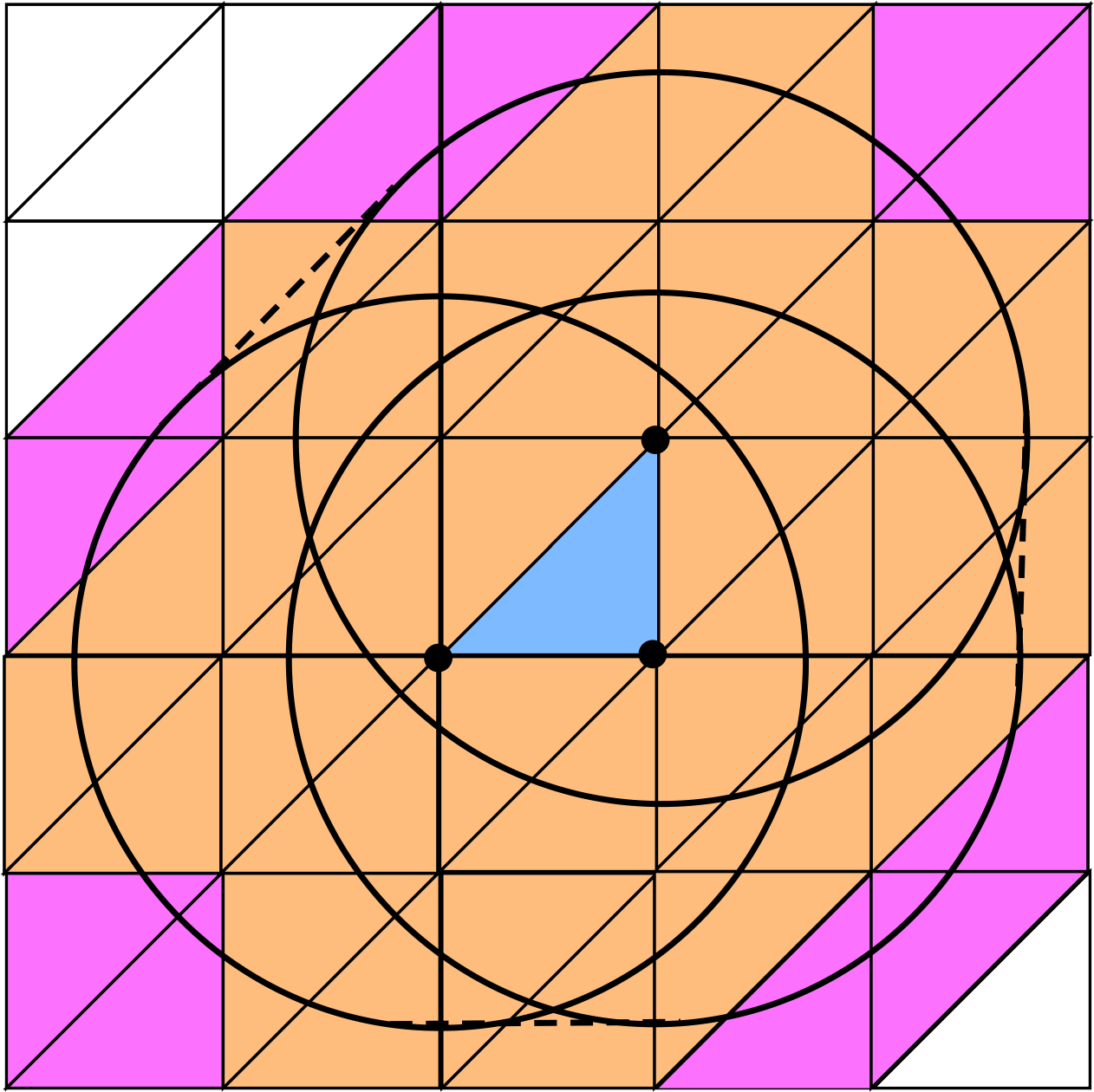}
&
\includegraphics[width=1.4in]{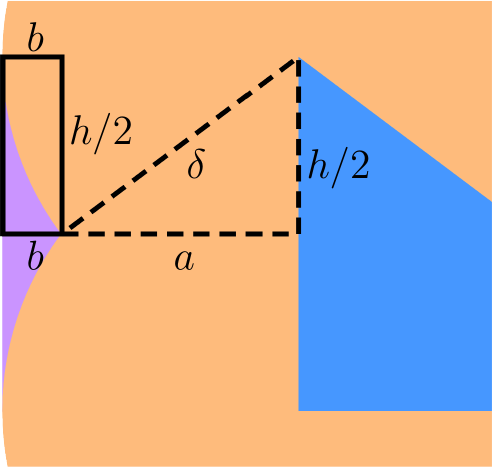}
\\
(g) & (h) & (i)
\end{tabular}
\end{center}
\caption{Illustrations related to missed triangles.}
\label{fig:missing-triangles}
\end{figure}

A simple computation shows that the difference between the violet and orange areas in Fig. \ref{fig:missing-triangles}c, and therefore also in Fig. \ref{fig:missing-triangles}f, is of order $\mcO(h^3)$. {Using the notation of Fig.  \ref{fig:missing-triangles}i, we have that $a+b=\delta$, $a^2 +\frac14 h^2 =\delta^2$, and the area $R$ of the rectangle is $\frac12bh$ so that, for fixed $\delta$ and small $h$,
$$
   a \approx \delta - \frac{h^2}{8\delta},
\qquad
 b \approx  \frac{h^2}{8\delta},
 \qquad\mbox{and}\qquad
 R  \approx \frac{h^3}{16\delta}.
$$
The area of each violet region in Fig. \ref{fig:missing-triangles}c is less than twice the area of the rectangle in Fig. \ref{fig:missing-triangles}i.} Clearly, the area of the missing triangle, depicted in black in Figs. \ref{fig:missing-triangles} (d) and (e), is then also of $\mcO(h^3)$. Of course, this means the violet area in Fig. \ref{fig:missing-triangles}f is also of $\mcO(h^3)$ but with a substantially smaller constant in the order relation. We note that for the configuration of \ref{fig:missing-triangles}b for which the quadrature points are usually at a distance of $\mcO(h)$ away from the vertices, the violet area is of $\mcO(h^2)$.

The barycenter based polytopial ball approximation misses additional inner integral triangles due to its definition. In fact, it misses precisely those $\mcE_{k'}$ for which the barycenter is not contained in the interaction set of $\mcE_k$, i.e., $\xb_{k'}^{barycenter} \notin \mcI^{barycenter}_{\del,h}$. Due to its dependence on $\xb$ it may miss even more interacting triangles due to an inconvenient choice of quadrature rules. However, by employing a composite quadrature rule on $S_{k,k'}^{barycenter}=\mcE_k \cap B_\del(\xb_{k'}^{barycenter})$, as proposed in the Sec. \ref{subsubsec:discontinuity_outer}, we do not only circumvent the discontinuity of $\mcK^{barycenter}_{k;j,j'}$ but also only neglect the conceptually missed interacting triangles.

Similarly, the shifted ball approximation misses interacting inner integral triangles due to its definition. In fact, the approximate interaction domain of $\mcE_k$ is given by
$$\mcI_k^{shifted} = \{ \yb \in \R^d\backslash \mcE_k\colon~ \yb \in B_{\del}(\xb_k^{barycenter}) \} = B_{\del}(\xb_k^{barycenter})\backslash \mcE_k  .$$
Therefore the set of missed triangles is composed of those $\mcE_{k'}$ for which $\mcE_{k'} \cap B_{\del}(\xb_k^{barycenter})=\emptyset$ and it cannot be affected by the choice of quadrature rules.

\subsubsection{\textbf{Heuristics about the choice of quadrature rules in Case 2}}\label{quadrules_outer}

Let us continue the reasoning in Sec. \ref{quade_inner} about the choice of quadrature rules. For this purpose, we suppose that the inner integrals in \eqref{cqr3} and \eqref{cqr1} are integrated exactly. Then, for piecewise-linear basis functions and again assuming that the kernel function $\psi(\xb,\yb)$ is constant, we have that the integrand $\mcK^\sharp_{k';j,j'}(\xb)$ of the outer integral is a polynomial of degree $3$ in the components of $\xb$. Thus, for a typical outer integral triangle $\mcE_{k}$, heuristically one should use a quadrature rule $\{{\xb^{outer}_{k,q}},{w^{outer}_{k,q}}\}_{q=1}^{Q^{outer}_k}$ of precision $3$ for the outer integral. A four-point symmetric Gaussian quadrature rule of precision three (see Ref. \refcite{AandS}) would suffice for this purpose.

Commensurate quadrature rules that result in an $\mcO(h^2)$ approximation use even fewer quadrature points, so they in general would result in the missing triangle syndrome.

\subsection{Final word on choosing a quadrature rule for the outer integral}

The discussion in Sec. \ref{quadrules_outer} focused only  on precision, but as we have seen, quadrature point placement also is important. Thus, in choosing the quadrature points for the outer integrand, not only do we have to guarantee a sufficiently accurate integration of the integrand, but also have enough well-placed quadrature points so that either we do not miss any inner integral triangles $\mcE_{k'}$ or such that the missed triangles have a negligible contribution to the integration.

A precision-three rule that includes the vertices of the outer integral triangle $\mcE_k$ seemingly can satisfy both the precision requirement stemming from the heuristic approach of Sec. \ref{quadrules_outer} and the point-placement requirement of Sec. \ref{missing}. Specifically, the seven-point rule having quadrature points at the barycenter, the vertices, and the mid-side points and the corresponding weights are $\frac{27}{60}\!\cdot\!\frac12$, $\frac{3}{60}\!\cdot\!\frac12$, and $\frac{8}{60}\!\cdot\!\frac12$, respectively, has precision $3$ (Ref.~\refcite{AandS}) and includes vertex points; see Fig. \ref{fig:7points}. Note that the factor $\frac12$ in the weights is the area of the reference triangle. This rule has the bonus feature of including mid-side quadrature points so that missing triangles are unlikely to affect the overall accuracy.

Note that the seven-point rule of Fig. \ref{fig:7points} is not optimal with respect to the number of points; 4-point
precision-three rules such as the one mentioned in Sec. \ref{quadrules_outer} are known to integrate cubics exactly and seven-point rules exist that integrate quintics exactly. It is not optimal even among quadrature rules that include vertex points because a six-point rule with three additional judiciously placed interior points can have precision 3.  However, the rule of Fig. \ref{fig:7points} is a precision-three rule having the minimum number of points, if vertices and midsides have to be included. The aforementioned bonus of having mid-side quadrature points leads us to the seven-point rule of Fig. \ref{fig:7points} as the quadrature rule of choice.
\begin{figure}[ht]
\centering
\includegraphics[width=0.2\textwidth]{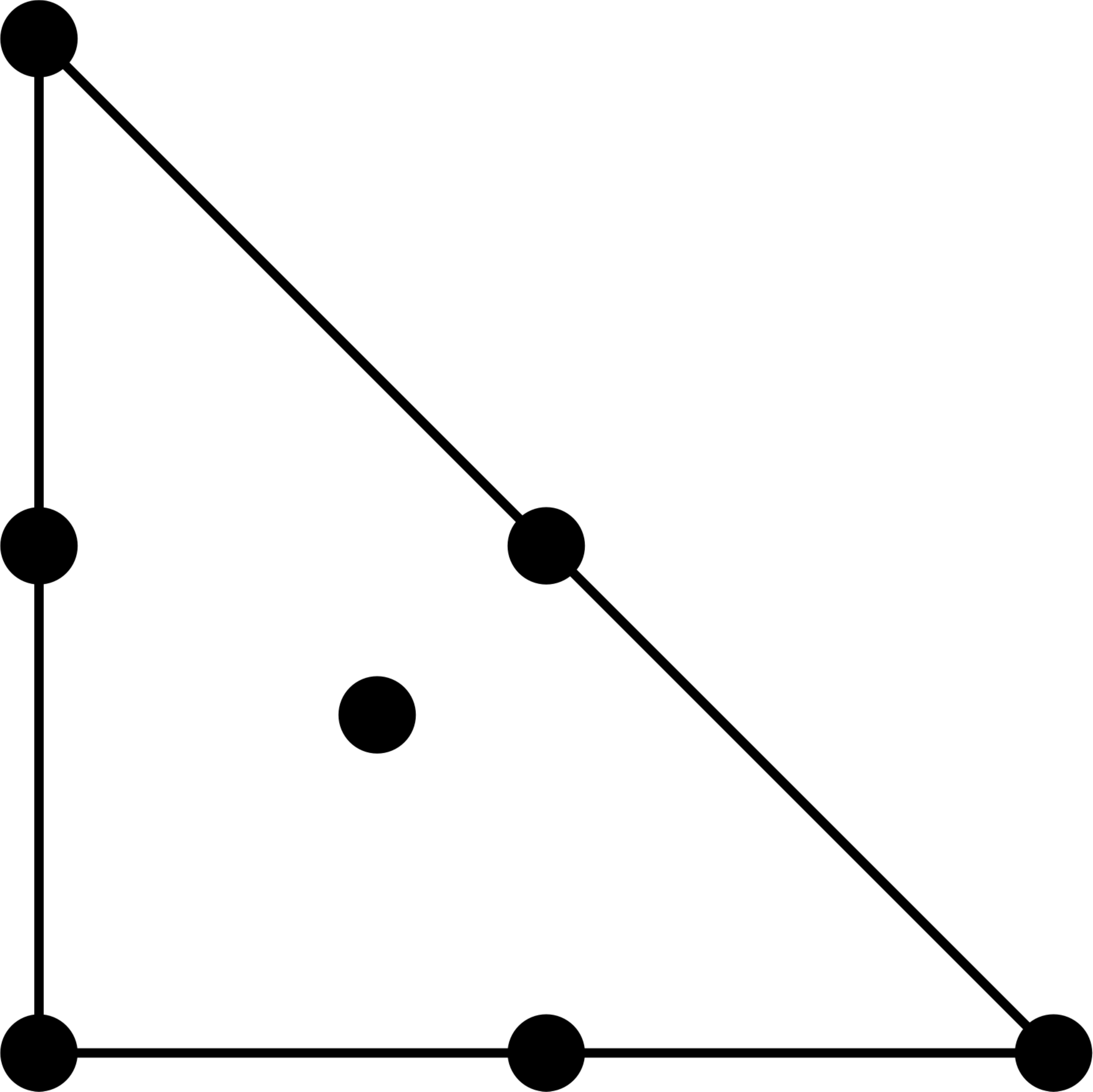}
\caption{Illustration of the nodes of a seven-point rule that integrates cubic polynomials exactly and includes the vertices and mid-sides of the triangle.}
\label{fig:7points}
\end{figure}

{\section{Efficient implementation}\label{efficient}

\subsection{\textbf{Tasks for polytopial approximate ball construction}}
\label{ballconstr}

In this section, we provide details about how the six tasks listed in Sec. \ref{polyballs} can be efficiently executed. We assume that we have in hand a finite element mesh (see Sec. \ref{femgrid}) having maximum grid size $h_{max}$ and minimum grid size $h_{min}$ and a ball $B_\del(\xb)$ having radius $\del$ and centered at a point $\xb\in\OuO$.

\hangfive 1.~{\em Determination of the location of the barycenter of an element.} This task is easily accomplished because the coordinates of the barycenter are simply the average of the coordinates of the vertices of the element.

\hangfive 2.~{\em Identification of elements that intersect the ball.} Let $\mcE_k$ denote a fixed outer integral triangle. Then, during the inner assembly loop, we consider all inner integral triangles $\mcE_{k'}$ for which  $|\xb^{barycenter}_k-\xb^{barycenter}_{k'}| < \delta + h_{max}$. Thus, there may be $\xb \in \mcE_k$ for which $B_\del(\xb)\cap\mcE_{k'} = \emptyset$. However, these cases are automatically identified by the following routines. Alternatively, one could also implement some type of breadth-first search.

\hangfive 3.~{\em Identification of elements wholly contained within a ball.} If all the vertices of an element are contained within the ball, then the whole element is contained within the ball, i.e., ${\mathcal E}_k\cap B_{\del}(\xb)={\mathcal E}_k$. {Thus in order to identify elements of this type we have to compute the Euclidean distance between the three vertices and the midpoint $\xb$ of the ball.}

\hangfive 4.~{\em Identification of elements that partially overlap with a ball.} If one or two but not three vertices of an element are inside the ball, that element only partially overlaps with the ball so that the identification of such elements is an easy matter; see Figs. \ref{balltriangle}a and \ref{balltriangle}b for examples of one and two vertices being inside the ball, respectively. However, it is possible for an element to intersect the ball without having an element vertex inside the ball, a situation that occurs when the boundary of the ball intersects a single element edge at two points; see Fig. \ref{balltriangle}c. {In order to identify when this situation occurs we also compute the set of intersection points resulting from intersecting the boundary of the ball with the boundary of the element (see next task). If there are two such intersection points but no element vertex inside the ball, then we have identified a partially covered triangle of the latter kind.}

\hangfive 5.~{\em Identification of the points at which the boundary of the ball intersects the boundary of the elements.} The boundary of a ball may intersect the boundary of an element in several different ways. For example, in Figs. \ref{balltriangle}a to \ref{balltriangle}c, there are two intersection points whereas in Fig. \ref{balltriangle}d there are four. There are other configurations for the intersection of balls and triangles; see Ref. \refcite{feifei2}; the ones depicted in Fig. \ref{balltriangle} are the possibilities that exist if the diameter $2\del$ of the ball is larger than the diameter of the triangle. {To identify the intersection points we intersect each side of the triangle with the boundary of the ball by solving the determining quadratic equations. More precisely, let $\{\vb_1, \vb_2, \vb_3 \}$ denote the vertices of a finite element triangle. Then by solving the quadratic equation $q(\lambda) = |\vb_i + \lambda (\vb_j-\vb_i) - \xb |^2 -\delta^2 = 0$, where $i\neq j \in \{1,2,3\}$, for $\lambda \in [0,1]$, we find the intersection points.}

\hangfive 6.~{\em Determining a subdivision of a polygon into triangles.} {We have determined the element vertices which lie inside the ball and the points at which the boundary of the ball intersects the boundary of the elements. If the union of these points is larger than two, then we can define a polytopial approximation to the convex intersection region. For this purpose we first order these points (counter-)clockwise which results in an ordered set of points $\{\pb_1,\ldots, \pb_n\}$ for $3\leq n\leq 6$. A subdivision into triangles is then given by $\{\{\pb_1,\pb_{i+1},\pb_{i+2}\}: ~\text{for}~  i=1,\ldots, n-2$\}.}
\begin{figure}[ht]
\centering
\begin{tabular}{cc}
\includegraphics[height=.9in]{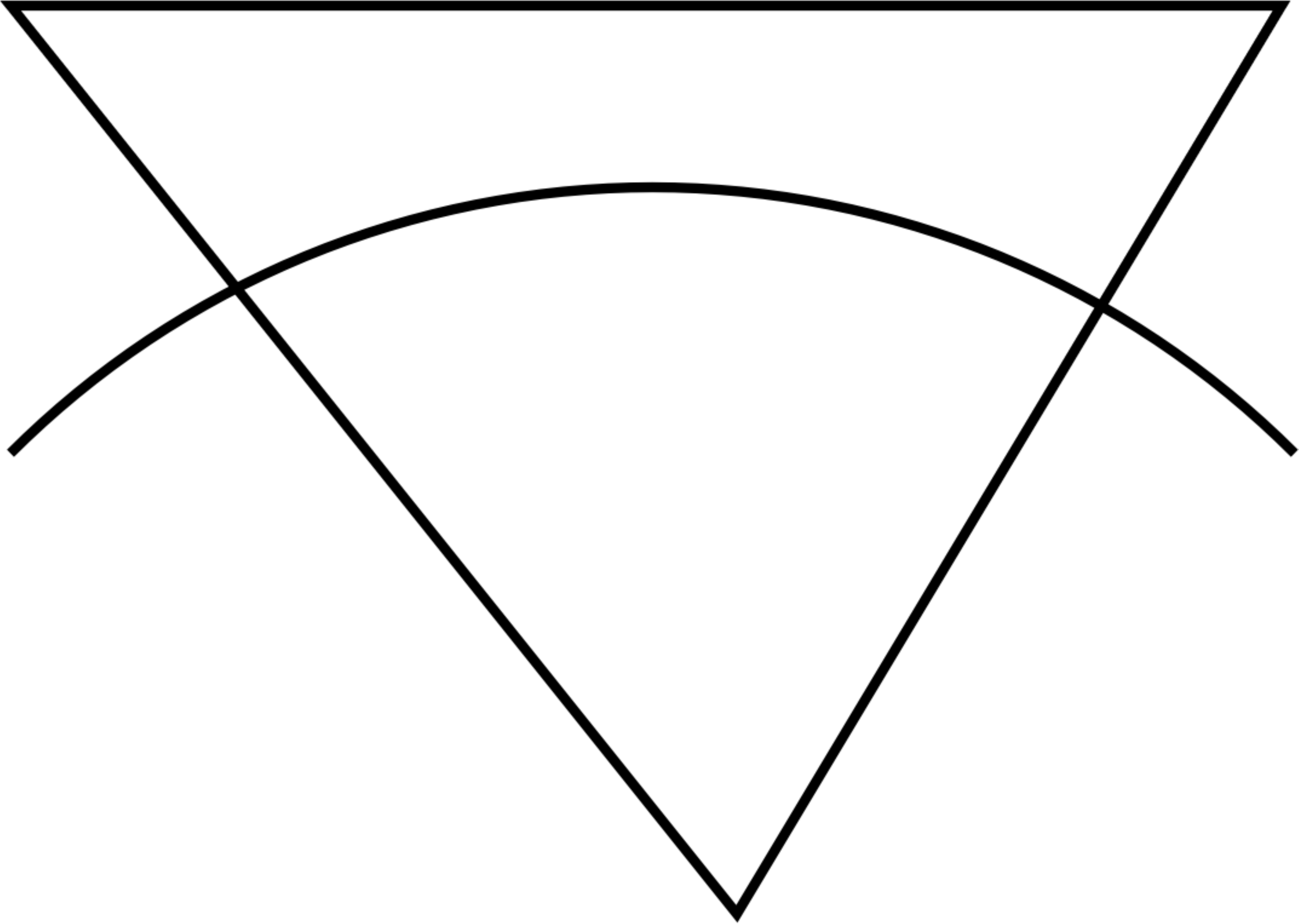} &
\includegraphics[height=.9in]{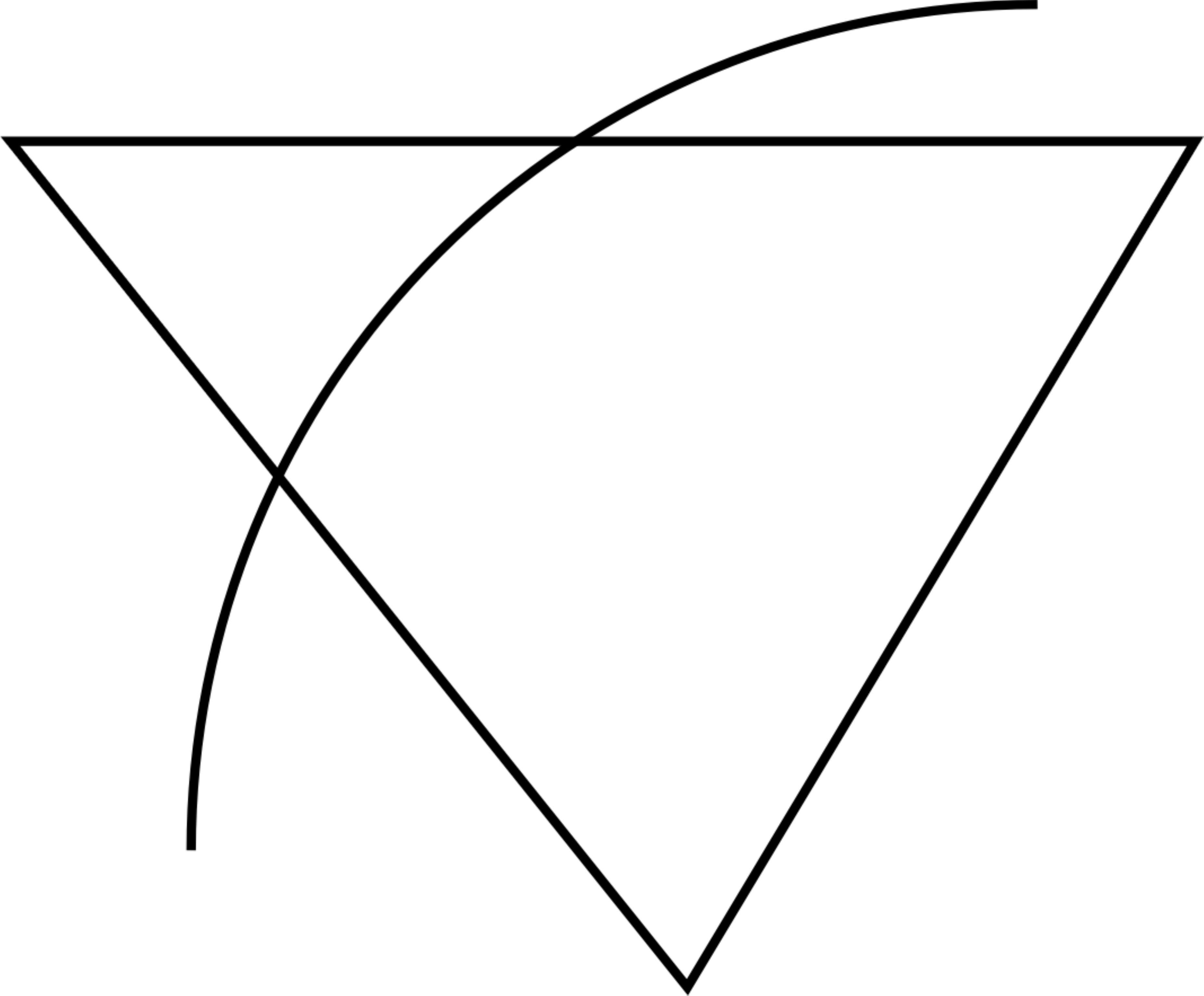} 
\\
(a) & (b) 
\\[1ex]
\includegraphics[height=.9in]{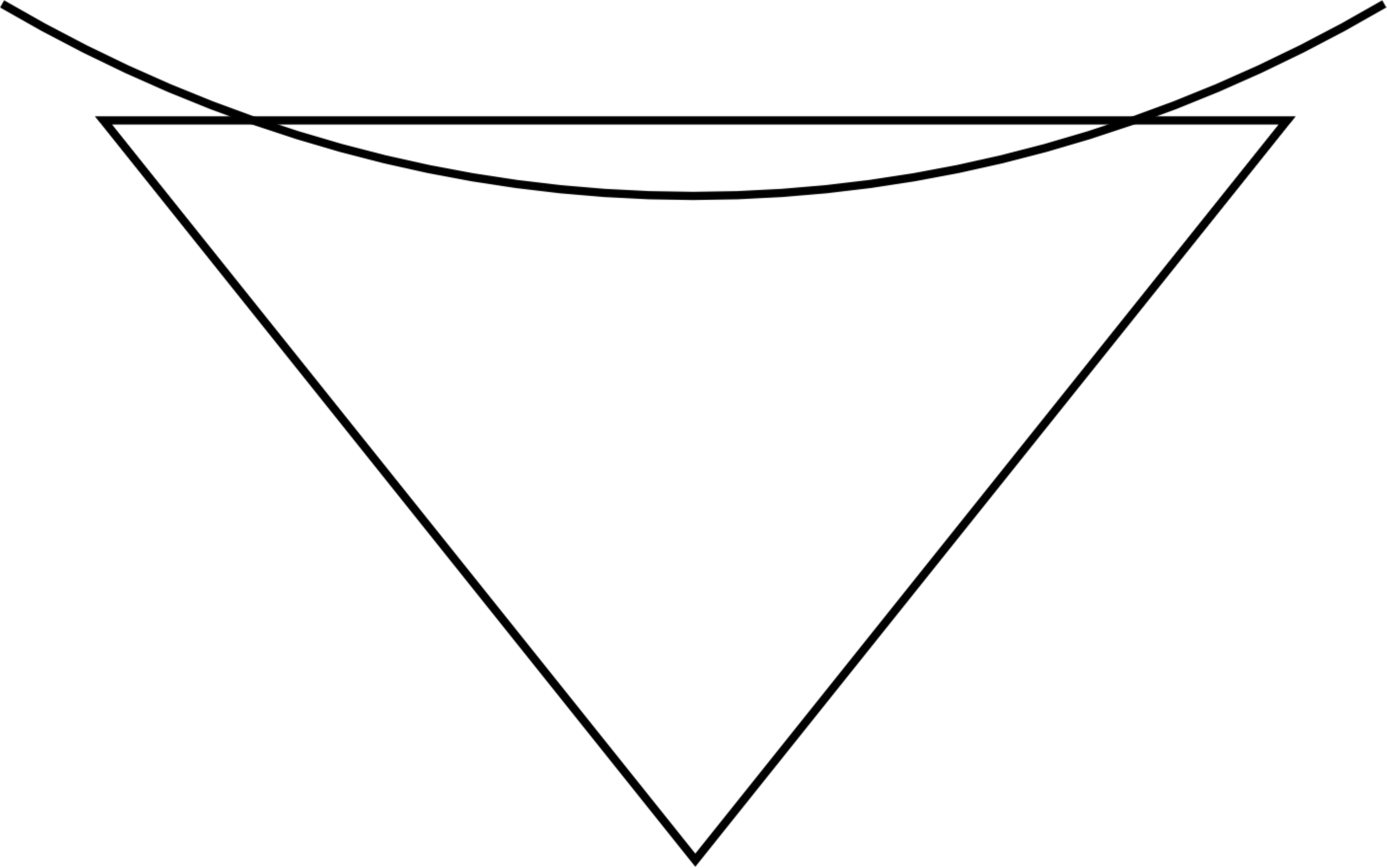} &
\includegraphics[height=.9in]{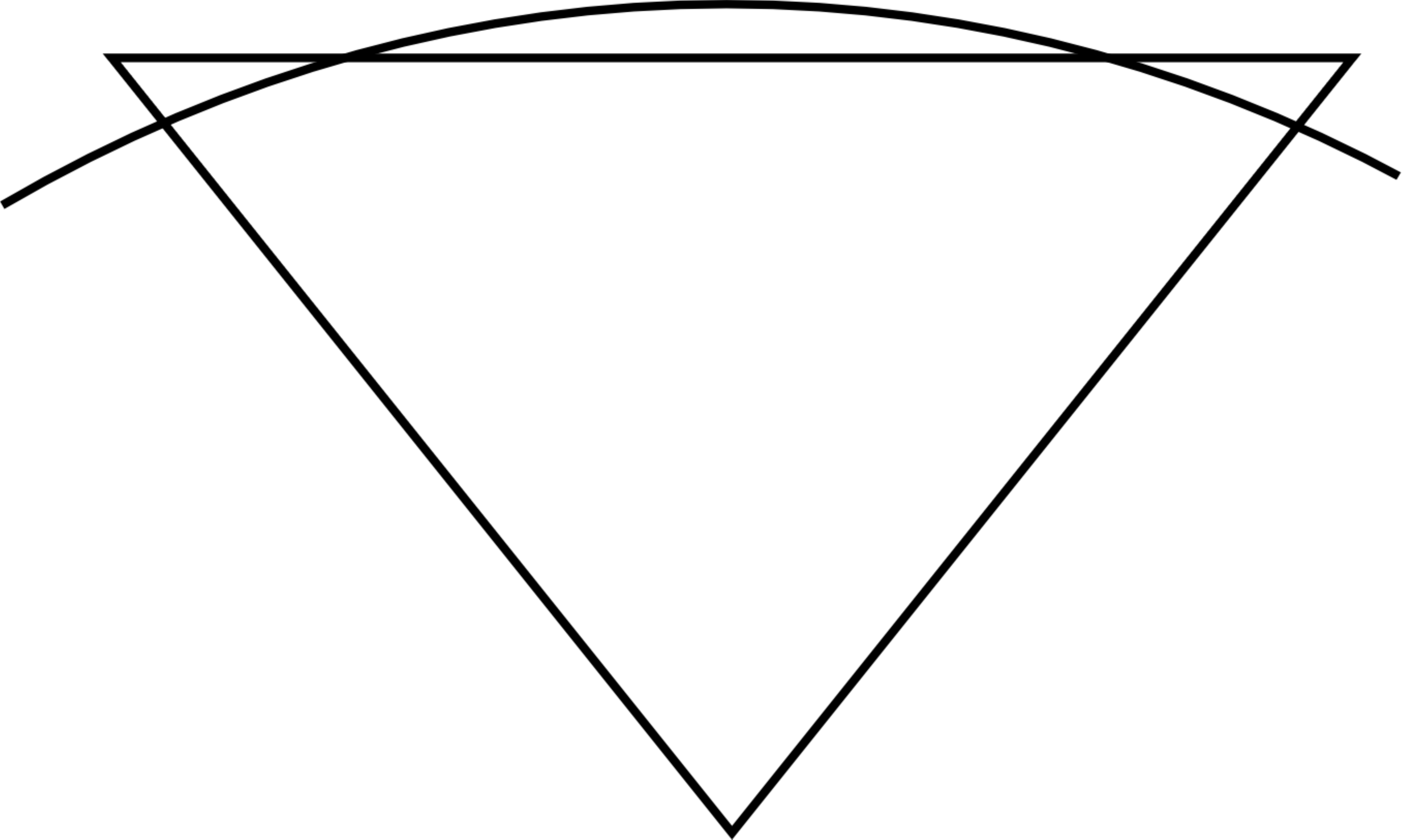}
\\
 (c) & (d)
\end{tabular} 
\caption{For (a), (b), and (c), the circle intersects the boundary of the triangle at two points whereas for (d), there are four such points. For (a) and (b), the overlap of the ball and the triangle is a three-sided, respectively four-sided, figure with one curved side. For (c), the overlap is a two-sided figure with one curved side whereas for (d), the overlap is a five-sided figure with two curved sides.}
\label{balltriangle}
\end{figure}

\subsection{\textbf{The efficient assembly of the stiffness matrix and right-hand side vector}}\label{femassemb}

In this section, we discuss the finite element assembly process for the linear system \eqref{D-Dh-weak} for the {\em approximate balls} introduced in Sec. \ref{approxballs}.

{For the sake of simplicity of exposition, we describe the assembly process for the stiffness matrix entries \eqref{bform3} and the components of right-hand side vector \eqref{lfunc3}. The assembly process for the fully-discrete system using quadrature rules and approximate balls follows along the same lines.}

The assembly of the entries $A(\phi_{j'},\phi_j)$ of the stiffness matrix and the components $F(\phi_{j'})$ of the right-hand side vector for nonlocal problems differs in several ways from that for local problems. Because the differences are substantial, in this section, we discuss, in some detail, the assembly process for nonlocal problems. Thus, the tasks in hand is to describe how to compute the entries of the stiffness matrix and the components of the right-hand side vector corresponding to the finite element discretization \eqref{weakja} of the nonlocal weak formulation \eqref{weak}. 

Of course, these tasks can be accomplished through the direct use of \eqref{bform3} and \eqref{lfunc3}. However, for the reasons we are about to remark on, an alternate approach results in a more efficient assembly process.

\hangfive~--~If $\xb\in\Omega$ is within a distance $\del$ of the boundary of $\Omega$, we have that $B_\del(\xb)=(\Omega\cap B_\del(\xb))\cup (\Omega_{\mathcal I}\cap B_\del(\xb))$ with both $\Omega\cap B_\del(\xb) \ne \emptyset$ and $\Omega_{\mathcal I}\cap B_\del(\xb) \ne \emptyset$, i.e., $\Omega\cap B_\del(\xb) \ne B_\del(\xb)$ and $\Omega_{\mathcal I}\cap B_\del(\xb) \ne B_\del(\xb)$, so that the domains of integration of all three inner integrals in \eqref{bform3} and \eqref{lfunc3} are {\em partial balls}. 

\hangfive~--~Having to define quadrature rules for partial balls certainly adds to the complexity of the stiffness matrix assembly process. For example, one is not only faced with the task of identifying intersections of the surface of the ball with finite elements, but one is also faced with the equally daunting task of identifying the intersection of finite elements and the boundary of $\Omega$ that separates the two partial balls. 

\hangfive~--~Thus, one would rather only deal, as much as possible, with integrations over {\em whole balls}, which, as seen in Sec. \ref{approxballs}, is in itself already a complex process. 

\hangfive~--~Fortunately, taking advantage of the fact that for any $\xb\in\Omega$, we have that $(\OuO)\cap B_\del(\xb) = B_\del(\xb)$, i.e. a whole ball, and also taking advantage of the equivalence \eqref{equival}, it is possible to only deal with whole balls by basing the assembly process not on \eqref{bform3} and \eqref{lfunc3}, but instead on the finite element discretization of \eqref{bform0} and \eqref{lfunc0}. 

\noindent Thus, we describe the assembly process using
\beq{bformhj}
\bal
&D(\phi_{j'},\phi_{j}) \\&\qquad= 
\sum_{k=1}^{K}\sum_{k'=1}^{K}
\int_{\mcE_k}\int_{\mcE_{k'}\cap B_\del(\xb)} \big(\phi_{j}(\xb)-\phi_{j}(\yb)\big)\big(\phi_{j'}(\xb)-\phi_{j'}(\yb)\big)\psi(\xb,\yb)d\yb d\xb
\eal
\eeq
for $j,j'=1,\ldots,J$, and
\beq{lfunchj}
G(\phi_{j'}) = \sum_{k=1}^{K_\Omega} \int_{\mcE_k} \phi_{j'}(\xb) f(\xb) d\xb \quad \mbox{for $j'=1,\ldots,J$},
\eeq
keeping in mind that the equivalence \eqref{equival} requires that $v_h(\xb)=0$ whenever $\xb\in\overline\Omega_{\mathcal I}$ so that that any term in \eqref{bformhj} involving a basis function $\phi_{j'}(\cdot)$ evaluated at any point in $\Omega_{\mathcal I}$ can be ignored, i.e., it does not contribute\footnote{We start the assembly process with a {$J_\Omega\times J_\Omega$ matrix} having all entries set to zero and a $J_\Omega$-dimensional vector having all components set to zero. Then, by ``contribute'' we mean that a computed quantity such as $D(\phi_{j'},\phi_{j})$ is, for example, added to whatever is already present in $A(\phi_{j'},\phi_{j})$ entry of the stiffness matrix.} to the stiffness matrix entry $A(\phi_{j'},\phi_{j})$. 

We reiterate that the task at hand is not to assemble the $J\times J$ matrix having entries \eqref{bformhj} and the $J$-dimensional vector having components \eqref{lfunchj}, but instead it is to use \eqref{bformhj} and \eqref{lfunchj} to compute the entries in \eqref{bform3} and the components \eqref{lfunc3}.

For $\xb\in\Omega$, the domain of integration of the inner integral in \eqref{bformhj} is the whole ball because in this case $(\OuO)\cap B_\del(\xb)=B_\del(\xb)$. On the other hand, for $\xb\in\Omega_{\mathcal I}$, the domain of integration of the inner integral is a partial ball because, in this case, $(\OuO)\cap B_\del(\xb)$ is a strict subset of $B_\del(\xb)$. This, however, does not cause a problem because points exterior to $\Omega\cap\Omega_{\mathcal I}$ are never accessed during the assembly process.
  
Thus, the remaining task is to assign the various terms appearing in \eqref{bformhj} to either contribute to the stiffness matrix entry $A(\phi_j,\phi_{j'})$ in \eqref{bform3} or to the right-hand side vector component $F(\phi_{j'})$ in \eqref{lfunc3}. It is important to note that {\em the assignment rules automatically take care of the fact that we have, in \eqref{bform3}, partial ball integrations.} These assignments are made as given in the boxed text below. We note that accounting for the contribution of $G(\phi_{j'})$ to the right-hand side vector, i.e., to the first term in \eqref{lfunc3}, is entirely identical to what is done in finite element methods for the local PDE setting so that we do not further consider this step.

\begin{table*}[ht]
\begin{center}
{
\begin{tabular}{|l|}
\hline
\multicolumn{1}{|l|}{}\\
A. If $\xb\in\mcE_{k}\in\mcT_{h,\Omega}=\Omega$ and $\yb\in\mcE_{k'}\cap B_\del(\xb)\in\mcT_{h,\Omega}=\Omega$ 
\\
\qquad\qquad\qquad[{\em occurs for all $\xb\in\Omega$}]
\\
\qquad$\Rightarrow$\quad $\phi_{j'}(\xb)\ne0$, $\phi_{j'}(\yb)\ne0$  
\\ 
\qquad$\Rightarrow$\quad each of the $\phi_j(\xb)\phi_{j'}(\xb)$,\, $\phi_j(\xb)\phi_{j'}(\yb)$,\, $\phi_j(\yb)\phi_{j'}(\xb)$, and $\phi_j(\yb)\phi_{j'}(\yb)$
\\\qquad\qquad
 terms in \eqref{bformhj} makes a nonzero
contribution to the stiffness matrix
\\\qquad\qquad
 entry $A(\phi_{j'},\phi_j)$.
\\[2ex]
B. If $\xb\in\mcE_{k}\in\mcT_{h,\Omega}=\Omega$ and $\yb\in\mcE_{k'}\cap B_\del(\xb)\in\mcT_{h,\Omega_{\mathcal I}}=\Omega_{\mathcal I}$  
\\\qquad\qquad\qquad
[{\em occurs only if the distance from $\xb\in\Omega$}
\\\qquad\qquad\qquad\quad
{\em to the boundary of $\Omega$ is less than $\del$}]
\\\qquad
$\Rightarrow$\quad $\phi_{j'}(\xb)\ne0$, $\phi_{j'}(\yb)=0$ \\
\qquad $\Rightarrow$\quad the $\phi_j(\xb)\phi_{j'}(\xb)$ term in \eqref{bformhj} makes a nonzero contribution to the
\\\qquad\qquad
 stiffness matrix entry $A(\phi_{j'},\phi_j)$;
\\ 
\qquad $\Rightarrow$\quad the $\phi_j(\yb)\phi_{j'}(\xb)$ term in \eqref{bformhj} makes nonzero contribution to the
\\\qquad\qquad
 right-hand side vector component $F(\phi_{j'})$.
\\[2ex]
C. If $\xb\in\mcE_{k}\in\mcT_{h,\Omega_{\mathcal I}}=\Omega_{\mathcal I}$ and $\yb\in\mcE_{k'}\cap B_\del(\xb)$
\\\qquad\qquad\qquad
[{\em occurs for all $\xb\in\Omega_{\mathcal I}$ but only if the distance from $\yb\in\Omega$}
\\\qquad\qquad\qquad\quad
{\em to the boundary of $\Omega$ is less than $\del$}]
\\
\qquad$\Rightarrow$\quad $\phi_{j'}(\xb)=0$, $\phi_{j'}(\yb)\ne0$ 
\\ 
\qquad $\Rightarrow$\quad the $\phi_j(\yb)\phi_{j'}(\yb)$ term in \eqref{bformhj} makes nonzero contribution to 
\\\qquad\qquad
the stiffness matrix entry $A(\phi_{j'},\phi_j)$;
\\ 
\qquad $\Rightarrow$\quad the $\phi_j(\xb)\phi_{j'}(\yb)=g(\xb)\phi_{j'}(\yb)$ term in \eqref{bformhj} makes a nonzero 
\\\qquad\qquad
contribution to the right-hand side vector component $F(\phi_{j'})$.
\\[2ex]
{All other combinations of $\xb$ and $\yb$ and domains result in zero contributions.}
\\[2ex]
\hline
\end{tabular}
}
\end{center}
\end{table*}

In Fig. \ref{fintint}, the white and orange regions are part of $\Omega$ and the yellow and magenta regions are part of $\Omega_{\mathcal I}$. Choice (A) in the box involves a whole ball lying completely within $\Omega$ if the distance from $\xb\in\Omega$ to the boundary of $\Omega$ is larger than $\delta$; see Fig. \ref{fintint}a for an illustration. On the other hand, if the distance from $\xb\in\Omega$ to the boundary of $\Omega$ is smaller than $\delta$ we again have a whole ball but Case (A) applies only to the partial ball lying within $\Omega$, as illustrated by the orange partial ball in Fig. \ref{fintint}b, and Case (B) applies to the partial ball lying within $\Omega_{\mathcal I}$, as illustrated by the magenta partial ball in Fig. \ref{fintint}b.
However, one does not have to explicitly deal with the partial balls; one simply cycles through all the triangles that intersect with the whole ball and let the assignment rules (A) and (B) automatically take care of which terms are assigned to make contributions to either the stiffness matrix or the right-hand side vector. 
Choice (C) involves three partial balls, i.e., the orange and magenta regions depicted in Fig. \ref{fintint}c and the missing part of the ball that lies outside of $\Omega\cup\Omega_{\mathcal I}$. As was the case for Cases (A) and (B), the colored regions need not be explicitly differentiated because the assignment rules (C) automatically take care of which terms are assigned to make contributions to either the stiffness matrix or the right-hand side vector. The part of the ball that lies outside of $\Omega\cup\Omega_{\mathcal I}$ is also automatically ``taken care of'' because at no step in the assembly process are points in that partial ball accessed. 

\begin{figure}[ht]
\centering
\begin{tabular}{ccc}
\includegraphics[height=1.5in]{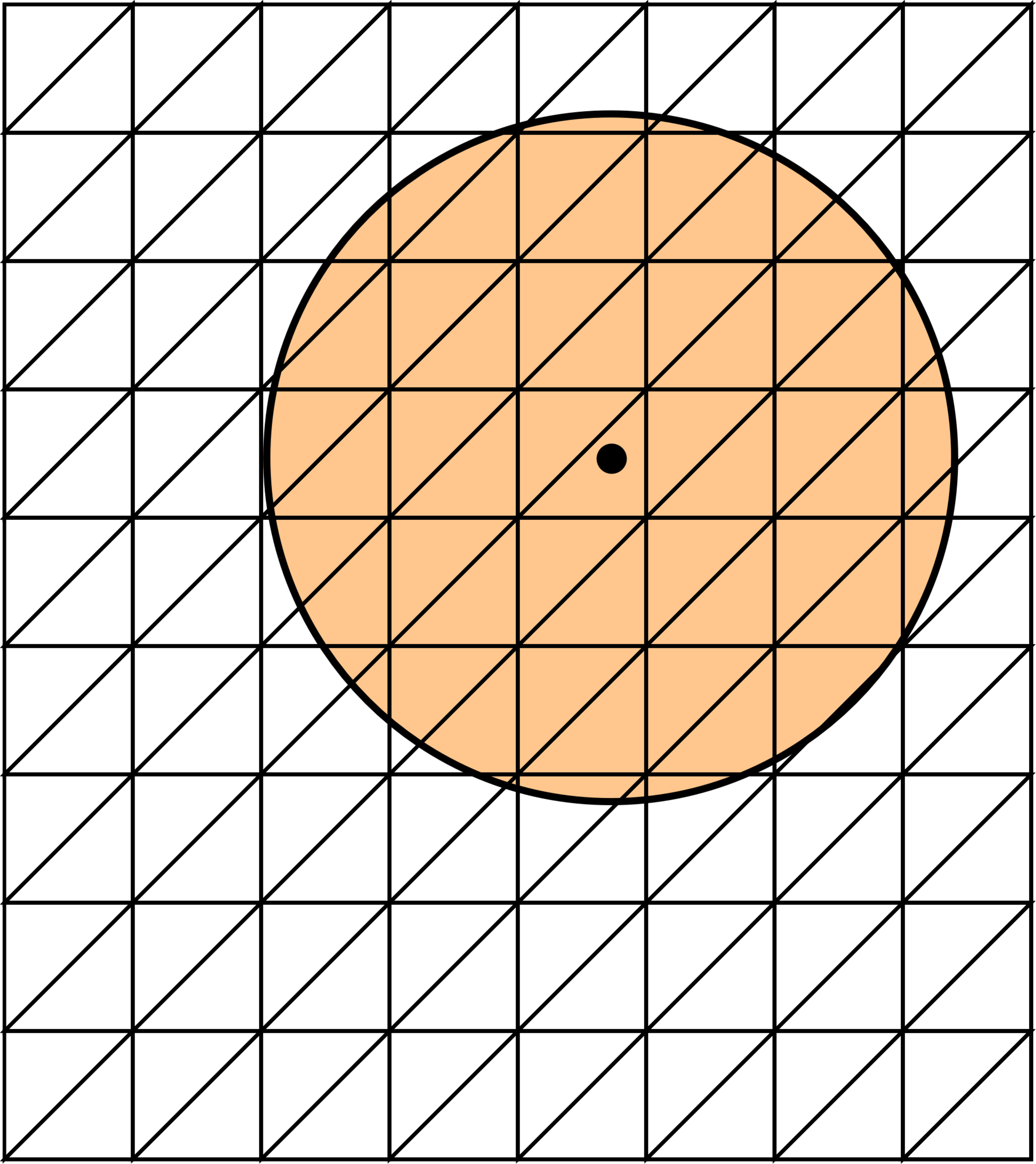} &
\includegraphics[height=1.5in]{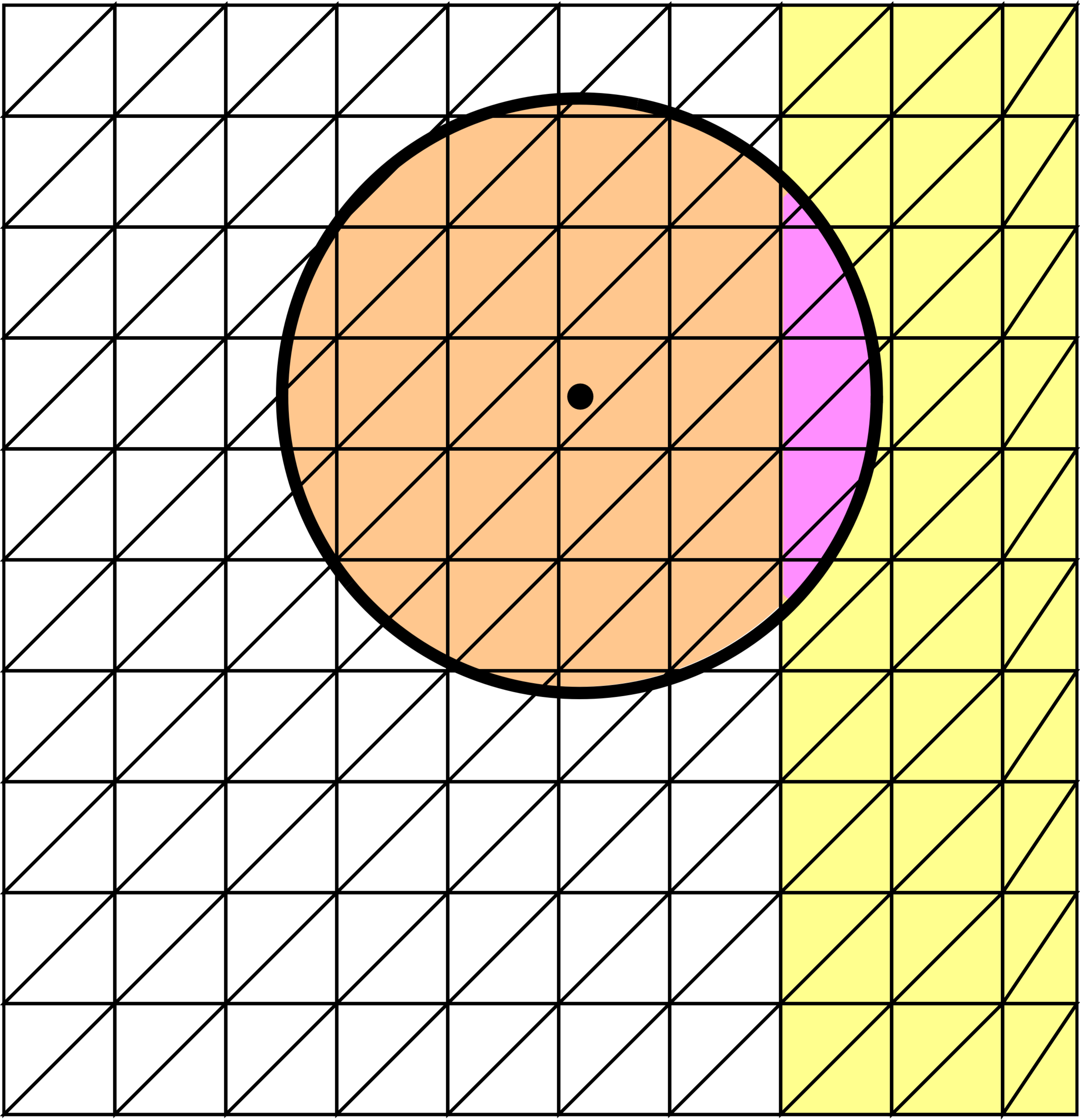} &
\includegraphics[height=1.5in]{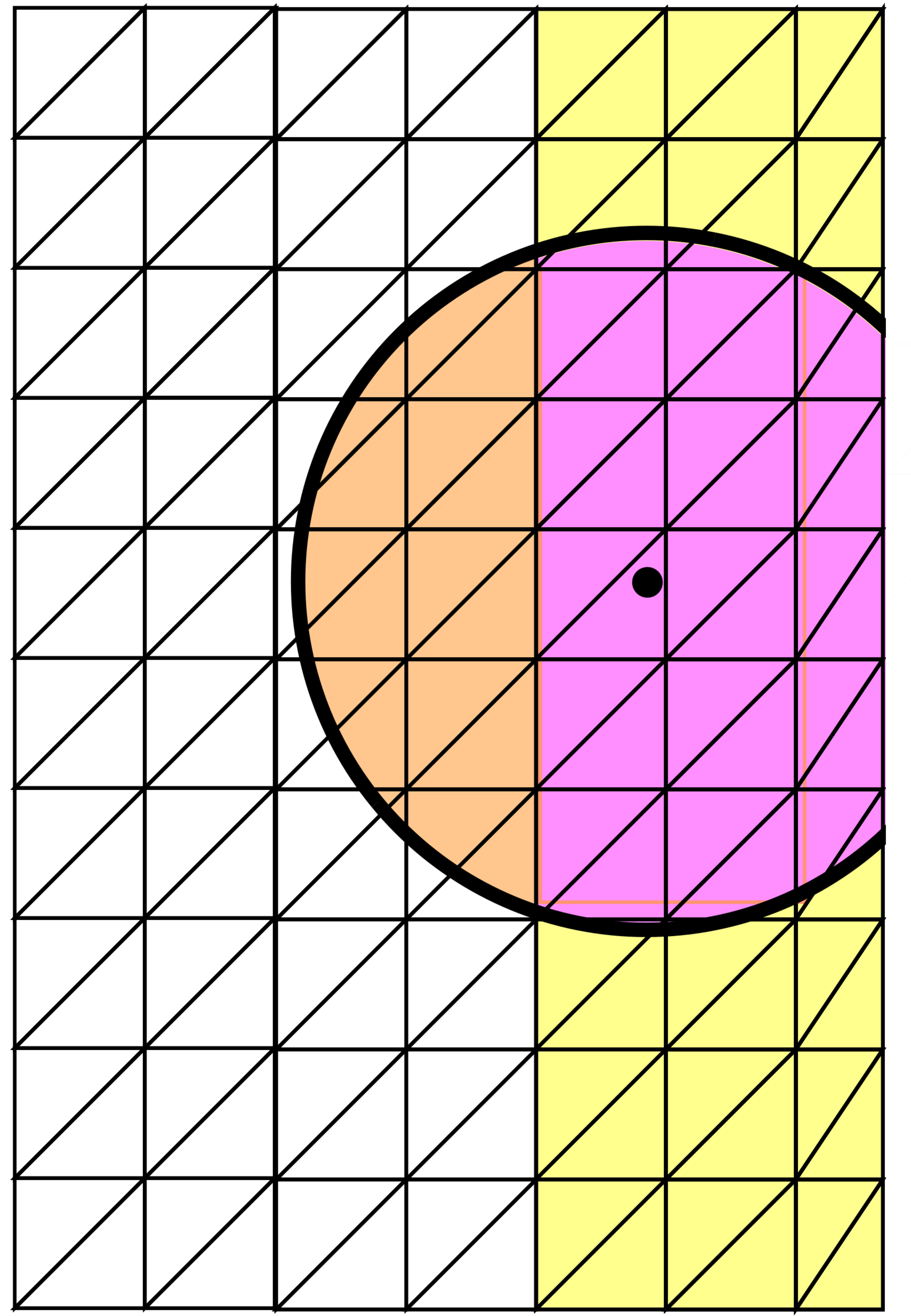}
\\
(a) & (b) & (c)
\end{tabular} 
\caption{Three of the possible configurations for balls $B_\del(\xb)$ relative to the position of their center and the domains $\Omega$ (the white and orange regions) and $\Omega_{\mathcal I}$ (the yellow and magenta regions). }
\label{fintint}
\end{figure} 

Note that the assignment recipe (A-B-C) is the analog of what is done in the local case for which, in the finite element assembly process, terms that correspond to boundary nodes are moved to the right-hand side whereas terms that involve interior nodes contribute to the stiffness matrix. 

Some remarks are in order.

\vskip5pt

{\em\textbf{Reduced sparsity due to nonlocality.}} 
During the finite element assembly process, one is faced with having to compute terms that contribute to the stiffness matrix, terms such as $A(\phi_j,\phi_{j'})$. In stark contrast with local models, for nonlocal models \eqref{bform3} implies that even a pair of basis functions $\{\phi_j,\phi_{j'}\}$ having non-overlapping support may \textit{interact}, i.e., may yield a nonzero entry in the stiffness matrix because interactions occur over a distance. Consequently, compared to that for finite element discretizations of local models that use the same grid and same finite element spaces, the nonlocal stiffness matrix is more densely populated which is the discrete realization of nonlocality. The source of reduced sparsity is illustrated in the left plot of Fig. \ref{sparsification}. In that figure, the triangles represent a portion of a domain $\Omega$. The filled circle is a point in the blue outer integral triangle and is the center of the orange ball. The open circle at a vertex of the blue triangle is a node $\widetilde\xb_j$. The blue and red patch of triangles represent the support of the basis function  $\phi_j(\xb)$ corresponding to that node. The pink patch of triangles represent the support of a basis function $\phi_{j'}(\xb)$ corresponding to the open circle node $\widetilde\xb_{j'}$ in that patch. Because both the blue/red and pink patches overlap with the orange ball, the pair $\{\phi_j(\xb),\phi_{j'}(\xb)\}$ makes a nonzero contribution to the stiffness matrix entry $A(\phi_j,\phi_{j'})$. The number of nonzero entries depends on the relations between the size of the interaction radius $\del$, the size of the domain $\Omega$, and the grid size. The reduced sparsity compared to that for stiffness matrices for local models results in greater assembly costs, in the need for additional memory storage, and in greater solution costs. The use of appropriate solvers for the linear systems is of fundamental importance. See, e.g., {Refs. \refcite{AkPa11,dt-burak13,DElia-ACTA-2020,vollman,WaTi12}} for further discussions about this issue.

\begin{figure}[ht]
\centerline{
\includegraphics[height=1.5in]{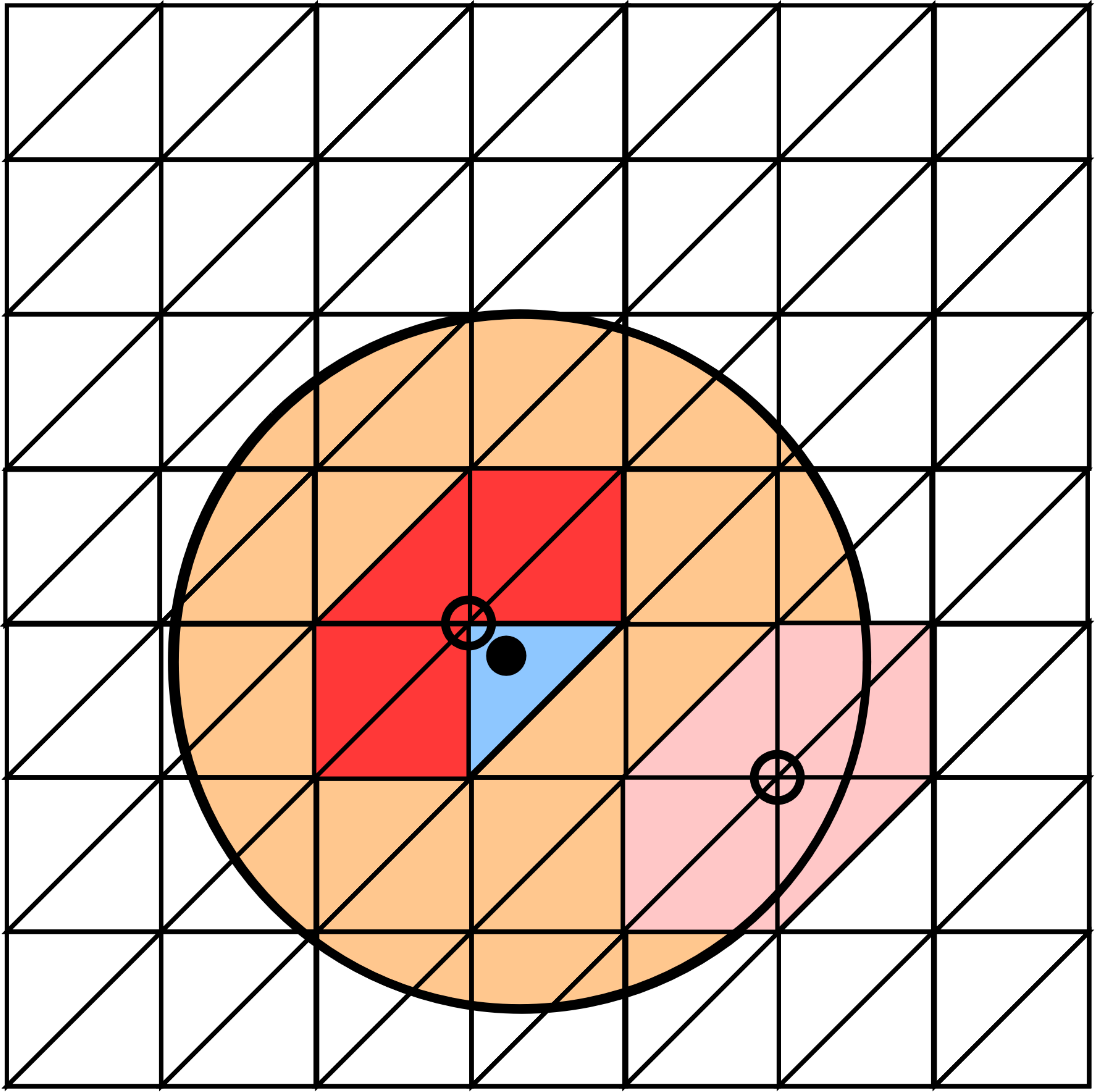}
\quad
\includegraphics[height=1.5in]{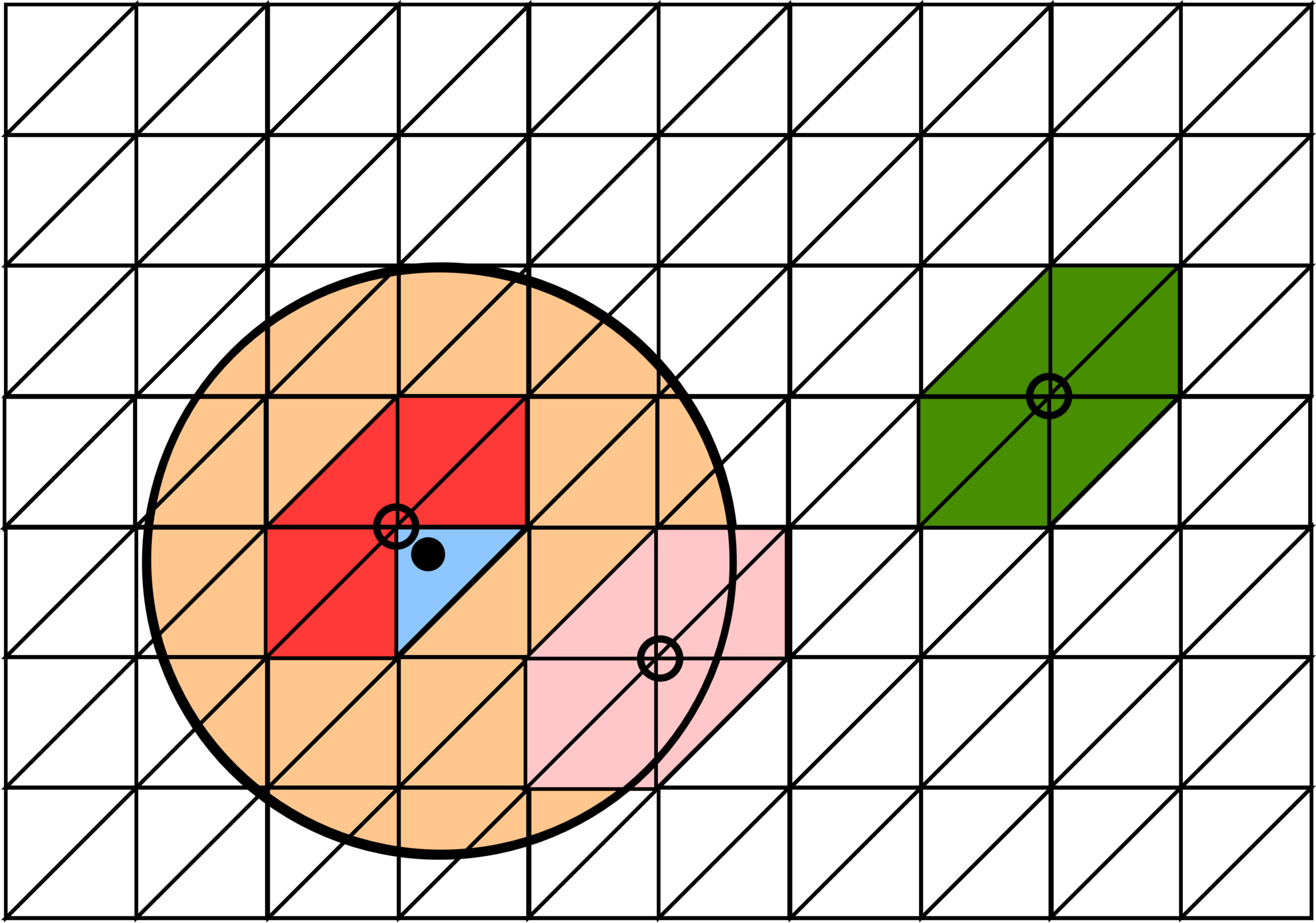}
} 
\caption{Left: the pair of blue/red and pink patches represent the support of two basis functions that make a nonzero contribution to the stiffness matrix. Right: the pair of blue/red and green patches do not make such a contribution.}
\label{sparsification}
\end{figure} 

{\em\textbf{Sparsification due to finite horizons.}}
It is clear from \eqref{bform3} that, for the nonlocal case, two finite element basis functions $\phi_j(\xb)$ and $\phi_{j'}(\xb)$ interact \textit{only if} both of their supports overlap with $B_\del(\xb)$.  Thus, if the diameter $2\del$ of the ball $B_\del(\xb)$ is larger than the diameter of $\Omega$, then the nonlocal stiffness matrix is a full matrix. On the other hand, if $\Omega\cap B_\del(\xb)\ne \Omega$, i.e., if the diameter $2\del$ of the ball $B_\del(\xb)$ is smaller than the diameter of $\Omega$, some entries in the stiffness matrix vanish. This situation is illustrated in the right plot of Fig. \ref{sparsification}. 
On the other hand, the green patch that is the support of a basis function $\phi_{j'}(\xb)$ that now corresponds to the open circle node in that patch, does not overlap with the $B_\del(\xb)$ so that, paired with $\phi_j(\xb)$, it does not contribute to the stiffness matrix. This leads to the sparsification we have been alluding to.

\subsection{\textbf{Fully-discrete weak formulation}}\label{fuldis}

The fully-discrete weak formulations we consider can involve the possible application of three approximations to the linear system \eqref{weakja}.

\hangseven--~{An approximate polytopial ball $B_{\del,h}(\xb)$ is used to approximate the ``exact'' ball $B_{\del}(\xb)$; see Sec. \ref{approxballs}.} 

\hangseven--~{A quadrature rule is used to approximate inner integrals; see Sec. \ref{apinner}.} 

\hangseven--~{A quadrature rule is used to approximate outer integrals; see Sec. \ref{apouter}.}

\noindent To define a fully-discrete stiffness matrix and right hand-side vector, we need to have ready the following mise en place about which detailed considerations are given in Sections \ref{approxballs} to \ref{apouter}.

\hangseven--~For each element $\mcE_k\in\mcT_{h,\Omega}$, the outer integrals in \eqref{bform3} and \eqref{lfunc3} are approximated using a quadrature rule with quadrature points $\xb^{outer}_{k,q}$ and corresponding weights $w^{outer}_{k,q}$, $q=1,\ldots,Q_k^{outer}$. 

\hangseven--~The approximate balls $B_{\del,h}(\xb^{outer}_{k,q})$ centered at each of the quadrature points $\xb^{outer}_{k,q}$ of the outer integral are subdivided into a set of subdomains ${\widetilde\mcT}_{h,\delta,k,q}$. 

\hangseven--~The integrals over each subdomain ${\widetilde\mcT}_{k'}\in{\widetilde\mcT}_{h,\delta,k,q}$ are approximated using a quadrature rule with quadrature points $\xb^{inner}_{k',q'}$ and corresponding weights $w^{inner}_{k',q'}$, $q'=1,\ldots,Q_k^{inner}$.

\noindent Then, applying the three approximations and three ingredients to \eqref{bform3} and \eqref{lfunc3} leads to the discrete approximation of the linear system \eqref{weakja} given by
\beq{weakjqqha}
\sum_{j=1}^{J_\Omega} A_{qh}(\phi_{j'},\phi_{j})U_{j;qh} = F_{qh}(\phi_{j'}) \quad  \mbox{for $j'=1,\ldots,J_\Omega$},
\eeq
where the fully-discrete stiffness matrix entries are given by
\beq{bformjqqh}
\bal
& A_{qh}(\phi_{j'},\phi_{j}) =
\\& \sum_{\mcE_k\in\mcT_{h,\Omega}} \sum_{q=1}^{Q^{outer}} {w^{outer}_{k,q}}
\sum_{ \widetilde\mcE_{k'}\in{\widetilde\mcT}_{h,\delta,k,q}} \,
\sum_{q'=1}^{Q^{inner}} {w^{inner}_{k',q'}}   
 \big(\phi_j(\yb^{inner}_{k',q'})-\phi_j(\xb^{outer}_{k,q})\big)
 \\&\qquad\qquad\qquad
 \times\,\,\big(\phi_{j'}(\yb^{inner}_{k',q'})-\phi_{j'}(\xb^{outer}_{k,q})\big)\psi(\xb^{outer}_{k,q},\yb^{inner}_{k',q'})
\\&\qquad+2
\sum_{\mcE_k\in\mcT_{h,\Omega}} \sum_{q=1}^{Q^{outer}} {w^{outer}_{k,q}} 
\phi_j(\xb^{outer}_{k,q})\phi_{j'}(\xb^{outer}_{k,q})
\\&\qquad\qquad\qquad
\times\,\sum_{ \widetilde\mcE_{k'}\in{\widetilde\mcT}_{h,\delta,k,q}}\, \sum_{q'=1}^{Q^{inner}}\, {w^{inner}_{k',q'}} 
 \psi(\xb^{outer}_{k,q},\yb^{inner}_{k',q'})
\eal
\eeq
for $j,j'=1,\ldots,J_\Omega$ and the fully-discrete right-hand side components are given by
\beq{lfunchq}
\bal
F_{qh}(\phi_{j'}) &= \sum_{\mcE_k\in\mcT_{h,\Omega}} \sum_{q=1}^{Q^{outer}} {w^{outer}_{k,q}}  \phi_{j'}({\xb^{outer}_{k,q}}) 
\Big(f({\xb^{outer}_{k,q}})
\\& +   
\sum_{ \widetilde\mcE_{k'}\in{\widetilde\mcT}_{h,\delta,k,q}} \sum_{q'=1}^{Q^{inner}}\, \, {w^{inner}_{k',q'}} g(\yb^{inner}_{k',q'})
 \psi(\xb^{outer}_{k,q},\yb^{inner}_{k',q'})  \Big)
\eal
\eeq
for $j'=1,\ldots,J_\Omega$. The assembly of the coefficient matrix entries \eqref{bformjqqh} and right-hand side vector components  \eqref{lfunchq} can be accomplished by using the (A-B-C) recipe of Sec. \ref{femassemb} with $B_{\del}(\xb)$ replaced by the approximate ball $B_{\del,h}(\xb)$ and, of course, with \eqref{bformhj} and \eqref{lfunchj} approximated by a fully-discrete approximation in much the same way as the pair \eqref{bform3} and \eqref{lfunc3} was approximated by the fully-discrete pair \eqref{bformjqqh} and \eqref{lfunchq}. 
 
\section{Numerical illustrations}\label{numerics}

We consider the nonlocal problem \eqref{nonlocalsetup} on the domain $\Omega = (0,1)^2$ with a constant kernel $\g(\xb,\yb) = \frac{4}{\pi\delta^4}{\mathcal X}_{B_\delta(\xb)}(\yb)$ with $\delta = 0.1$. Of course, for this kernel, the kernel function $\psi(\xb,\yb)=\frac{4}{\pi\delta^4}$ is integrable and translationally invariant so that the error estimate \eqref{eq:FEMaccuracy} holds with $m=1$ for piecewise-linear finite element basis functions. The scaling constant $\frac{4}{\pi\delta^4}$ guarantees that $\Lu u= \Delta u$ for polynomials $u$ with order up to three; see, e.g., Ref. \refcite{Vollmann2019}. We make use of the manufactured solution  
$$
u(\xb) = x_1^2x_2 + x_2^2
$$ 
for which we obtain the corresponding source term $f(\xb)= -\Delta u(\xb) = -\Lu u=-2(x_2+1)$ for $\xb \in \Omega$ and nonlocal Dirichlet volume constraint data $g(\xb) = u(\xb)$ for $\xb\in \Omega_{\mathcal I}$.
 
We use piecewise-linear finite element basis functions on triangular grids and report on the convergence rates of the finite element approximation $u_h$ to the given exact solution $u$ with respect to the $L^2$-norm on $\Omega$. Examples of the types of grids used in the numerical illustrations are given in Fig. \ref{fig:nugrid}.

\begin{figure}[ht]
\centering
\begin{tabular}{cc}
\includegraphics[height=0.45\textwidth]{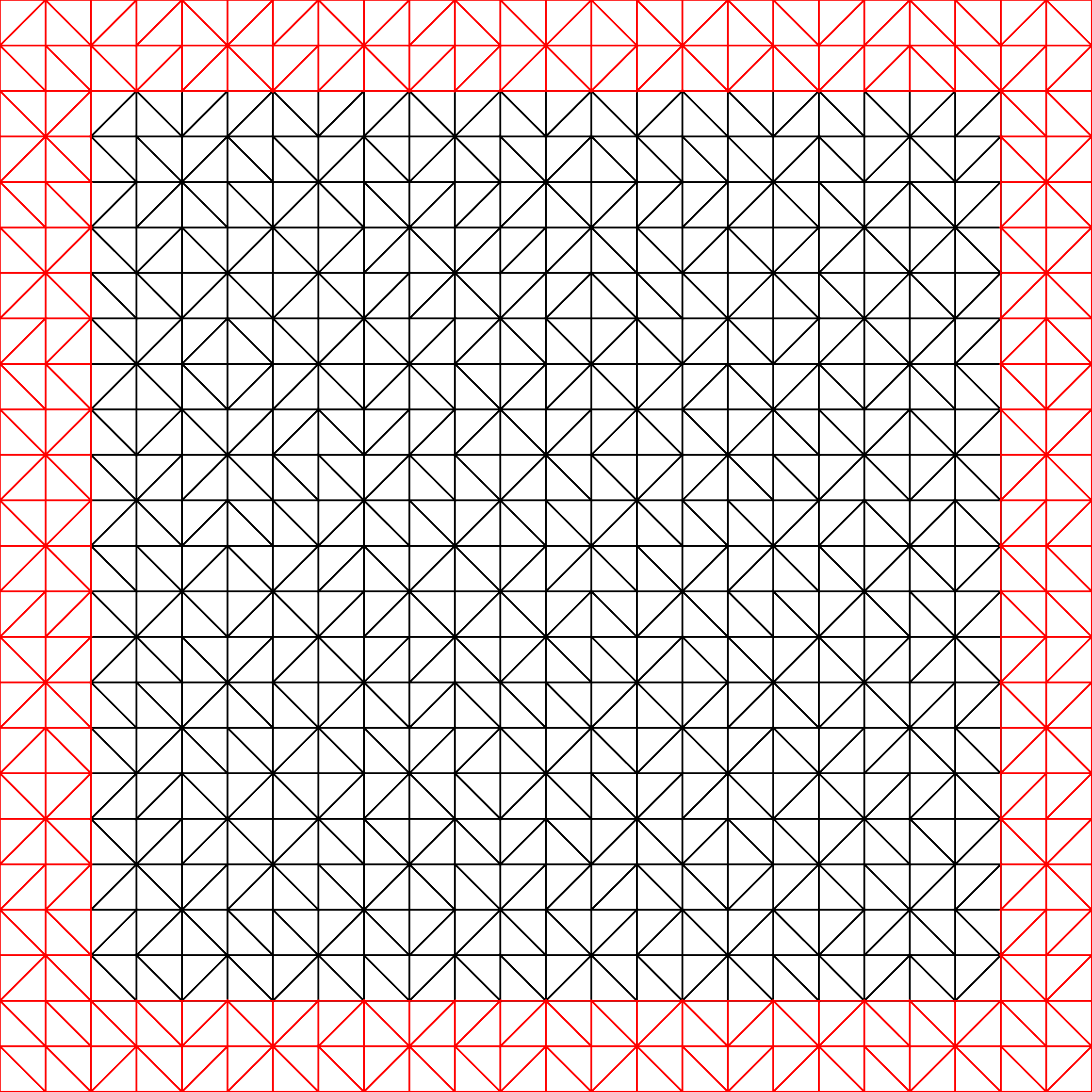}
&
\includegraphics[height=0.45\textwidth]{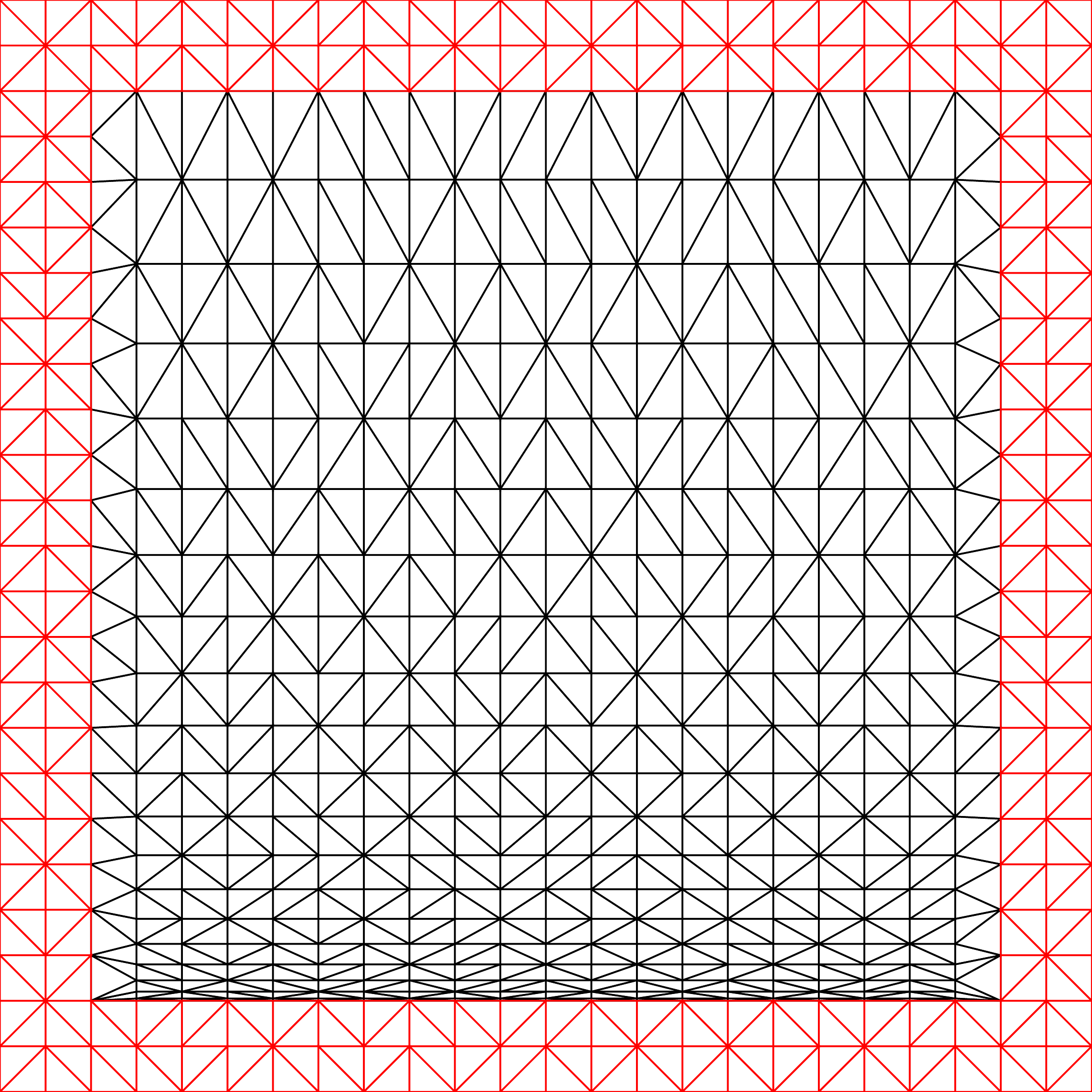}
\\
(a) & (b) 
\\[1ex]
\includegraphics[height=0.45\textwidth]{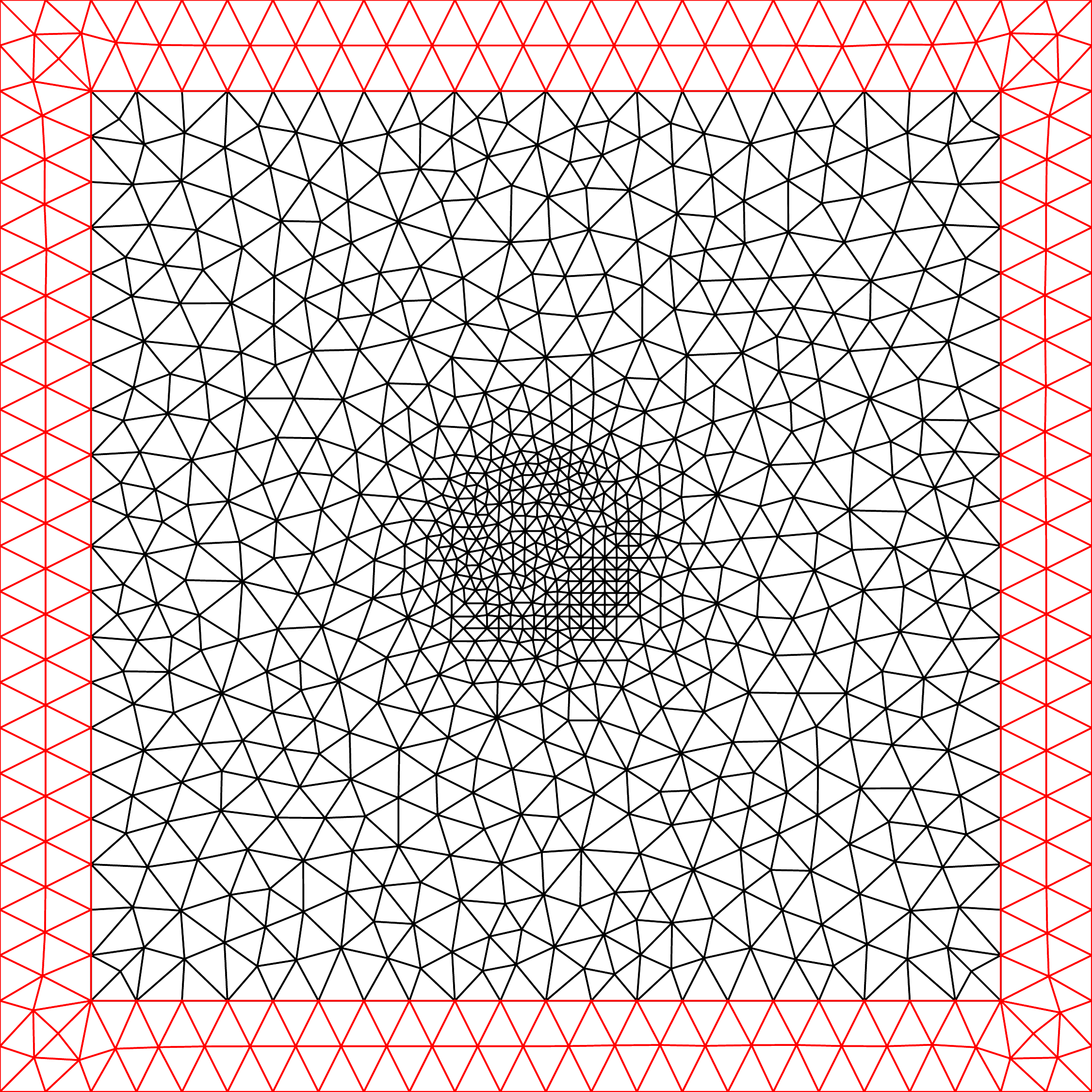}
&
\includegraphics[height=0.45\textwidth]{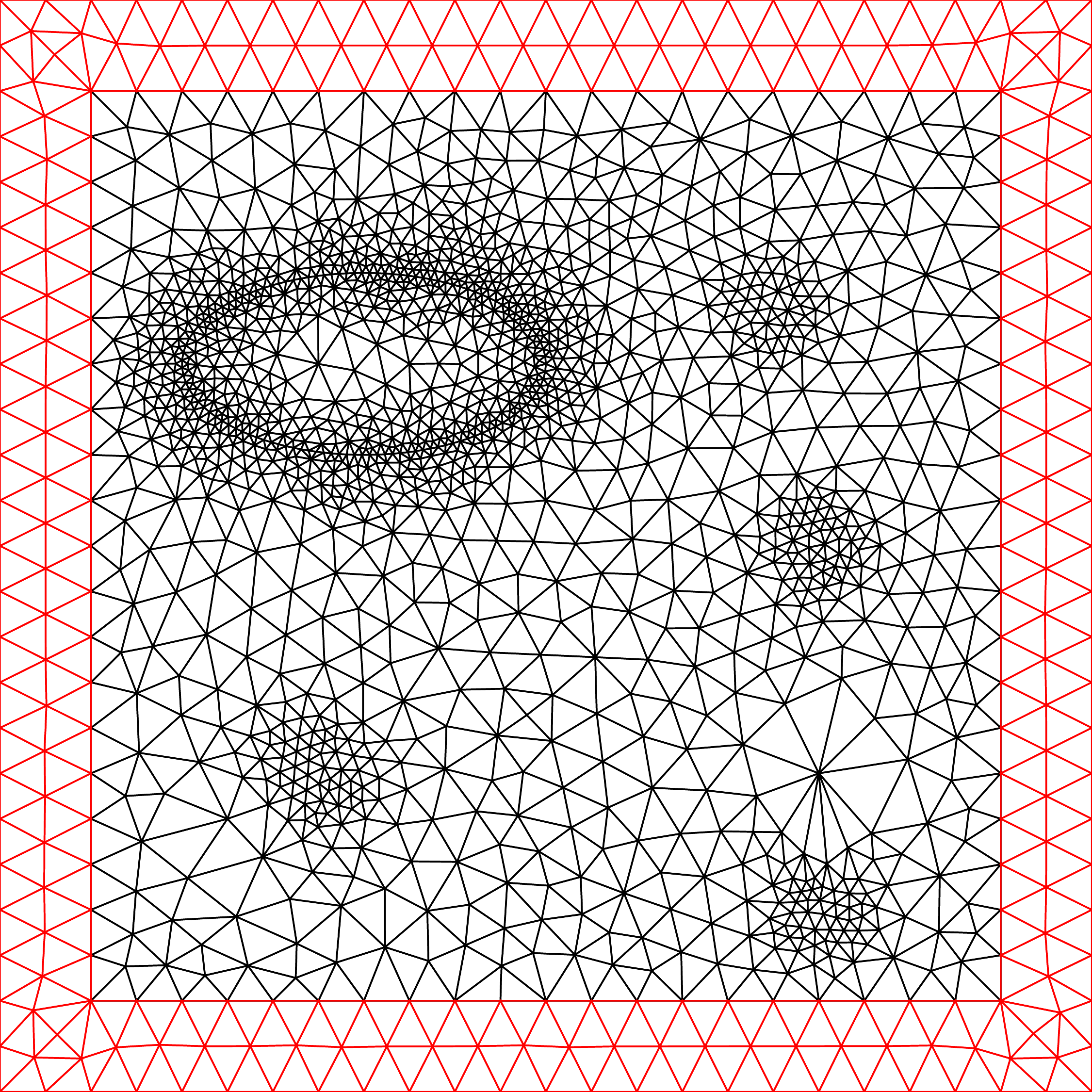} 
\\
(c) & (d) 
\end{tabular}
\caption{
(a) The uniform grid used for experiments (1) to (8). %
(b) A Cartesian but nonuniform grid used for experiments (9) and (10). It is obtained from the uniform grid (a) by applying the transformation $(x_1, x_2) \mapsto (x_1, x_2^2)$ to the interior vertices of $\Omega$. %
(c) A nonuniform and non-Cartesian grid with smooth element size transition used for experiments (11) and (12). %
(d) A highly nonuniform grid with abrupt changes in the element size used for experiments (13) and (14). %
The meshes (c) and (d) have been generated with gmsh (http://gmsh.info/). Note that nested grid refinement is used for the meshes of type (a) and (b), whereas this is not the case for the meshes of type (c) and (d).%
}%
\label{fig:nugrid}
\end{figure}

We apply the approximations for the inner and outer integrals as described in Sections \ref{apinner} and \ref{apouter}. More precisely, for inner integrals, we use a three-point, precision-two symmetric Gaussian rule for the finite element triangles and subtriangles resulting from subdividing polygonal intersection regions  and we use use a one-point centroid rule for circular caps. For outer integral triangles, in Case 1 we use a four-point precision-three symmetric Gaussian rule (see Sec. \ref{subsec:case1}) whereas for Case 2, we use the seven-point, precision-three quadrature rule introduced in Fig. \ref{fig:7points} (see Sec. \ref{subsec:case2}). In order to identify these two cases for a given pair of outer and inner integral triangles $(\mcE_k,\mcE_{k'})$, we use the following approximate criterion: we apply the four-point Gaussian rule if $|\xb^{barycenter}_k-\xb^{barycenter}_{k'}| < \delta - h$, where $h>0$ denotes the largest diameter of all finite element triangles; otherwise the aforementioned seven-point rule is used. This criterion is not sharp in the sense that we may be applying the seven-point rule to pairs of triangles for which the issues of Sec. \ref{subsec:case2} do not arise.

\subsection{Uniform grid results}
	
As can be seen in Table \ref{num:approx_inner} and Fig. \ref{fig:loglogplots}, we observe second-order convergence rates for the \textit{exactcaps} ball, as well as the ball approximations $\{nocaps,approxcaps\}$. In contrast, the $barycenter$ ball approximation produces rather erratic rates which is due to the fact that the integrand of the outer integral is discontinuous for certain pairs of outer and inner integral triangles $(\mcE_k,\mcE_{k'})$ (see Sec. \ref{subsec:case2}) so that a quadrature rule for polynomials results in inaccurate approximations. The $overlap$ approximation yields nearly first-order rates. Furthermore, the \textit{exactcaps} ball (quadrature rule for caps) and the approximation $approxcaps$ (one triangle per cap) have comparable absolute errors. Due to the decreasing approximation quality we observe higher absolute errors for the ball approximations $\{nocaps, barycenter, overlap\}$. 

For uniform grids, we also provide the relative computational time needed to assemble the respective nonlocal stiffness matrices. Therefore, for such grids, we compare all computation times relative to the largest one across results for the different ball approximations and over the grid sizes used, thus providing comparable insights into the computational effort required by the use of different ball approximations. The costliest computation was for the $exactcaps$ case with the finest grid size $0.00625$, so that, e.g., the computational cost for the $approxcaps$ case with a grid size $0.0125$ was $8.36\%$ of the highest computational cost. We observe that the determination of the centroid and the area of a circular cap requires similar steps as those to approximate the cap by a single triangle, thus computation times are comparable. Because the $nocaps$ variant is a subroutine of the $approxcaps$ variant, we observe lower costs for the latter variant, although the savings are small. Also, because the \textit{barycenter} and \textit{overlap} methods do not require the computation of intersections, they are even cheaper.

\begin{table}[ht]
\centering	
\footnotesize
\begin{tabular}{ r| ccr | ccr } 
\multicolumn{1}{c }{}& \multicolumn{3}{c }{(1) \em exactcaps} &\multicolumn{3}{c }{(2) \em  approxcaps}  \\ 		
\noalign{\smallskip}\hline\noalign{\smallskip}
$h$   & $\|u-u_h\|_{L^2}$ & rate & time [\%] & $\|u-u_h\|_{L^2}$ & rate & time [\%]\\
\noalign{\smallskip}\hline\noalign{\smallskip}
        0.1 & 9.01e-03          & -    & 0.01     & 3.78e-03          & -    & 0.01    \\
        0.05 & 1.58e-03          & 2.51 & 0.12    & 5.84e-04          & 2.70 & 0.10   \\
        0.025 & 4.43e-04          & 1.84 & 1.10    & 1.67e-04          & 1.81 & 0.87   \\
        0.0125 & 1.11e-04          & 1.99 & 10.31    & 4.24e-05          & 1.98 & 8.36    \\
		0.00625   & 2.81e-05  & 1.98	    & 100.00			&  1.09e-05	 & 1.96		    &88.09  	\\
		\noalign{\smallskip}\hline
	\end{tabular}
	    ~\\[0.2cm]
	    \begin{tabular}{ r | ccr } 
		\multicolumn{1}{c }{}&  \multicolumn{3}{c }{(3) \em nocaps}\\ 	
		\noalign{\smallskip}\hline\noalign{\smallskip}
		$h$ 	 	&  $\|u-u_h\|_{L^2}$   	&rate 	&time [\%] 
		\\	
		\noalign{\smallskip}\hline\noalign{\smallskip}
		0.1	   & 2.80e-02	&- 	    &0.01	\\		
		0.05	 & 3.92e-03	&2.84	&0.09	\\
		0.025   & 1.04e-03	&1.91 	&0.81	\\
		0.0125  & 2.57e-04	&2.02 	&7.76	\\
		0.00625  & 6.45e-05	&2.00	&86.25	\\			\noalign{\smallskip}\hline
	\end{tabular}
	    ~\\[0.2cm]
		\begin{tabular}{ r   | c c r| c cr  } 
		\multicolumn{1}{c }{}& \multicolumn{3}{c }{(4) \em barycenter} &\multicolumn{3}{c }{(5) \em overlap} \\	
		\noalign{\smallskip}\hline\noalign{\smallskip}
		$h$ 	 	& $\|u-u_h\|_{L^2}$   	&rate 		& time [\%]	& $\|u-u_h\|_{L^2}$   	&rate 	&time [\%]	\\
				\noalign{\smallskip}\hline\noalign{\smallskip}
		0.1   & 1.71e-01		&-&   0.00	        & 1.54e-01	&- 	    &0.01	\\ 
		0.05  & 6.00e-02		&1.51&0.03      & 9.88e-02	&0.65	&0.05	\\ 
		0.025  & 1.51e-02		&1.99&0.37      & 6.49e-02	&0.60 	&0.50	\\
		0.0125  & 2.34e-03		&2.69&4.24    	& 3.71e-02	&0.81 	&5.23	\\
	    0.00625  & 4.64e-04		&2.33 &54.49  	& 1.95e-02	&0.92	&63.27	\\
		\noalign{\smallskip}\hline
	\end{tabular}
	\caption{
		Errors and relative assembly costs for the exact ball and for different ball approximations for the inner integrals corresponding to uniform grids of type (a) in Fig. \ref{fig:nugrid}. The relative assembly costs are obtained by dividing the absolute assembly time of the respective run by the largest assembly time for Experiment (1) with the finest grid. Also note that $h$ corresponds to the uniform grid sizing in each dimension, so that the diameter of each element is given by $\sqrt{2}h$.
		}
	\label{num:approx_inner}
\end{table}

For the most part, our predictions concerning the rates corresponding to the different balls as well as the heuristics concerning the choice of quadrature rules are confirmed by the numerical results. However, there are two apparent anomalies. 

The first is that the convergence rate for the $barycenter$ ball approximation is better than what we are able to prove, in fact it is of second-order. The better than linear convergence rate obtained using the $barycenter$ ball gives credence to the possible explanation for this behavior given in Sec. \ref{bary}.

The second anomaly is that although both converge at the expected second-order rate, the errors for the \textit{exactcaps} ball are larger than that for the $approxcaps$ ball. This behavior is due to our use of a one-point quadrature for caps for the former whereas we use a three-point Gauss quadrature formula for the latter. As a consequence, the constant in the $\mcO(h^2)$ relation is smaller for the $approxcaps$ case compared to that for the \textit{exactcaps} case. This comparison shows how using quadrature rules that are more accurate than needed to achieve optimal convergence rates can result in smaller constants in the $\mcO(h^2)$ relation.
\begin{figure}[ht]
\centering
\includegraphics[height=0.4\textwidth]{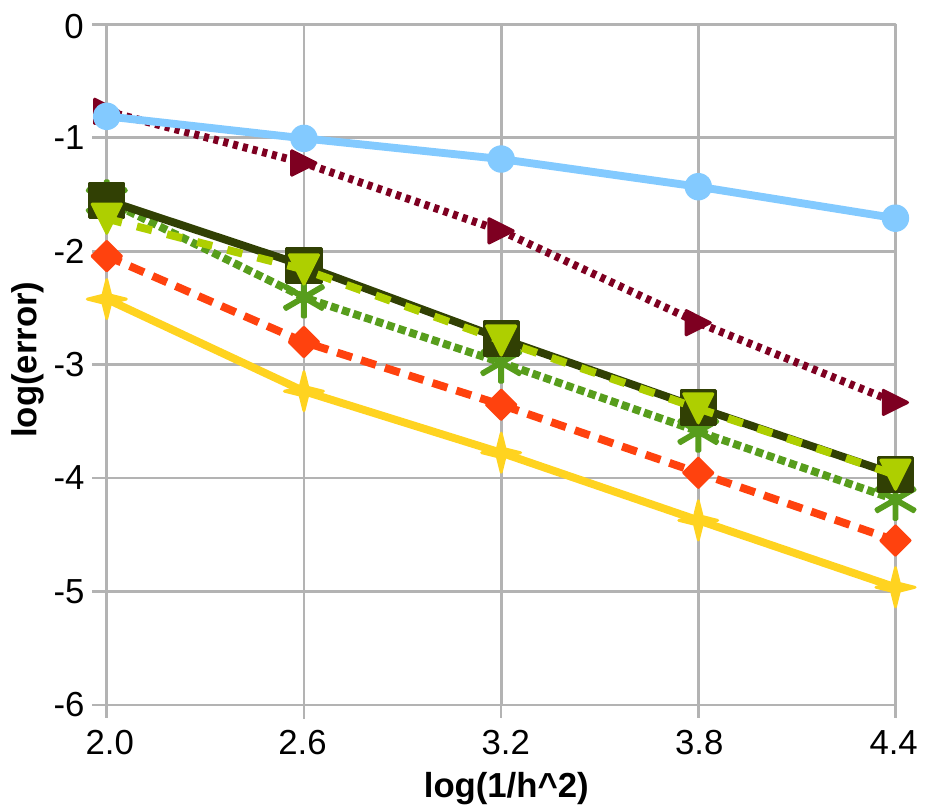}\quad
\includegraphics[height=0.4\textwidth]{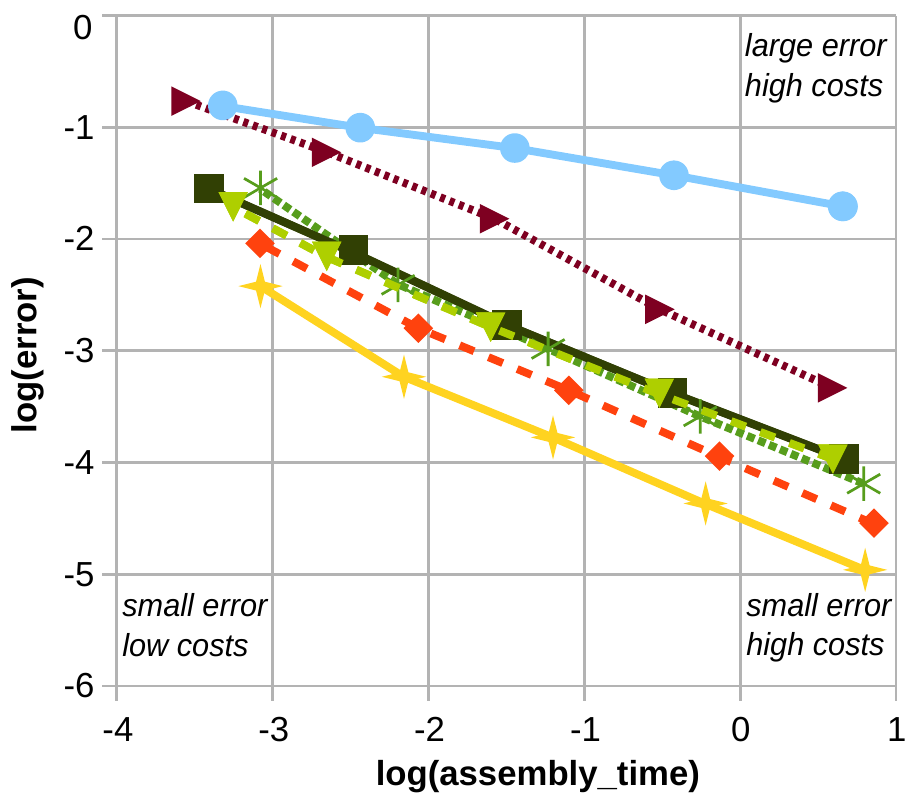}\vskip5pt \includegraphics[height=0.2\textwidth]{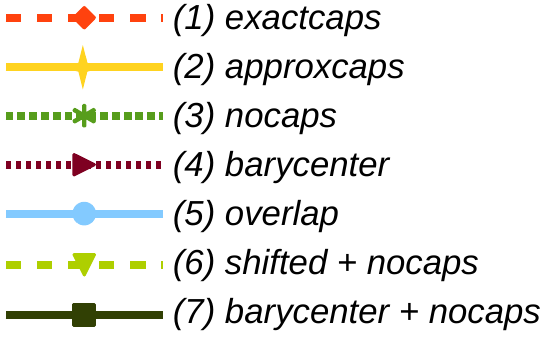} 
\vspace{-0.3cm}
\caption{Plots of errors vs. grid sizes (left) and assembly times (right) of the results that are given in Tables \ref{num:approx_inner}, \ref{num:bary_shifted} and \ref{num:bary} with the legend numbers corresponding to the numbering of columns in those tables. }
\label{fig:loglogplots}
\end{figure}
 
\subsubsection{\bf Approximate shifted ball}
 As proposed in Sec. \ref{subsec:shifted}, for all quadrature points of the outer integral triangle one could \textit{shift} the corresponding ball to the barycenter; one could then approximate this ball by choosing any of $\{exactcaps, approxcaps, nocaps, barycenter\}$.
 In Table \ref{num:bary_shifted} column (6), related numerical results are presented by using the $nocaps$ variant to approximate the shifted ball; see also Fig. \ref{fig:loglogplots}. For outer integral triangles we use a four-point, precision-three symmetric Gaussian rule and for inner integral triangles and potential subelements a three-point, precision-two symmetric Gaussian rule. The results are comparable to that in column (7) of Table \ref{num:bary}. However, the errors are higher compared to column (3) in Table \ref{num:bary_shifted} for the same ball but with no shift. However, we cannot yet explain the observed second-order convergence rate for shifted ball approximations, although the conjecture about this anomaly  given in Sec. \ref{subsec:shifted} is supported by the results of Table \ref{num:bary_shifted}.
\begin{table}[h]
\centering	
\footnotesize
\begin{tabular}{ r  | c cr } 
\multicolumn{1}{c }{}& \multicolumn{3}{c }{(6) {\em shifted + nocaps}}\\
\noalign{\smallskip}\hline\noalign{\smallskip}
$h$	& $\|u-u_h\|_{L^2}$   	&rate 	& time [\%]\\
\noalign{\smallskip}\hline\noalign{\smallskip}
0.1	 		& 1.93e-02	 & - 	    & 0.01 \\  
0.05 				& 6.89e-03	 & 1.48 	& 0.03 \\
0.025 			& 1.62e-03 	 & 2.08		& 0.35 \\	
0.0125 			& 4.11e-04	 & 1.98	 	& 4.15 \\
0.00625 	& 1.06e-04	 & 1.94		& 54.93 \\
\noalign{\smallskip}\hline
\end{tabular}
\caption{Errors and relative assembly costs for the $shifted$ ball approximated by the $nocaps$ variant. The results correspond to the use of uniform grids of type (a) in Fig. \ref{fig:nugrid}.}
\label{num:bary_shifted}
\end{table}

\subsubsection{\bf Improving outer integral approximation for barycenter approximate balls}\label{sec:barynumbers}

As observed in Sec. \ref{subsubsec:discontinuity_outer}, when using the $barycenter$ ball approximation, the integrands for some outer integral triangles have a jump discontinuity for certain pairs of outer and inner integral triangles ({Case 2} in Section \ref{subsec:case2}). In fact, the support of such integrands of an outer integral triangle $\mcE_k$ is given by $S^{barycenter}_{k,k'} = \mcE_k \cap B_\del(\xb^{barycenter}_{k'})$ as is illustrated in the first plot in Fig. \ref{fig:outerq}e. The results in column (4) of Table \ref{num:approx_inner} for the $barycenter$ approximate ball were obtained by applying a quadrature rule to the whole outer integral triangle $\mcE_k$, even when the integrand is discontinuous over such triangles. Although we observe the conjectured improved convergence rates over that which is proved for this case, we also observe erratic behavior in those rates.

Here we consider the question of possible improvements in convergence behaviors accruing from using a quadrature rule not over the whole outer integral triangle, but just over the support region $S^{barycenter}_{k,k'}$ which is illustrated in red in the first plot of Fig. \ref{fig:outerq}e. More precisely, we integrate over polygonal approximations to that support region. We consider two such geometric approximations. The first, which we refer to as the $barycenter + nocaps$ case, is illustrated in green in the second plot of Fig. \ref{fig:outerq}e. An improved geometric approximation is obtained by adding approximate caps as is illustrated in the third plot of Fig. \ref{fig:outerq}e.

Numerical results are presented in Table \ref{barycaps} and Fig. \ref{fig:loglogplots}. For outer integral triangles and potential subelements we use a four-point, precision-three symmetric Gaussian rule and for inner integral triangles a three-point, precision-two symmetric Gaussian rule. We compare the results to those of column (4) of Table \ref{num:approx_inner}. Because we improve the quadrature quality for the outer integrals by taking care of the discontinuity, we not only produce smaller errors but also less erratic rates. %
Computation times slightly increase due to the additional intersection task needed for some outer integral triangles. In fact, because similar tasks are required, they are comparable to those of the $shifted+nocaps$ method presented in Table \ref{num:bary_shifted}. Furthermore, we observe a second-order convergence rate although we have only proven a first-order rate (this is also the case for the results of column (4) in Table \ref{num:approx_inner}). As already alluded to, we conjecture that this is due to a cancellation effect.  Also, among all ball approximations, the $nocaps$ variant for the outer integral triangle performed best in terms of errors (slightly better than using the $approxcaps$ variant for the outer integral triangle). We do not have an explanation for this behavior, but conjecture that the cancellation effect may again be in play.

\begin{table}
\centering	
\footnotesize
\begin{tabular}{ r | cc | cc  } 
\multicolumn{1}{c }{}& \multicolumn{2}{c }{(7) {\em barycenter + nocaps}} & \multicolumn{2}{c }{(8) \em barycenter + approxcaps}  \\ 	
\noalign{\smallskip}\hline\noalign{\smallskip}
$h$   &  $\|u-u_h\|_{L^2}$ & rate & $\|u-u_h\|_{L^2}$ & rate  \\
\noalign{\smallskip}\hline\noalign{\smallskip}
0.1      &  2.81e-02   & -    & 3.50e-02 & -    \\
0.05     &  7.47e-03   & 1.91 & 9.68e-03 & 1.85 \\
0.025    &  1.70e-03   & 2.13 & 2.47e-03 & 1.96 \\
0.0125   &  4.18e-04   & 2.02 & 6.26e-04 & 1.98 \\
0.00625  &  1.07e-04   & 1.96 & 1.58e-04 & 1.98 \\ 
\noalign{\smallskip}\hline
\end{tabular}
\caption{Errors for two barycenter-based approximate balls that (approximately) respect discontinuities in the integrand of the outer integral. The results correspond to the use of uniform grids of type (a) in Fig. \ref{fig:nugrid}.}\label{barycaps}
\end{table}

\subsection{Nonuniform grids}\label{sec:nonuniform}

One naturally asks if the better than provable rates given in Tables \ref{num:approx_inner} and \ref{barycaps} for barycenter ball approximations are an artifact due to the use of uniform Cartesian grids with the same grid size in the both directions. In this section we address this question.
\begin{table}
\centering	
\footnotesize
\begin{tabular}{ rcc | cc | cc  } 
\multicolumn{3}{c}{\text{(b)}}& \multicolumn{2}{c }{(9) {\em barycenter + nocaps}} & \multicolumn{2}{c }{(10) \em barycenter + approxcaps}  \\ 	
\noalign{\smallskip}\hline\noalign{\smallskip}
$J_\Omega$ & $h_{min}$ & $h_{max}$  &  $\|u-u_h\|_{L^2}$ & rate & $\|u-u_h\|_{L^2}$ & rate  \\
\noalign{\smallskip}\hline\noalign{\smallskip}
81    & 0.1004 & 0.2147 & 2.65e-02 & -    & 4.22e-02 & -     \\
361   & 0.0500 & 0.1095 & 1.03e-02 & 1.37 & 1.41e-02 & 1.58  \\
1521  & 0.0250 & 0.0553 & 2.56e-03 & 2.01 & 3.85e-03 & 1.88  \\
6241  & 0.0125 & 0.0278 & 6.37e-04 & 2.01 & 9.92e-04 & 1.96  \\
25281 & 0.0062 & 0.0139 & 1.58e-04 & 2.00 & 2.51e-04 & 1.98  \\    	
\noalign{\smallskip}\hline ~\\  
\multicolumn{3}{c}{\text{(c)}}& \multicolumn{2}{c }{(11) {\em barycenter + nocaps}} & \multicolumn{2}{c }{(12) \em barycenter + approxcaps}  \\ 	
\noalign{\smallskip}\hline\noalign{\smallskip}
$J_\Omega$ & $h_{min}$ & $h_{max}$  &  $\|u-u_h\|_{L^2}$ & rate & $\|u-u_h\|_{L^2}$ & rate  \\
\noalign{\smallskip}\hline\noalign{\smallskip}
166   & 0.0259 & 0.1412 & 8.27e-03 & -    & 1.56e-03 & -   \\
712   & 0.0119 & 0.0699 & 1.65e-03 & 2.32 & 4.06e-03 & 1.95 \\
2924  & 0.0060 & 0.0353 & 3.98e-04 & 2.05 & 1.06e-04 & 1.94 \\
11750 & 0.0028 & 0.0184 & 1.19e-04 & 1.74 & 2.84e-04 & 1.90 \\
\noalign{\smallskip}\hline~\\
\multicolumn{3}{c }{\text{(d)}}& \multicolumn{2}{c }{(13) {\em barycenter + nocaps}} & \multicolumn{2}{c }{(14) \em barycenter + approxcaps}  \\ 	
\noalign{\smallskip}\hline\noalign{\smallskip}
$J_\Omega$ & $h_{min}$ & $h_{max}$   &  $\|u-u_h\|_{L^2}$ & rate & $\|u-u_h\|_{L^2}$ & rate  \\
\noalign{\smallskip}\hline\noalign{\smallskip}
418   & 0.0143 & 0.1335 & 5.38e-03 & -    & 1.14e-03 & -    \\
1400  & 0.0070 & 0.1075 & 1.79e-03 & 1.58 & 3.65e-03 & 1.65 \\
5477  & 0.0032 & 0.0816 & 6.35e-04 & 1.50 & 1.07e-03 & 1.77 \\
20755 & 0.0018 & 0.0534 & 2.94e-04 & 1.11 & 3.79e-04 & 1.50 \\    
		\noalign{\smallskip}\hline
	\end{tabular}
\caption{{Results for the $barycenter$ ball approximation for the inner integrals combined with the variants $nocaps$ and $approxcaps$ to approximate the support of the integrand of the outer integral. Each table (b)--(d) corresponds to the grids (b)--(d) depicted in Fig. \ref{fig:nugrid}. Here, $J_\Omega$ denotes the number of FEM nodes (degrees of freedom) inside $\Omega$  and $h_{min}$ ($h_{max}$) denote the minimum (maximum) diameter over all triangles.
}}
\label{num:bary}
\end{table}

The results given in Table \ref{num:bary} indicate that the approximate $\mcO(h^2)$ convergence is achieved for nonuniform grids with smooth transitions in the grid size such as that depicted in Figs. \ref{fig:nugrid}b (for columns (9) and (10) in that table) and \ref{fig:nugrid}c (for columns (11) and (12)). We conjecture that this effect holds for such grids because the beneficial cancellation effect is a ``localized'' phenomenon. By this we mean that if changes in the grid are sufficiently smooth and if the grid size is small enough, then for any sufficiently short arc of the boundary of the ball the cancellation effect occurs because the grid along that arc is quasi-uniform. 

We observe that even for the very highly nonuniform grid illustrated in Fig. \ref{fig:nugrid}d, convergence rates higher than the proven first-order rate are obtainable, at least for the ${barycenter + approxcaps}$ case. Another observation is that the erratic convergence behavior seen in Table \ref{num:bary}(d) also occurs for other approximate balls, even for those for which convergence rates are provably $\mcO(h^2)$. Such erratic behavior is to be expected because clearly $h_{max}$ is nowhere near small enough for computed errors to be in the asymptotic range required for error estimates to hold. Another cause for the erratic convergence behavior when using grids of type (d) (and also, to a somewhat lesser extent, for grids of type (c)) is that grid refinement is effected  using non-nested grids; erratic behaviors are often observed for such refinements, especially for relative coarse grid sizes.

\section{Closing remarks and recommendations}\label{conclusion}

\subsection{Closing remarks}

\textbf{\em Higher-order FEMs.}
An often-stated advantage of finite element methods is the relative ease with which higher-order discretizations can be constructed. Unfortunately, the geometric errors incurred by the approximate balls we have considered (that are best of ${\mathcal O}(h^2)$) would dominate over the approximation capabilities of higher-order polynomial finite element bases. {\em Using exact caps would clearly be useful in this context because no geometric error would be incurred.}

Alternately, one could approximate the cap by many small triangles. This approach would not change the rate of convergence of the geometric error but would render much smaller the constant in the ${\mathcal O}(h^2)$ relation. In practical computations, one often selects a desired value of the grid size, so that making that constant small enough would, for that fixed grid size, make the geometric error commensurate with the other errors incurred. Of course, this approach incurs greater assembly and solution costs relative to that for exact caps.}  

Another approach along these lines is to use a higher-precision quadrature rule, i.e., a quadrature rule for caps that employs many quadrature points; see Ref. \refcite{sectorrule}. Again, rates of convergence would not be improved, but constants in order relations may be significantly smaller so that again, for a fixed grid size and for a sufficient number of quadrature points, geometric errors may be significantly lessened.

\vskip5pt
\noindent\textbf{\em Other approximate balls.}
Other approximations to Euclidean balls come to mind. 
{For example, as an alternative to the inscribed triangle-based polygon of Fig. \ref{fappoxball}a, one could instead use a {\em regular inscribed polygon}; see {Ref. \refcite{Bond2015}}. To preserve accuracy, the sides of the regular polygon would have to be of ${\mathcal O}(h)$. The advantage of doing so is that the definition of a regular polygon is independent of the finite element triangulation, i.e., to construct a regular polygon one does not have to determine intersections of the boundary of the ball with triangle edges. However, there are disadvantages in using regular polygons. For example, for the purpose of defining a composite quadrature rule over a regular polygon, one can easily subdivide the ball into triangles that are not finite element triangles; however, in this case the finite element approximation would be a piecewise polynomial that, because they are not finite element triangles, is merely continuous over those triangles which compromises the accuracy of the quadrature rule over that triangle. On the other hand, the construction of a triangulation of a regular polygon so that all triangles are contained within finite element triangles (so that the integrand is smooth) becomes a substantially more cumbersome task compared to that for polygons such as that depicted in Fig. \ref{fappoxball}a. {Not only does one have to now determine the points of intersection of the boundary of the regular polygon and the sides of the finite elements,} but one also has to deal with the fact that the vertices of the polygon are generally in the interior and not at the edges of the finite elements.}

Another possibility that is a whole-triangle alternative to the $\sharp=\{barycenter\}$ case is to keep whole triangles whenever the overlap with the Euclidian ball is greater or equal to half the area of the triangle; otherwise, a triangle is not included in the approximate ball. This approach requires the same steps as does the $\sharp=\{nocaps\}$ case, i.e, the determination of circle-triangle intersection points and the subsequent subdivision of quadrilaterals into triangles. However, it also requires the additional step of determining the area of the overlap. We do not study this type of approximate ball because it is much more difficult to implement compared to the $\sharp=\{barycenter\}$ case and more difficult to implement than even the $\sharp=\{nocaps,approxcaps\}$ cases and, also, it does not yield provably better rates of convergence than the latter two cases. 

Other examples are provided by, e.g., isoparametric, isogeometric, and extension approximations of the circle in much the same way as are used for finite element approximations for PDE problems posed on domains with curved boundaries; see, e.g., Refs. \refcite{brenner,cheung,ciarlet,isogeo}. For example, in such a method, the curved boundary of an element is often approximated by a polynomial. We do not study this type of approximate ball because, for a circle, it is much more efficient to use circular caps and, as a bonus, no geometric error is incurred.

\vskip5pt
\noindent\textbf{\em Towards three-dimensional finite element approximations.}
All ball approximations and the attendant quadrature rules used in our two-dimensional studies can be extended to the three-dimensional setting. However, some of the construction steps used such as determining  intersections of spheres and tetrahedrons are substantially more complicated to implement in three dimensions. Furthermore, the error vs. cost criterion that is used to select the ``best'' recipe could result in a different outcome in three dimensions. 

\subsection{Recommendations}

As we have repeatedly seen in the paper, the implementation of finite element methods for nonlocal models with a finite range of interaction is particularly challenging when the diameter of the interaction set is smaller than that of the domain. In fact, one has to compute integrals over the intersection between the interaction set (typically a Euclidean ball) and the elements of the mesh. For the two-dimensional case, we investigated several approaches to approximate this intersection through the use of the ball approximations $\sharp=\{nocaps,~approxcaps,~barycenter,~overlap,~shifted\}$ and the mixtures $\{barycenter+nocaps,~ barycenter+approxcaps,~shifted+nocaps\}$, all of which incur a geometric error. We also compared the use of these approximate balls to an approach that, through the use of quadrature rules for circular caps, uses the \textit{exactcaps} ball so that no geometric error is incurred.

All in all, comparing the error-to-cost ratio in Table \ref{num:approx_inner} and Fig. \ref{fig:loglogplots}-right, we conclude that, at least in the two-dimensional setting, the ${approxcaps}$ approximation is preferred over all other methods investigated in this paper, with the caveat that if the quadrature rules having the same precision are used for both the ${exactcaps}$ and ${approxcaps}$ cases, the error for the former would be lower and may in fact render the ${exactcaps}$ approximation to be superior. Also, as noted above, if higher-order finite element methods are used, the ${exactcaps}$ approach has the singular advantage over all the other methods because it does not incur any geometric error.

It should be noted that an approach that ``wins'' in two dimensions may or may not ``win'' in three dimensions. For this reason, other methods such as the ${shifted + nocaps}$ and ${barycenter + nocaps}$ approaches may warrant study in three dimensional settings. 

\vskip5pt
\noindent In follow-up work, we will delve deeply into the three-dimensional setting and also into higher-order finite element methods.

\section*{Acknowledgements}

The research of MD is supported by Sandia National Laboratories (SNL). SNL is a multimission laboratory managed and operated by National Technology and Engineering Solutions of Sandia, LLC., a wholly owned subsidiary of Honeywell International, Inc., for the U.S. Department of Energy's National Nuclear Security Administration under contract DE-NA-0003525. Note that this paper describes objective technical results and analysis. Any subjective views or opinions that might be expressed in the paper do not necessarily represent the views of the U.S. Department of Energy or the United States Government. Report number SAND2020-4407.

The research of MG is partially supported by the US Department of Energy Office of Science through a subcontract issued by the Oak Ridge National Laboratory.

The research of CV has been supported by the German Research Foundation (DFG) within the Research Training Group 2126: “Algorithmic Optimization”. 

CV is grateful to Professor Volker Schulz for all his valuable advice.

\appendix
{
\section{Proof of Proposition \ref{prop:error1}}\label{appendixa}
 
From \eqref{lfunc0}, \eqref{lfunc2}, \eqref{weakha}, and \eqref{D-Dh-weak}, we have that
$$
\left\{
\bal
&A(u_h,v_h) = F(v_h) = G(v_h) + G_{g}(v_h)
\\
&A_h(\wuh,v_h) = F_h(v_h) = G(v_h) + G_{g,h}(v_h)
\eal
\right.
\qquad\mbox{for $v_h \in V_c^h$}
$$
so that
$$
     A(u_h,v_h) = A_h(\wuh, v_h) - G_{g,h}(v_h) + G_{g}(v_h) ,
$$
where
$$
 G_{g}(v_h) = 2\int_{\Omega}v_h(\xb) 
\bigg(\int_{\Omega_{\mathcal I}\cap B_{\delta}(\xb)} g(\yb)\psi(\xb,\yb)d\yb  
  \bigg)d\xb
$$
and
$$
  G_{g,h}(v_h) = 2\int_{\Omega}v_h(\xb) 
\bigg(\int_{\Omega_{\mathcal I}\cap B_{\delta,h}(\xb)} g(\yb)\psi(\xb,\yb)d\yb  
  \bigg)d\xb.
$$
Then,
\beq{uwuh}
\bal
   | A( u_h - \wuh, v_h)| &= | A( u_h , v_h) - A(\wuh, v_h)|
  \\& = | A_h(\wuh, v_h) -  G_{g,h}(v_h)+  G_{g}(v_h) - A(\wuh, v_h)  |
  \\&
   \le | A_h(\wuh, v_h) - A(\wuh, v_h)| + |  G_{g,h}(v_h) - G_{g}(v_h) |.
\eal
\eeq

From \eqref{bform2} and \eqref{bform3hh}, we have that 
$$ 
\bal
  &   | A_h(w, z) - A(w, z)| 
\\&\quad \leq     \int_\Omega\int_\Omega \big|w(\xb)-w(\yb)\big|\big|z(\xb)-z(\yb)\big|\psi(\xb,\yb)
       \big| \mcX_{B_\delta(\xb)}(\yb)-\mcX_{B_{\delta,h}(\xb)}(\yb)\big|\,d\yb\,d\xb
\\
&\qquad\qquad
+ 2\int_{\Omega}\big|w(\xb)z(\xb)\big|\bigg(\int_{\Omega_{\mathcal I}}\psi(\xb,\yb)\big|\mcX_{B_\delta(\xb)}(\yb)-\mcX_{B_{\delta,h}(\xb)}(\yb)\big| \,d\yb\bigg) \,d\xb \, 
 \\[2mm]
& \quad =   \int_\Omega\int_{\Omega\cap\Delta B_{\delta,h}(\xb)} \big|w(\xb)-w(\yb)\big|\big|z(\xb)-z(\yb)\big| \psi(\xb,\yb)\,d\yb\,d\xb
\\
&\qquad\qquad
+ 2\int_{\Omega}\big|w(\xb)z(\xb)\big|\bigg(\int_{\Omega_{\mathcal I}\cap\Delta B_{\delta,h}(\xb)}\psi(\xb,\yb)\,d\yb\bigg) \,d\xb \, 
\eal
$$
so that
\beq{eq:d-estimate1}
\bal
 | A_h&(w, z) - A(w, z)|
\\& \leq 
\underbrace{\left(\int_\Omega\int_{\Omega\cap\Delta B_{\delta,h}(\xb)}
       \big(w(\xb)-w(\yb)\big)^2\psi(\xb,\yb)\,d\yb\,d\xb\right)^\frac12}_{I}
\\
&\qquad\qquad\qquad
\times \underbrace{\left(\int_\Omega\int_{\Omega\cap\Delta B_{\delta,h}(\xb)} 
       \big(z(\xb)-z(\yb)\big)^2\psi(\xb,\yb)\,d\yb\,d\xb\right)^\frac12}_{II}
\\
&\qquad 
+2\left( 
\int_{\Omega}w^2(\xb)\bigg(\int_{\Omega_{\mathcal I}\cap\Delta B_{\delta,h}(\xb)}\psi(\xb,\yb)\,d\yb\bigg) \,d\xb  \, \right)^\frac12
\\
&\qquad
\underbrace{\phantom{\qquad\qquad}\times \left( 
\int_{\Omega}z^2(\xb)\bigg(\int_{\Omega_{\mathcal I}\cap\Delta B_{\delta,h}(\xb)}\psi(\xb,\yb)\,d\yb\bigg) \,d\xb  \, \right)^\frac12}_{III},
\eal
\eeq
where $\Delta B_{\delta,h}=(B_\delta\setminus(B_\delta\cap B_{\delta,h}))\cup(B_{\delta,h}\setminus(B_\delta\cap B_{\delta,h}))$ and where we have used the Cauchy-Schwarz inequality. Also, {\em III}\, refers to the last two lines of \eqref{eq:d-estimate1}.

For the ({\em I}\;) term, we have that
$$
\bal
I^2 &= \int_\Omega\int_{\Omega\cap\Delta B_{\delta,h}(\xb)} \big(w(\xb)-w(\yb)\big)^2\psi(\xb,\yb)\,d\yb\,d\xb 
\\&\leq
2\int_\Omega\int_{\Omega\cap\Delta B_{\delta,h}(\xb)}\big(w^2(\xb)+w^2(\yb)\big)\psi(\xb,\yb)\,d\yb\,d\xb.
\eal
$$
For the  $w^2(\xb)$ term, we obtain
$$
\bal
     &\int_\Omega\int_{\Omega\cap\Delta B_{\delta,h}(\xb)}w^2(\xb)\psi(\xb,\yb)\,d\yb\,d\xb 
\\&\qquad
\leq \|w\|^2_{L^2(\Omega)} \sup_{\xb\in\Omega}\,\Big(
     \int_{\Omega\cap\Delta B_{\delta,h}(\xb)} \psi(\xb,\yb)\,d\yb\Big)
\leq   \QQQ_\Omega \,\sup_{\xb\in\Omega}|\Delta B_{\delta,h}(\xb)|\, \|w\|^2_{L^2(\Omega)},
\eal
$$
where 
$
\QQQ_\Omega=\sup_{\xb\in\Omega}\,\sup_{\yb\in\Omega} \psi(\xb,\yb).
$

Following the same arguments for the remaining term in ({\em I}\,) and the two analogous terms in ({\em II}\,),  we have
\beq{ineq-equiv}
\bal 
\mbox{({\em I}\,)({\em II}\,)} & \leq \left(4\QQQ_\Omega \,\sup_{\xb\in\Omega}|\Delta B_{\delta,h}(\xb)|\,\|w\|^2_{L^2(\Omega)}\right)^\frac12 
               \left(4\QQQ_\Omega \,\sup_{\xb\in\Omega}|\Delta B_{\delta,h}(\xb)|\,\|z\|^2_{L^2(\Omega)}\right)^\frac12
       \\&   \leq 4\QQQ_\Omega\, \sup_{\xb\in\Omega}|\Delta B_{\delta,h}(\xb)|\, \|w\|_{L^2(\Omega)} \|z\|_{L^2(\Omega)}.
\eal
\eeq
Also proceeding in a similar manner for the ({\em III}\,) term in \eqref{eq:d-estimate1}, we have that
\beq{ineq-equiv2}
\bal
III & \leq 2\left(\QQQ_{\Omega_{\mathcal I}} \,|\Delta B_{\delta,h}(\xb)|\,\|w\|^2_{L^2(\Omega)}\right)^\frac12 
               \left(\QQQ_{\Omega_{\mathcal I}} \,\sup_{\xb\in\Omega}|\Delta B_{\delta,h}(\xb)|\,\|z\|^2_{L^2(\Omega)}\right)^\frac12
        \\& \leq 2\QQQ_{\Omega_{\mathcal I}}\,  \sup_{\xb\in\Omega}|\Delta B_{\delta,h}(\xb)|
\, \|w\|_{L^2(\Omega)} \|z\|_{L^2(\Omega)},
\eal
\eeq
where
$
\QQQ_{\Omega_{\mathcal I}}=\sup_{\xb\in\Omega}\,\sup_{\yb\in\Omega_{\mathcal I}} \psi(\xb,\yb).
$ 
Substituting \eqref{ineq-equiv} and \eqref{ineq-equiv2} into \eqref{eq:d-estimate1} results in 
\beq{ttttt1} 
| A_h(w, z) - A(w, z)| \le (4\QQQ_\Omega + 2\QQQ_{\Omega_{\mathcal I}})\, \sup_{\xb\in{\Omega}}|
\Delta B_{\delta,h}(\xb)|\, \|w\|_{L^2(\Omega)} \|z\|_{L^2(\Omega)}.
\eeq
Next, we have that
\beq{tttttt} 
\bal
&|  G_{g,h}(z) - G_{g}(z) | \\  
&\qquad \leq  2\int_{\Omega}|z(\xb)| 
 \bigg(\int_{\Omega_{\mathcal I}} |g(\yb)|\psi(\xb,\yb)
 \big|\mcX_{B_\delta(\xb)}(\yb)-\mcX_{B_{\delta,h}(\xb)}
 (\yb)\big|d\yb\bigg)d\xb  \\
&\qquad= 2\int_{\Omega}
   \int_{\Omega_{\mathcal I}\ }
  |z(\xb)| |g(\yb)|\psi(\xb,\yb)\mcX_{\Delta B_{\delta,h}(\xb)}(\yb) d\yb d\xb   \\ 
&\qquad\leq 2  \|z\|_{L^2(\Omega)} \sqrt{\QQQ_{\Omega_{\mathcal I}}  \sup_{\xb\in{\Omega}}|
\Delta B_{\delta,h}(\xb)|} \cdot\|g\|_{L^2(\Omega_{\mathcal I})} \sqrt{\QQQ_{\Omega_{\mathcal I}} \sup_{\xb\in{\Omega}}|\Delta B_{\delta,h}(\xb)|}\\
&\qquad = 2  \|z\|_{L^2(\Omega)} \|g\|_{L^2(\Omega_{\mathcal I})} \QQQ_{\Omega_{\mathcal I}} \sup_{\xb\in{\Omega}}|\Delta B_{\delta,h}(\xb)|,
\eal
\eeq
where we have used the Cauchy-Schwarz inequality on $L^2(\Omega \times \Omega_{\mathcal I})$
for the second inequality

Setting $w=\wuh$ and $z=v_h=u_h - \wuh$ and substituting \eqref{ttttt1} and \eqref{tttttt} into \eqref{uwuh} results in
\beq{uwuh111}
\bal
   &| A( u_h - \wuh, u_h - \wuh)| 
   \\
   &\le \Big( (4\QQQ_\Omega + 2\QQQ_{\Omega_{\mathcal I}})\|\wuh\|_{L^2(\Omega)}+2\|g\|_{L^2(\Omega_{\mathcal I})}\QQQ_{\Omega_{\mathcal I}})\Big)\, \sup_{\xb\in\Omega}(|\Delta B_{\delta,h}(\xb)|)\;\|u_h - \wuh\|_{L^2(\Omega)}.
\eal
\eeq
Because the well posedness of the problem \eqref{D-Dh-weak} implies that $\|\wuh\|_{L^2(\Omega)}$ can be bounded by norms of the data $f$ and $g$, \eqref{eq:energy-diff1} follows from \eqref{uwuh111} and the definition of the energy norm.\quad$\Box$}

\bibliographystyle{amsplain}
\bibliography{literature}

\end{document}